\documentclass[12pt,a4paper]{report}
\usepackage[utf8]{inputenc}
\usepackage[T1]{fontenc}
\usepackage{amsmath}
\usepackage{amsfonts}
\usepackage{fancyhdr}
\usepackage{amssymb}
\usepackage{color} 
\usepackage{mdframed}
\usepackage{multirow} 
\usepackage{multicol} 
\usepackage{scrextend} 
\usepackage{tikz}
\usepackage{graphicx}
\usepackage[absolute]{textpos} 
\usepackage{colortbl}
\usepackage{array}
\usepackage{geometry}

\usepackage{amsfonts}
\usepackage{mathtools}
\usepackage{amssymb}
\usepackage[noadjust]{cite}
\usepackage{hyperref}
\usepackage{mathrsfs}
\usepackage{tikz-cd}
\usepackage{amsthm}

\usepackage{titlesec}
\usepackage{lipsum}

\titleformat{\chapter}[display]
  {\normalfont\bfseries}{}{0pt}{\Large}

\hypersetup{
	colorlinks = true,
	linkcolor = black,
	citecolor=black
}

\tikzset{
  symbol/.style={
    draw=none,
    every to/.append style={
      edge node={node [sloped, allow upside down, auto=false]{$#1$}}}
  }
}

\newtheorem{theorem}{Theorem}[section]
\newtheorem*{theorem*}{Theorem}
\newtheorem*{theoremfr*}{Théorème}

\newtheorem{innercustomthm}{Theorem}
\newenvironment{customthm}[1]
  {\renewcommand\theinnercustomthm{#1}\innercustomthm}
  {\endinnercustomthm}

\newtheorem{innercustomthmfr}{Théorème}
\newenvironment{customthmfr}[1]
  {\renewcommand\theinnercustomthmfr{#1}\innercustomthmfr}
  {\endinnercustomthmfr}

\newtheorem{corollary}[theorem]{Corollary}

\newtheorem{lemma}[theorem]{Lemma}
\newtheorem{proposition}[theorem]{Proposition}
\theoremstyle{definition}
\newtheorem{definition}[theorem]{Definition}
\theoremstyle{remark}
\newtheorem{remark}[theorem]{Remark}
\newtheorem{notation}[theorem]{Notation}

\renewenvironment{proof}[1][Proof.]{\begin{trivlist}\item[\hskip \labelsep {\bfseries #1}]}{\end{trivlist}}
\newenvironment{example}[1][Example.]{\begin{trivlist}
  \item[\hskip \labelsep {\bfseries #1}]}{\end{trivlist}}
\newenvironment{exemple}[1][Exemple.]{\begin{trivlist}
  \item[\hskip \labelsep {\bfseries #1}]}{\end{trivlist}}

\DeclareMathOperator{\Ima}{\mathrm{Im}}
\DeclareMathOperator{\Pic}{\mathrm{Pic}}
\DeclareMathOperator{\rk}{\mathrm{rk}}
\DeclareMathOperator{\Hom}{\mathrm{Hom}}
\DeclareMathOperator{\Aut}{\mathrm{Aut}}
\DeclareMathOperator{\Spec}{\mathrm{Spec}}

\renewcommand{\qed}{\hfill\blacksquare}

\newcommand\restr[2]{{
  \left.\kern-\nulldelimiterspace 
  #1 
  \vphantom{\big|} 
  \right|_{#2} 
  }}

\definecolor{Prune}{RGB}{99,0,60}

\begin{document}

\begin{titlepage}

\newgeometry{left=7.5cm,bottom=2cm, top=1cm, right=1cm}
\tikz[remember picture,overlay] \node[opacity=1,inner sep=0pt] at (-28mm,-135mm){\includegraphics{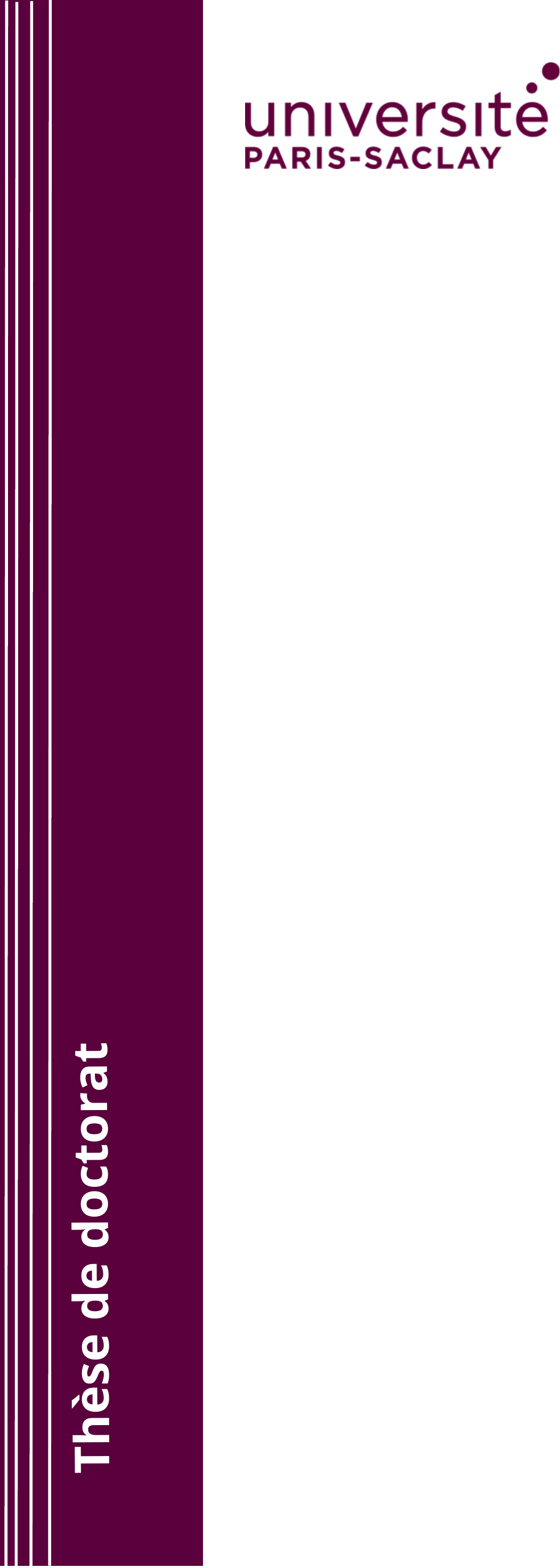}};

\fontfamily{fvs}\fontseries{m}\selectfont


\color{white}
\begin{picture}(0,0)
\put(-150,-735){\rotatebox{90}{NNT: 2020UPASA000}}
\end{picture}
 

\flushright
\vspace{10mm} 
\color{Prune}
\fontfamily{fvs}\fontseries{m}\fontsize{22}{26}\selectfont
  Courbes rationnelles sur les compactifications magnifiques des espaces symétriques.  


\normalsize
\vspace{1.5cm}

\color{black}
\textbf{Thèse de doctorat de l'Université Paris-Saclay}

\vspace{15mm}

\hspace*{-0.7cm}Ecole Doctorale de Mathématique Hadamard (EDMH) n$^{\circ}$ 574\\
\small Spécialité de doctorat~:
Mathématique fondamentale
\\
\footnotesize Unité de recherche~:

Laboratoire de math\'ematiques de Versailles \\(UVSQ), UMR 8100 CNRS\\

\footnotesize Référent~:
  Faculté des sciences d'Orsay\\
\vspace{15mm}

\textbf{Thèse présentée et soutenue à Paris, le 24 juin 2021, par}\\
\bigskip
\Large {\color{Prune} \textbf{Arsen Shebzukhov}}

\vfill
\flushleft \small
\textbf{Au vu des rapports de :}
\bigskip

{\scriptsize
\begin{tabular}{p{8cm}l}
\textbf{Michel Brion} & Rapporteur \\ 
Universit\'e Grenoble Alpes  &   \\ 
\textbf{Baohua Fu } & Rapporteur \\ 
Academie Chinoise des Sciences, Pekin  &   \\ 
\end{tabular}}

\bigskip
\textbf{Composition du jury :}
\bigskip

\scriptsize
\begin{tabular}{|p{8cm}l}
\arrayrulecolor{Prune}
\textbf{Michel Brion} & Rapporteur \\ 
Universit\'e Grenoble Alpes  &   \\ 
\textbf{Baohua Fu} & Rapporteur \\ 
Academie Chinoise des Sciences, Pekin  &   \\ 
\textbf{Dimitry Zvonkine} &   Président\\ 
Laboratoire de Mathématiques de Versailles & \\
\textbf{Nicolas Perrin} &   Directeur\\ 
Laboratoire de Mathématiques de Versailles & \\
\end{tabular} 
\end{titlepage}


\hbox{\includegraphics[width=6cm]{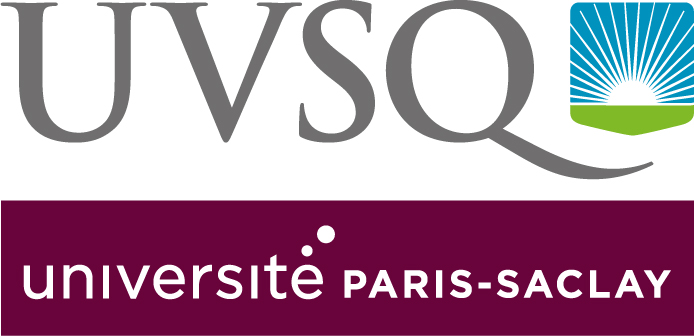}}\vspace*{0.5cm}
\hbox{\includegraphics[width=6cm]{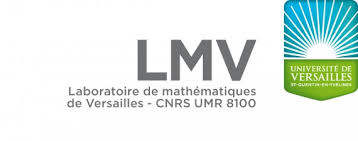}}\vspace*{0.5cm}

\tableofcontents

\chapter{Résumé en français}

\section{Motivation}

Les courbes rationnelles ont toujours joué un rôle important en géométrie algébrique complexe, 
que ce soit pour le programme du modèle minimal (MMP), 
pour les variétés des tangentes rationnelles minimales (VMRT), pour la nouvelle notion 
de simple connexité rationnelle ou pour la géométrie énumérative, comme par exemple la cohomologie quantique.
Dans plusieurs cas, le premier pas est de comprendre l'espace de modules de courbes rationnelles 
sur une variété donnée. Beaucoup de résultats ont été déjà obtenus sur ce sujet, 
par exemple \cite{fulton_pandharipande_nosmaqc} explique la construction 
d'espaces de modules d'applications stables et le fait que ces espaces sont projectifs
dans le cas des variétés projectives.

L'objectif principal de cette thèse est l'étude des courbes rationnelles sur les $G$-variétés magnifiques. 
Ce choix est basé sur le fait
que ces variétés ont une structure d'orbites simple. Cela nous permet de bien contrôler 
le groupe des classes de diviseurs $\Pic(X)$ et son dual, le groupe des courbes $N_1(X)$ modulo
équivalence numérique.

\section{Variétés magnifiques}

Nous allons travailler sur  $\mathbb{C}$.
Soit $G$ un groupe algébrique réductif connexe de rang $\rk(G) = r$. 
Fixons un sous-groupe de Borel $B \subset G$ et un tore maximal $T$ de $G$ avec $T \subset B$. 
Notons $R = R_G$ le système de racines de $(G,T)$, $R_G^+$ les racines de $B$, $R_G^- = -R_G^+$
et $S = S_G$ la base des racines simples.
Un sous-groupe parabolique $P$ de $G$ est dit standard, si $B \subset P$,
notons $P^-$ son sous-groupe parabolique opposé, tel que $T \subset P^-$.
Soit $S' \subset S$, $R_{G,S'}$ est l'ensemble des racines engendrées par $S'$ 
et $P_{S'}$ est le sous-groupe parabolique avec l'ensemble de racines $R^+_G \cup R_{G,S'}$.
Ce parabolique $P_{S'}$ est dit associé à $S'$.

Une $G$-variété normale $X$ avec une $B$-orbite ouverte est dite \textit{sphérique}, si de plus $X$ ne contient 
qu'une seule $G$-orbite fermée, $X$ est dite \textit{simple}. 
Une telle variété possède une $G$-orbite ouverte, isomorphe à $G/H$ pour un sous-groupe $H \subset G$,
qui est dit aussi \textit{sphérique}.
Alors, il y a un nombre fini de $B$ et $G$-orbites dans $X$. 
Nous pouvons définir l'ensemble des diviseurs premiers stables par $B$, mais pas par $G$, appelés \textit{couleurs}
et notés $D \in \Delta_X$.
Chaque $G$-variété sphérique simple projective est définie uniquement à isomorphisme près par le cône des valuations 
$G$-invariantes et le sous-ensemble des couleurs qui contiennent la $G$-orbite fermée de $X$.
Cette donnée est appelée le \textit{cône colorié} de $X$.
Un espace homogène sphérique est dit \textit{sobre}, si son cône des valuations $G$-invariantes
est strictement convexe ou, de manière équivalente, ne contient pas de sous-espace linéaire.
Une variété sphérique simple $X$ est dite \textit{toroïdale}, si la $G$-orbite fermée $Y$
n'est contenue dans aucune couleur $D \in \Delta_X$.
Une combinaison de cônes colorés (vérifiant certaines conditions), appelé l'éventail colorié
peut être utilisé pour classifier les variétés sphériques non-simples dans l'esprit de 
la classification des variétés toriques.

Une variété \textit{magnifique} est une $G$-variété lisse projective ayant une $G$-orbite ouverte $G/H$, telle que 
le complémentaire de cette orbite noté $\partial X = X \setminus (G/H)$ est un diviseur à croisements normaux lisses.
De plus les adhérences des $G$-orbites de $X$ sont exactement les intersections partielles 
de composantes irréductibles de $\partial X$. En particulier, $X$ a une unique $G$-orbite fermée:
l'intersection de toutes les composantes irréductibles.
Le nombre de composantes irréductibles de $\partial X$ est le \textit{rang} de $X$, noté $\rk(X)$.
On peut écrire $\partial X = \cup_{i = 1}^{\rk(X)} X_i$, où tout $X_i$
est un diviseur premier $G$-stable. Un tel diviseur premier $G$-stable est dit 
\textit{diviseur de bord}. Nous allons écrire $X_I = \cap_{i \in I} X_i$ pour 
toute sous-variété $G$-stable de $X$, où $I \subset \{1,\dots, \rk(X)\}$.

Domingo Luna a prouvé dans \cite{luna_tvmes}, que toute variété magnifique est sphéri\-que toroïdale.
Nous décrivons la combinatoire, utilisée pour la classification des variétés magnifiques.
On dit qu'une couleur $D \in \Delta_X$ est \textit{déplacée} par une racine simple $\alpha$, 
si $D$ n'est pas stable sous l'action du sous-groupe parabolique minimal $P_\alpha$, associé à $\alpha$.
L'ensemble des couleurs déplacées par $\alpha$ est noté $\Delta_X(\alpha)$.
On note $S_X^p$ l'ensemble des racines simples $\alpha \in S_G$, tel que $\Delta_X(\alpha) = \emptyset$.
Les \textit{racines sphériques} $\Sigma_X$ sont les $T$-poids de l'espace normal
\[
  \frac{T_z X}{T_z G z},
\]
où $z$ est l'unique point fixe par $B^-$ de la $G$-orbite fermée de $X$.
Soit $A_{G,X} \subset \Delta_X$ l'ensemble des couleurs $D \in \Delta_X$ déplacées par une racine simple 
$\alpha \in S_G \cap \Sigma_X$, qui est aussi une racine sphérique de $X$.
Toute racine sphérique $\gamma \in \Sigma_X$ est une combinaison linéaire à
coefficients positifs de racines simples de $G$. On définit 
le support de $\gamma$, noté $\mathrm{supp}(\gamma)$, comme les racines simples avec des coefficients
non-nuls dans la présentation de $\gamma$. 
Le \textit{système sphérique} de $X$ est défini comme la donnée $(S_X^p, \Sigma_X, A_{G,X})$.

Domingo Luna a conjecturé dans \cite{luna_vsdta}, que toute variété magnifique
est déterminée par son système sphérique à isomorphisme près et il a établi cette conjecture en type $A$. Récemment 
Paolo Bravi et Guido Pezzini l'ont prouvée dans \cite{bravi_pezzini_pwv}:

\begin{theoremfr*}[{\cite{bravi_pezzini_pwv}}]
  Il existe une bijection entre les classes d'isomorphisme de $G$-variétés magnifiques 
  et les $G$-systèmes sphériques.
\end{theoremfr*}

Ils ont introduit la notion de $G$-système abstrait et ont montré que 
tout $G$-système abstrait est le système sphérique d’une $G$-variété magnifique.

\begin{exemple}[Exemple 1]
  Les variétés projectives homogènes sont magnifiques de rang nul. 
  Soit $P_{S'}$ un sous-groupe parabolique standard de $G$ pour un $S' \subset S$. Alors la variété
  projective homogène $G/P_{S'}^-$ possède une $B$-orbite ouverte, où 
  $P_{S'}^-$ est le sous-groupe parabolique, opposé de $P_{S'}$.
  La $B$-orbite ouverte est donnée par $BP_{S'}^-/P_{S'}^-$, elle est affine, donc son complémentaire 
  est de codimension $1$ et les composantes irréductibles de ce complémentaire sont les diviseurs premiers $B$-stables, 
  qui sont exactement les couleurs de $G/P_{S'}^-$.
  La décomposition de Bruhat nous dit, que ces couleurs sont de la forme 
  \[
    D_\alpha(G/P_{S'}^-) = \overline{Bs_\alpha P_{S'}^-/P_{S'}^-},
  \]
  où $\alpha$ est une racine simple de $G$, mais pas de $P_{S'}^-$ et $s_\alpha$ est 
  un représentant de la réflexion correspondante dans le groupe de Weyl de $G$.

  Le parabolique $P_{S'}$ fixe les diviseurs stables par $B$, donc $S_{G/P_{S'}^-}^p = S'$.
  Comme $T_z G/P_{S'}^- = T_z G z$ pour tout point fixe $z$ de $B$,
  nous avons que $\Sigma_{G/P_{S'}^-} = \emptyset$, donc $A_{G, G/P_{S'}^-} = \emptyset$ aussi.
  Finalement, le système sphérique de $G/P_{S'}^-$ est donné par $(S', \emptyset, \emptyset)$.
\end{exemple}

Toute sous-variété $G$-stable fermée d'une variété magnifique $X$ est encore magnifique.
Supposons, que $X_I \subset X$ est une sous-variété $G$-stable fermée. On a
$S_{X_I}^p = S_X^p$. Comme $\partial X = \cup_{i = 1}^{\rk(X)} X_i$ est un diviseur à croisements normaux lisses,
$\Sigma_X$ est aussi l'ensemble des $T$-poids sur $T_z X/T_z X_i$ pour tout $i \in \{1, \dots, \rk(X)\}$.
Cela signifie que $\Sigma_{X_I} = \{\gamma_i \in \Sigma_X\ |\ i \notin I\}$
et $A_{G, X_I}$ est formé des couleurs déplacées par les racines simples dans $S_G \cap \Sigma_{X_I}$.
De plus il existe une application $X_I \to G/P_{S_I}$ et si $X^I$ est la fibre au point de base, nous avons 
$X_I = G \times^{P_{S_I}} X^I$, où $P_{S_I}$ est le parabolique associé à
$S_I = S_G \setminus \big(S_X^p \cup \bigcup_{i \in I}\mathrm{supp}(\gamma_i)\big)$.

\begin{exemple}[Exemple 2]
  Supposons que $G$ est un groupe algébrique simple de type adjoint de rang $\rk(G) = r$.
  Le groupe $G$ peut être considéré comme une variété sphérique, sous l'action de $G \times G$,
  donnée par $(g,g')  x = gx(g')^{-1}$.
  Le groupe $B \times B^-$ est un sous-groupe de Borel de $G \times G$ et l'ensemble $BB^- \subset G$ est
  la $B \times B^-$-orbite ouverte de $id \in G$.
  Cela signifie que $G$ est sphérique. Par ailleurs, il admet une compactification magnifique $X$.

  Soit $S_{G\times G} = \{\alpha_1, \dots, \alpha_r, \beta_1, \dots, \beta_r\}$
  une base de racines simples de $G$, cor\-res\-pondant à $B\times B^-$, 
  telle que $\alpha_i$ sont les racines simples du premier 
  facteur de $G\times G$ et les $\beta_i$ du deuxième pour $i \in \{1, \dots, r\}$.
  Soit $z \in Y \subset X$ l'unique point fixe de $B\times B^-$ dans l'orbite fermée de $X$,
  qui est isomorphe à $G/B^- \times G/B$. Alors 
  $T_z X / T_z (G\times G  z) \simeq \mathfrak{t}$,
  où $\mathfrak{t}$ est l'algèbre de Lie du tore maximal $T = B\cap B^-$.
  Le tore $T\times T$ agit sur $\mathfrak{t}$ linéairement avec des poids $\gamma_i = \alpha_i - \beta_i$.
  Ce sont les racines sphériques $\gamma_i \in \Sigma_X$ de $X$.

  Les diviseurs premiers $B\times B^-$-stables sont les adhérences des cellules de Bruhat de codimension un.
  Ce sont les couleurs de $X$ et chaque couleur est déplacée par précisément deux racines simples de $G\times G$:
  $\alpha_i$ et $\beta_i$. On a $S^p_X = \emptyset$. Comme 
  $\Sigma_X \cap S_{G\times G} = \emptyset$, on a que $A_{G\times G, X}$ est vide aussi.
  Le système sphérique de la compactification magnifique de $G$ est donné par
  $(\emptyset, \{\alpha_i - \beta_i\ |\ 1 \leq i \leq r\}, \emptyset)$.
\end{exemple}

\section{Diviseurs, courbes et courbes rationnelles}

Pour une variété sphérique simple projective $X$, les couleurs qui ne contiennent 
pas l'orbite fermée $Y$ forment une base du groupe de Picard $\Pic(X)$:
\[
  \Pic(X) = \bigoplus_{D \in \Delta_X, Y \not\subset D} \mathbb{Z}D
\] 
On note $N_1(X)$ le groupe de classes de courbes modulo l'équivalence numérique.
Si $X$ est sphérique complète, alors $N_1(X)$ est le dual de $\Pic(X)$, voir Propositions 3.1 et 3.2 \cite{brion_candisv}.
Nous allons utiliser la notation
\begin{align*}
  \langle\ ,\ \rangle: \Pic(X) \times N_1(X) \to \mathbb{Z}
\end{align*}
pour le produit d'intersection entre courbes et diviseurs.
Sur une variété magnifique il existe une courbe $C_D$ pour toute couleur $D \in \Delta_X$,
qui est irréductible, $B$-stable et passe par l'unique point fixe de $B^-$ dans l'orbite fermée de $X$
et telle que pour toute couleur $D$ on a $\langle D, [C_D] \rangle = 1$ et
$\langle D, [C_D'] \rangle = 0$ pour tout $D' \neq D$, voir \cite{luna_gcplvs}.

\begin{exemple}[Exemple 3]
  Soit $G/P^-$ une variété projective homogène. Nous savons déjà, que les couleurs de $G/P^-$
  sont données par les cellules de Bruhat de codimension $1$, donc les couleurs de $G/P^-$ sont les diviseurs de Schubert.
  Les courbes $C_D$ sont les courbes de Schubert de $G/P^-$, données par $\overline{P_\alpha P^- /P^-}$,
  où $P_\alpha$ est le parabolique minimal, correspondant à une racine simple $\alpha$ de $G$,
  mais pas de $P$.
\end{exemple}

Nous sommes intéressés par les courbes rationnelles sur les variétés magnifiques.
Ces courbes ont déjà été considérées. Par exemple, Michel Brion et Baohua Fu dans \cite{brion_fu_mrc}
ont étudié les courbes rationnelles minimales, passant par l'identité $x \in X$ dans la
compactification magnifique d'un groupe semisimple de type adjoint. 

Nous introduisons les définitions comme dans \cite{brion_fu_mrc}.
Soit $\mathrm{RatCurves}^n(X)$ la normalisation de l'espace des courbes rationnelles sur $X$.
Il existe une famille universelle $\mathcal{U}$ avec des projections $\rho: \mathcal{U} \to \mathcal{K}$
et $\mu: \mathcal{U} \to X$, où $\mathcal{K} \subset \mathrm{RatCurves}^n(X)$ est une composante irréductible.
Pour tout $x \in X$ définissons $\mathcal{U}_x = \mu^{-1}(x)$ et $\mathcal{K}_x = \rho(\mathcal{U}_x)$.
La famille $\mathcal{K}$ est dite famille des courbes rationnelles minimales si $\mathcal{K}_x$ est non-vide et projective 
pour un point $x$ général.
Il existe une application rationnelle $\tau : \mathcal{K}_x \dashrightarrow \mathbb{P}T_x(X)$ qui envoie toute courbe
lisse en $x$ sur sa direction tangente. L'adhérence de l'image de $\tau$ est notée $\mathcal{C}_x$ 
et appelée la \textit{variété des tangentes rationnelles minimales} (VMRT) en $x$.
Si $X = \overline{G\times G/G}$ est la compactification magnifique d'un groupe algébrique $G$ simple de type adjoint avec 
algèbre de Lie $\mathfrak{g}$, alors $T_x(X) = \mathfrak{g}$ pour le point de base $x$ de la $G \times G$-orbite ouverte de $X$.
Le groupe $G$ agit sur $T_x(X)$ par le plongement diagonal et l'action induite sur $\mathbb{P} \mathfrak{g}$ stabilise
$\mathcal{C}_x$.

Le résultat suivant est prouvé dans \cite{brion_fu_mrc}:
\begin{theoremfr*}[{\cite[Theorem 1.1]{brion_fu_mrc}}]
  Soit $X$ la compactification magnifique d'un groupe algébrique simple $G$ de type adjoint.
  Alors 
  \begin{enumerate}
    \item Il existe une famille unique de courbes rationnelles minimales
      $\mathcal{K}$ sur $X$. De plus, la famille $\mathcal{K}_x$ des courbes passant par $x$ est lisse et
      l'application rationnelle $\tau: \mathcal{K}_x \dashrightarrow \mathcal{C}_x$
      est un isomorphisme.
    \item $\mathcal{C}_x$ est la $G$-orbite fermée dans $\mathbb{P}\mathfrak{g}$,
      si $G$ n'est pas du type $A$.
    \item Quand $G$ est de type $A_l$ $(l \geq 2)$, c'est-à dire $G = \mathrm{PGL}(V)$
      pour un espace vectoriel $V$ de dimension $l+1$, la VMRT $\mathcal{C}_x$ est l'image de
      $\mathbb{P}V \times \mathbb{P}V^*$ sous le plongement de Segre
      $\mathbb{P}V \times \mathbb{P}V^* \to \mathbb{P}\mathrm{End}(V)$, suivi par la projection
      $\mathbb{P}\mathrm{End}(V)\dashrightarrow \mathbb{P}(\mathrm{End}(V)/ \mathbb{C}id) = \mathbb{P}\mathfrak{g}$
      de centre $x = [\rm{id}]$.
  \end{enumerate}
\end{theoremfr*}

Un autre objet important, que nous présentons dans ce texte, est la compactification de Deligne-Mumford 
de l'espace de modules de courbes rationnelles de classe $\eta \in N_1(X)$ avec $n$ points marqués. 
Cet espace est noté $\overline{M_{0,n}}(X, \eta)$, on l'appelle aussi l'espace de Kontsevich de modules 
d'applications stables. Cette compactification contient un ensemble ouvert composé des morphismes 
$\{f:\mathbb{P}^1 \to X\ |\ f_*[\mathbb{P}^1] = \eta\}$ et l'étend avec des applications stables.
Une application stable de genre $0$ et de classe $\eta$ avec $n$ points marqués est la donnée
$(C, p_1, \dots, p_n, f)$,
où $C$ est une courbe rationnelle connexe nodale de genre $0$, où $p_1, \dots, p_n$ sont 
des points réguliers distincts sur $C$ et où $f: C \to X$ est un morphisme, tel que
$f_*[C] = \eta$ et qui vérifie la condition suivante pour chaque composante irréductible $E \subset C$:
si $E \simeq \mathbb{P}^1$ et $f(E)$ est un point, alors $E$ doit contenir au moins trois points spéciaux 
(soit les points marqués soit les intersections avec les autres composantes de $C$).

\begin{theoremfr*}[{\cite[Theorem 1]{fulton_pandharipande_nosmaqc}}]
  Soit $X$ une variété complexe projective lisse. Il existe un espace de module grossier projectif
  $\overline{M}_{0,n}(X,\eta)$ contenant $M_{0,n}(X, \eta)$ comme un sous-ensemble ouvert.
\end{theoremfr*}

\begin{theoremfr*}[{\cite[Theorem 2]{fulton_pandharipande_nosmaqc}}]
  Soit $X$ une variété projective lisse telle que $H^1(\mathbb{P}^1, \mu^*(T_X)) = 0$ pour toute application $\mu: \mathbb{P}^1 \to X$.
  \begin{itemize}
    \item $\overline{M}_{0,n}(X, \eta)$ est une variété normale, projective de la dimension pure
    \[
      \dim(X) - \langle K_X, \eta \rangle + n - 3.
    \]
    \item $\overline{M}_{0,n}(X, \eta)$ est localement un quotient d'une variété lisse par un groupe fini.
  \end{itemize}
\end{theoremfr*}

Il est déjà connu, que l'espace de module est irréductible dans le cas d'espaces homogènes projectifs.

\begin{theoremfr*}[{\cite[Theorem 1]{thomsen_im0ngp}}]
  Soit $X$ une variété des drapeaux, soit $\eta \in N_1(X)$ une classe effective. Alors 
  l'espace $\overline{M}_{0,n}(X,\eta)$ est irréductible pour tout $n \in \mathbb{N}$ entier positif.
\end{theoremfr*}

Une de nos motivations dans cette thèse était de voir si ce résultat s'étend au cas des variétés magnifiques.

\section{Application limite de Luna}

Nous décrivons les outils combinatoires, utilisés pour prouver notre résultat.
Soit $X$ une $G$-variété projective lisse.
Soit $\lambda:\mathbb{C}^* \to G$ un sous-groupe à un paramètre de $G$.
Notons $X^\lambda$ la composante connexe de l'ensemble des points fixes $X^{\lambda(\mathbb{C})^*}$, 
telle que l'ensemble $X_\lambda = \{x \in X\ |\ \lim_{t \to 0}\lambda(t)x \in X^\lambda\}$ 
est ouvert dans $X$. 
L'ensemble des points fixes $X^\lambda$ est lisse.
Définissons aussi le sous-groupe parabolique 
$G_\lambda = \{g \in G\ |\ \lim_{t\to 0}\lambda(t)g\lambda(t)^{-1} \in G^\lambda\}$,
où $G^\lambda$ est le centralisateur dans $G$ de $\lambda(\mathbb{C}^*)$ et est un sous-groupe 
de Levi de $G_\lambda$.
Sur la cellule ouverte $X_\lambda$ de cette décomposition il existe
une application limite
\begin{align*}
  \varphi_\lambda: X_\lambda &\to X^\lambda \\
  x &\mapsto \lim_{t \to 0}\lambda(t)x 
\end{align*}
qui est une fibration localement triviale en espaces affines.
Cela produit une application rationnelle
$\varphi_\lambda: X \dashrightarrow X^\lambda \subset X$. 
Pour toute application rationnelle entre variétés projectives lisses
le lieu d'indétermination est de codimension au moins deux, donc pour une courbe 
assez positive $C$ sur $X$, nous pouvons supposer, qu'elle évite ce lieu d'indétermination
et nous pouvons considérer la courbe $\varphi_\lambda(C)$ dans $X^\lambda$.
Cela va nous permettre de procéder par récurrence sur la dimension des variétés magnifiques.


\section{Variétés magnifiques symétriques}

Les résultats de cette thèse portent sur une classe spéciale de variétés magnifiques, les 
compactifications magnifiques des espaces symétriques, introduites par Corrado De Concini et Claudio Procesi 
dans \cite{deconcini_procesi_csv}. Un espace symétrique est une variété algébrique, isomorphe à
un quotient $G/G^\sigma$ pour une involution $\sigma: G \to G$, où $G$ est est un groupe réductif.
Les espaces symétriques sont sphériques.

Soit $G$ un groupe algébrique semisimple de type adjoint, alors pour tout espace symétrique
$G/G^\sigma$ il existe une compactification magnifique $X = \overline{G/G^\sigma}$.
Un tel $X$ est appelé une variété magnifique symétrique.

Dans ce cas les racines sphériques de $X$ sont données par $\bar{\alpha}_i = \alpha_i - \sigma(\alpha_i)$
pour $1 \leq i \leq \rk(X)$. 
Il faut noter qu'il est possible que $\rk(X) < \rk(G)$ et
qu'il existe $1 \leq i \neq j \leq \rk(X)$, tels que $\bar{\alpha}_i = \bar{\alpha}_j$.

Pour une variété symétrique magnifique $X$ il existe un sous-ensemble ouvert affine $V$,
tel que $V \simeq U \times \mathbb{A}^l$, où $U = R_u(P_{S_X^p})$.
Nous écrivons $(u, x) \in V$ avec $x = (x_1, \dots, x_l)$.
Pour ce $V$ tout diviseur de bord $X_i$ de $X$ est localement donné par 
\[
  X_i \cap V = \{(u, x)\ |\ x_i = 0\}.
\]

Soit $I \subset \{1, \dots, l\}$ et soit $X_I$ une variété $G$-stable fermée de $X$.
La variété $X_I$ est de la forme $G\times^{P^-} X^I$, où $P$ est le sous-groupe parabolique,
associé à $S_I = \{\alpha \in S\ |\ \bar{\alpha} = \bar{\alpha}_i \text{ for some $i \notin I$}\}$.

Le groupe de Picard $\Pic(X) \simeq \mathbb{Z}^{l+s}$, où $s$ est le nombre de racines sphériques $\bar{\alpha}$, 
telles qu'il existe deux racines simples distinctes $\alpha, \alpha' \in S$, avec
soit $\alpha \neq -\sigma(\alpha')$ soit $\alpha \neq -\sigma(\alpha')$ et $(\alpha, \alpha') \neq 0$.

De plus pour toute variété $G$-stable fermée $X_I$ il existe une suite exacte courte
\[
  \begin{tikzcd}
    0 \arrow[r] & \Pic(G/P_{S_I}^-) \arrow[r, "\pi_I^*"] & \Pic(X_I) \arrow[r,"\iota^*"] & \Pic(X^I) \arrow[r] & 0
  \end{tikzcd}
\]
où $\pi_I : X_I \to G/P_{S_I}^-$ est la projection $G$-équivariante et $\iota: X^I \to X_I$ est
l'inclusion de la fibre au point $P_{S_I}^-$-fixe de $G/P_{S_I}^-$.

\section{Réductibilité de l'espace de modules}

Soit $X$ une variété symétrique magnifique et soit $\eta\in N_1(X)\setminus \{0\}$ une classe de courbes mobiles.
Pour les variétés magnifiques symétriques, cette condition est équivalente
au fait, que $\eta$ a une intersection positive avec tous les diviseurs effectifs de $X$.
Nous considérons l'espace de module $\overline{M}_{0,0}(X, \eta)$ des applications stables de genre $0$ sans point marqué
et définissons le sous-ensemble ouvert $M_{0,0}^\circ(X,\eta)$ de courbes irréductibles rencontrant
la $G$-orbite ouverte de $X$.
Il se trouve que l'ensemble $M_{0,0}^\circ(X,\eta)$ est non-vide si $\eta = 2\eta'$ pour 
une classe $\eta' \in N_1(X)$, positive sur les diviseurs effectifs de $X$.
Nous prouvons aussi le théorème suivant, qui montre que dans beaucoup de cas l'espace $\overline{M}_{0,0}(X,\eta)$ est réductible.

\begin{customthmfr}{A}\label{thm_a_fr}
  Soit $X$ une variété magnifique symétrique.
  Soit $\eta \in N_1(X)\setminus \{0\}$ une classe de courbes telle que $M_{0,0}^\circ(X, \eta) \neq \emptyset$.
  Supposons qu'il existe une décomposition $\eta = \eta_1 + \eta_2$, telle que
  $\eta_1$ et $\eta_2$ sont des classes de courbes sur $X$ non-triviales effectives
  et il existe un diviseur de bord $D$ tel que
  \[
    \langle D, \eta_2 \rangle \leq -2.  
  \]

  Alors l'espace de modules $\overline{M}_{0,0}(X, \eta)$ est réductible.
\end{customthmfr}

En particulier nous obtenons le résultat suivant:
\begin{customthmfr}{B}\label{thm_b_fr}
  Soit $X$ la compactification magnifique d'un groupe simple de type adjoint $G$ de rang $\rk(G) \geq 3$.
  Alors $\overline{M}_{0,0}(X, \eta)$ est réductible pour toute classe de courbes mobiles $\eta$.
\end{customthmfr}
Ce résultat constraste avec le fait, que la VMRT
est irréductible pour une compactification magnifique d'un groupe.

Ce résultat motive aussi la partie finale de la thèse et la recherche d’une composante 
irréductible spéciale de $\overline{M}_{0,0}(X, \eta)$. 

\section{Irréductibilité d'une composante spéciale}

Soit $X$ la compactification magnifique d'un groupe $G$ algébrique semisimple de type adjoint
de rang $\rk(G) = r$.
Tout co-poids fondamental $\omega_{i_0}$ pour un $i_0 \in \{1, \dots, r\}$ de $G$ définit un sous-groupe 
à un paramètre de $G$ et aussi de $G\times G$ avec l'inclusion $g \mapsto (g, 1_G)$.
Il existe une application limite $\varphi_{i_0}:X \dashrightarrow X$ avec $\varphi_{i_0}(X) \subset X_{i_0}$. 
Soit $\eta$ une classe d'équivalence numérique de courbes mobiles.
On montre que $M_{0,0}^\circ(X,\eta)$ n'est pas vide
et qu'il contient un sous-ensemble ouvert de courbes, évitant le lieu d'indétermination de $\varphi_{i_0}$.
Donc il existe une application rationnelle
\begin{align*}
  \varphi_{i_0}: M_{0,0}^\circ(X, \eta) &\dashrightarrow M_{0,0}^\circ(X_{i_0}, \bar{\eta})\\
    C &\mapsto \varphi_{i_0}(C),
\end{align*}
où $\bar{\eta}$ est la classe de $\varphi_{i_0}(C)$ dans $X_{i_0}$.
Une assertion similaire est vraie pour les sous-variétés $G$-stables fermées de $X$.
Nous procédons par récurrence pour montrer:

\begin{customthmfr}{C}\label{thm_c_fr}
  Soit $X$ la compactification magnifique d'un groupe algébrique semi\-simple $G$ de type adjoint.
  Soit $\eta \in N_1(X)$ une classe de courbes mobiles. 
  Alors l'espace $\overline{M_{0,0}^\circ(X,\eta)}$ est irréductible.
\end{customthmfr}

Une suite possible à ce travail serait d'élargir la classe de variétés, pour les\-quelles l'asser\-tion sur
l'irréductibilité est vraie, mais nous pourrions aussi con\-si\-dé\-rer la 
connexité rationnelle simple des variétés magnifiques.
Une variété $X$ est dite \textit{rationnellement connexe}, s'il existe une courbe rationnelle irréductible,
passant par deux points généraux quelconques. Un premier exemple est le cas de l'espace projectif $\mathbb{P}^N$. 

La \textit{simple connexité rationnelle} est plus complexe, il faut que l'espace de courbes rationnelles soit
rationnellement connexe. Puisqu'il y a beaucoup d'espa\-ces de modules (pour chaque classe $\eta \in N_1(X)$)
et que ces espaces ont de nombreuses composantes irréductibles, comme nous l'avons vu au Théorème \ref{thm_a_fr},
une bonne définition de la simple connexité rationnelle a besoin de l'existence d'une 
composante canonique de l'espace de modules d'applications stables d'une classe $\eta$ 
assez positive, et que celle-ci soit 
rationnellement connexe. La composante étudiée dans le Théorème \ref{thm_b_fr} est une 
bonne candidate pour être cette "composante canonique".

\chapter{Introduction}

\section{Motivation}
Rational curves have always played an important role in complex algebraic geometry, be it the minimal model program (MMP),
the variety of minimal rational tangents (VMRT), the new notion of simple rational connectedness
or enumerative geometry, like quantum cohomology. In many of these cases one of the
steps to results would be understanding the moduli space of rational curves on a given variety.
A lot has been done on this topic already, for example \cite{fulton_pandharipande_nosmaqc}
explains the construction of the moduli space of stable maps and that it is projective
in case the underlying variety is a projective variety.

Our main focus in this PhD thesis are wonderful $G$-varieties. 
This choice is based on the fact that these varieties have a very nice structure of $G$-orbits.
This allows us to control very well the divisor class group $\Pic(X)$ and its dual, 
the group of curves $N_1(X)$ modulo numerical equivalence.

\section{Wonderful varieties}

We will work over $\mathbb{C}$.
Let $G$ be a connected reductive algebraic group of rank $\rk(G) = r$. 
Fix a Borel subgroup $B \subset G$ and a maximal torus $T$ of $G$ with $T \subset B$. Denote by $R = R_G$ 
the root system of $(G,T)$, by $R_G^+$ the roots of $B$, by $R_G^- = -R_G^+$ and by $S = S_G$ 
the basis of simple roots.
A parabolic subgroup $P$ of $G$ is standard if it contains $B$,
its opposite parabolic containing $T$ is denoted by $P^-$.
Let $S' \subset S$, then $R_{G,S'}$ denotes the set of roots, generated by $S'$
and $P_{S'}$ is the parabolic subgroup of $G$ with roots $R^+_G \cup R_{G,S'}$.
The parabolic $P_{S'}$ is said to be associated to $S'$.

A normal $G$-variety $X$ with an open $B$-orbit is called \textit{spherical}.
If $X$ contains exactly one closed $G$-orbit, it is also called \textit{simple}.
Such a variety has an open $G$-orbit, that is isomorphic to the spherical 
homogeneous space $G/H$ for some subgroup $H \subset G$,
which is also called \text{spherical}.
In the setting of algebraic geometry this implies that there are finitely many 
$B$ and $G$-orbits in $X$.
We can define the set of $B$- and not $G$-stable prime divisors, called \textit{colors}
and denoted by $D \in \Delta_X$.
It turns out that a projective simple spherical $G$-variety
$X$ is uniquely defined up to isomorphism by the cone of $G$-invariant valuations
and by the subset of colors that contain the closed $G$-orbit of $X$. 
This datum is called the \textit{colored cone} of $X$. 
A spherical homogeneous space is called \textit{sober}, if the cone of $G$-invariant
valuations is strictly convex or, equivalently, does not contain a linear subspace.
A simple spherical variety $X$ is called \textit{toroidal} if its $G$-closed orbit $Y$
is not contained in any color $D \in \Delta_X$. 
A combination of colored cones with special conditions, 
called colored fan can be used to classify nonsimple spherical varieties
in the same spirit as the classification of toric varieties.

A \textit{wonderful} variety is a smooth complete $G$-variety with an open $G$-orbit $G/H$,
such that its complement $\partial X = X \setminus (G/H)$
is a smooth normal crossing divisor. Furthermore the closures of $G$-orbits of $X$ are exactly the partial 
intersections of irreducible components of $\partial X$. In particular there is a unique 
closed $G$-orbit, the intersection of all components. The number 
of irreducible components of $\partial X$ is called the \textit{rank} of $X$,
denoted by $\rk(X)$. Then we can write $\partial X = \cup_{i = 1}^{\rk(X)} X_i$,
where every $X_i$ is a $G$-stable prime divisor. Such $G$-stable prime 
divisors on wonderful varieties will be called \textit{boundary divisors}.
We will write $X_I = \cap_{i \in I} X_i$
for any $G$-stable subvariety of $X$, where $I \subset \{1,\dots, \rk(X)\}$.

Domingo Luna proved in \cite{luna_tvmes} that wonderful varieties are toroidal spherical varieties.
We are going to describe the combinatorics, used in the classification of wonderful varieties.
We say that a color $D \in \Delta_X$ is \textit{moved} by a simple root $\alpha$ if
$D$ is not stable under the action of the minimal parabolic subgroup $P_\alpha$.
The set of colors moved by a simple root $\alpha$ is denoted by $\Delta_X(\alpha)$.
We will denote by $S_X^p$ the set of all simple roots of $G$, such that $\Delta_X(\alpha) = \emptyset$.
The \textit{spherical roots} $\Sigma_X$ are the $T$-weights of the normal space
\[
  \frac{T_z X}{T_z G z},
\]
where $z$ is the $B^-$ fixed point of the closed orbit of $X$. 
Let $A_{G,X} \subset \Delta_X$ be the set of colors $D \in \Delta_X$, which are
moved by some simple root $\alpha \in S_G \cap \Sigma_X$, that is also a spherical root of $X$.
Any $\gamma \in \Sigma_X$ is a linear combination of simple roots of $G$ with nonnegative coefficients. 
We define the support of $\gamma$, denoted by $\mathrm{supp}(\gamma)$ to be the simple 
roots having a nonzero coefficient in the presentation of $\gamma$.

The \textit{spherical system} of $X$ is defined as $(S_X^p, \Sigma_X, A_{G,X})$.
It was previously conjectured and proved in type $A$ by Domingo Luna in \cite{luna_vsdta}, that any 
wonderful variety is defined uniquely up to isomorphism by its spherical system.
Recently Paolo Bravi and Guido Pezzini proved it in \cite{bravi_pezzini_pwv}:

\begin{theorem*}[{\cite[Theorem 1.2.3]{bravi_pezzini_pwv}}]
  There is a bijection between isomorphism classes of wonderful $G$-varieties and 
  spherical $G$-systems.
\end{theorem*}

They have introduced the notion of an abstract spherical $G$-system 
and shown that every abstract spherical $G$-system is a spherical system of some wonderful $G$-variety.

\begin{example}[Example 1]
  Projective rational homogeneous varieties are wonderful of trivial rank.
  Let $P_{S'}$ be a standard parabolic subgroup of $G$, associated to some $S' \subset S_G$. Then 
  the projective homogeneous variety $G/P_{S'}^-$ has an open $B$-orbit, where 
  $P_{S'}^-$ denotes the parabolic subgroup, opposite of $P_{S'}$.
  The open $B$-orbit of $G/P_{S'}^-$ is given by $BP_{S'}^-/P_{S'}^-$, it is 
  affine, so its complement is of codimension one and its irreducible components
  are the $B$-stable prime divisors, which are precisely the colors of $G/P_{S'}^-$.
  Bruhat decomposition tells us that these colors are of the form 
  \[
    D_\alpha(G/P_{S'}^-) = \overline{Bs_\alpha P_{S'}^-/P_{S'}^-},
  \]
  where $\alpha$ is a simple root of $G$, which is not a simple root of $P_{S'}^-$ and $s_\alpha$ is 
  a representative of the corresponding reflection in the Weyl group of $G$.

  The parabolic $P$ fixes the $B$-stable divisors,
  so $S_{G/P^-}^p = S'$. Since $T_z G/P_{S'}^- = T_z G z$ for 
  any fixed point $z$ of $B$, we get that $\Sigma_{G/P_{S'}^-} = \emptyset$
  so $A_{G, G/P_{S'}^-} = \emptyset$ as well.
  Altogether, the spherical system of $G/P_{S'}^-$
  is given by $(S', \emptyset, \emptyset)$.
\end{example}

Any closed $G$-stable subvariety of a wonderful variety $X$ is again wonderful. 
Suppose that $X_I \subset X$ is a closed $G$-stable subvariety. Then 
$S_{X_I}^p = S_X^p$. As $\partial X = \cup_{i = 1}^{\rk(X)} X_i$ is a simple normal crossings divisor, 
$\Sigma_X$ is also the set of $T$-weights on $T_z X/T_z X_i$ for all $i \in \{1, \dots, \rk(X)\}$.
This means that $\Sigma_{X_I} = \{\gamma_i \in \Sigma_X\ |\ i \notin I\}$
and $A_{G,X_I}$ are the colors moved by simple roots in $S_G \cap \Sigma_{X_I}$.
Additionally there exists a $G$-equivariant map $X_I \to G/P_{S_I}$ and if $X^I$ is the fiber
at the base point of $G/P_{S_I}$,
we have $X_I = G \times^{P_{S_I}} X^I$, 
where $P_{S_I}$ is the parabolic associated to 
$S_I = S_G \setminus\big(S_X^p \cup \bigcup_{i \in I}\mathrm{supp}(\gamma_i)\big)$.

\begin{example}[Example 2]
  Suppose that $G$ is a simple algebraic group of adjoint type of rank $\rk(G) = r$.
  The group $G$ can be considered as a spherical variety itself, under the action of $G \times G$,
  given by $(g,g')  x = gx(g')^{-1}$.
  Then $B\times B^-$ is a Borel subgroup of $G\times G$ and the set $BB^- \subset G$ is 
  the open $B \times B^-$-orbit of $id \in G$. 
  That means that we can consider the wonderful compactification $X$ of $G$. 

  Let $S_{G\times G} = \{\alpha_1, \dots, \alpha_r, \beta_1,\dots, \beta_r\}$ 
  be a basis of simple roots of $G$ corresponding to $B\times B^-$, 
  such that $\alpha_i$ denote the simple roots of the first factor of $G\times G$ 
  and $\beta_i$ those of the second for $i \in \{1,\dots, r\}$.
  Let $z \in Y \subset X$ be the unique fixed point of $B\times B^-$ in the closed orbit of $X$,
  which is isomorphic to $G/B^- \times G/B$. Then 
  $T_z X / T_z (G\times G  z) \simeq \mathfrak{t}$,
  where $\mathfrak{t}$ is the Lie algebra of the maximal torus $T = B\cap B^-$.
  Then $T\times T$ acts on $\mathfrak{t}$ linearly with the weights $\gamma_i = \alpha_i - \beta_i$.
  They are precisely the spherical roots $\gamma_i \in \Sigma_X$.
  
  The $B\times B^-$ stable divisors of $X$ are the closures of the codimension one Bruhat cells.
  These are precisely the colors of $X$ and every color is moved by exactly two orthogonal simple roots 
  of $G\times G$, namely $\alpha_i$ and $\beta_i$. This means that $S^p_X = \emptyset$.
  Notice that $\Sigma_X \cap S_{G\times G} = \emptyset$, the set $A_{G\times G,X}$ is empty as well.
  Then the spherical system of the wonderful compactification of $G$ is given by 
  $(\emptyset, \{\alpha_i - \beta_i\ |\ 1 \leq i \leq r\}, \emptyset)$.
\end{example}

\section{Divisors, curves and rational curves}
Let $X$ be a simple spherical projective variety, 
the Picard group $\Pic(X)$ is the free abelian group 
on the set of colors $\Delta_X$, which do not contain the closed orbit $Y\subset X$:
\[
  \Pic(X) = \bigoplus_{D \in \Delta_X, Y \not\subset D}\mathbb{Z}D.
\]
We denote by $N_1(X)$ the group of classes of curves modulo numerical equivalence.
If $X$ is spherical and complete, then $N_1(X)$ is dual to $\Pic(X)$, see Propositions 3.1 and 3.2 \cite{brion_candisv}. 
From here on we will use the notation 
\begin{align*}
  \langle\ ,\ \rangle: \Pic(X) \times N_1(X) \to \mathbb{Z}
\end{align*} for the intersection pairing between curves and divisors.
On wonderful varieties there exists a curve $C_D$ for any $D\in \Delta_X$, which is irreducible, 
$B$-stable and passes through the unique 
$B^-$ fixed point in the closed orbit of $X$
and such that for any color $D$ we have that $\langle D, [C_D] \rangle = 1$ and
$\langle D, [C_D'] \rangle = 0$ for any $D' \neq D$, see \cite{luna_gcplvs}.

\begin{example}[Example 3]
  Let $G/P^-$ be a projective homogeneous variety, where $G$ is reductive. We have seen that 
  colors of $G/P^-$ are given by Bruhat cells of codimension $1$, which means that
  the colors of $G/P^-$ are the Schubert divisors. Then $C_D$ defined 
  previously are the Schubert curves of $G/P^-$, given by 
  $\overline{P_\alpha P^- /P^-}$, where $P_\alpha$ are the minimal parabolics,
  corresponding to the simple roots $\alpha$ of $G$, that are not simple roots of $P^-$.
\end{example}

We are interested in rational curves on wonderful varieties. 
These curves have been already considered. For example, 
Michel Brion and Baohua Fu in \cite{brion_fu_mrc} have closely studied minimal 
rational curves passing through the identity in wonderful compactifications
of groups. 

We introduce the definitions as in \cite{brion_fu_mrc}.
Let $\mathrm{RatCurves}^n(X)$ denote the normalization of the space of rational curves on $X$.
There exists a universal family $\mathcal{U}$ together with projections $\rho: \mathcal{U} \to \mathcal{K}$
and $\mu: \mathcal{U} \to X$, where $\mathcal{K} \subset \mathrm{RatCurves}^n(X)$ is an irreducible component.
The family $\mathcal{K}$ is called a family of minimal rational curves if $\mathcal{K}_x$ is nonempty and projective for a general point $x$.
For any $x \in X$ define $\mathcal{U}_x = \mu^{-1}(x)$ and $\mathcal{K}_x = \rho(\mathcal{U}_x)$.
There is a rational map $\tau : \mathcal{K}_x \dashrightarrow \mathbb{P}T_x(X)$ that sends any curve, which 
is smooth at $x$ to its tangent direction. The closure of the image of $\tau$ is denoted by $\mathcal{C}_x$ 
and called the \textit{variety of minimal rational tangents} (VMRT) at the point $x$.
If $X = \overline{G\times G/G}$ is the wonderful compactification of a simple algebraic group $G$ of adjoint type with 
Lie algebra $\mathfrak{g}$, then $T_x(X) = \mathfrak{g}$ for the base point $x$ in the open orbit of $X$.
Then $G$ acts on $T_x(X)$ via the diagonal embedding and the induced action on $\mathbb{P} \mathfrak{g}$ stabilizes 
$\mathcal{C}_x$.

The following important result is proved in \cite{brion_fu_mrc}:
\begin{theorem*}[{\cite[Theorem 1.1]{brion_fu_mrc}}]
  Let $X$ be the wonderful compactification of a simple algebraic group 
  $G$ of adjoint type. Then 
  \begin{enumerate}
    \item There exists a unique family of minimal rational curves 
      $\mathcal{K}$ on $X$. Moreover, $\mathcal{K}_x$ is smooth and 
      the rational map $\tau: \mathcal{K}_x \dashrightarrow \mathcal{C}_x$
      is an isomorphism.
    \item $\mathcal{C}_x$ is the closed $G$-orbit in $\mathbb{P}\mathfrak{g}$,
      if $G$ is not of type $A$.
    \item When $G$ is of type $A_l$ $(l \geq 2)$, so that $G = \mathrm{PGL}(V)$
      for a vector space $V$ of dimension $l+1$, the VMRT $\mathcal{C}_x$ is the image
      of $\mathbb{P}V \times \mathbb{P}V^*$ under the Segre embedding 
      $\mathbb{P}V \times \mathbb{P}V^* \to \mathbb{P}\mathrm{End}(V)$, followed by the projection
      $\mathbb{P}\mathrm{End}(V)\dashrightarrow \mathbb{P}(\mathrm{End}(V)/ \mathbb{C}id) = \mathbb{P}\mathfrak{g}$
      with center $x = [\rm{id}]$.
  \end{enumerate}
  In particular the VMRT is irreducible.
\end{theorem*}

Another important object we introduce is the Deligne-Mumford compactification of 
the moduli space of rational curves of class $\eta  \in N_1(X)$ with $n$ marked points.
This space is denoted by $\overline{M}_{0,n}(X, \eta)$, sometimes also called the Kontsevich 
moduli space of stable maps.
This compactification takes as a base open subset 
the set of morphisms $\{f:\mathbb{P}^1 \to X\ |\ f_*[\mathbb{P}^1] = \eta\}$ 
and expands it with stable maps.
A stable map of genus $0$ and class $\eta$ with $n$ marked points is the data 
$(C, p_1, \dots, p_n, f)$,
where $C$ is a rational nodal curve of genus $0$, $p_1, \dots, p_n$ 
are distinct smooth points on $C$ and $f: C \to X$ is a morphism, such that 
$f_*[C] = \eta$ with the following condition for any irreducible component $E \subset C$:
if $E \simeq \mathbb{P}^1$ and $E$ is mapped to a point by $f$, then $E$ must 
contain at least three special points (either markings or intersections with other 
components of $C$). 

\begin{theorem*}[{\cite[Theorem 1]{fulton_pandharipande_nosmaqc}}]
  Let $X$ be a smooth projective complex variety. There exists a coarse moduli space 
  $\overline{M}_{0,n}(X,\eta)$ containing $M_{0,n}(X, \eta)$ as a open subset.
\end{theorem*}

\begin{theorem*}[{\cite[Theorems 1 and 2]{fulton_pandharipande_nosmaqc}}]
  Let $X$ be a projective, smooth variety such that $H^1(\mathbb{P}^1, \mu^*(T_X)) = 0$ for every map $\mu: \mathbb{P}^1 \to X$.
  \begin{itemize}
    \item $\overline{M}_{0,n}(X, \eta)$ is a normal, projective variety of pure dimension
    \[
      \dim(X) - \langle K_X, \eta \rangle + n - 3.
    \]
    \item $\overline{M}_{0,n}(X, \eta)$ is locally a quotient of a smooth variety by a finite group.
  \end{itemize}
\end{theorem*}

A projective homogeneous space $G/P$ satisfies the conditions of this theorem. Additionally it is known that 
$\overline{M}_{0,n}(G/P, \eta)$ is irreducible.

\begin{theorem*}[{\cite[Theorem 1]{thomsen_im0ngp}}]
  Let $X$ be a flag variety and $\eta \in N_1(X)$ an effective class. The space
  $\overline{M}_{0,n}(X,\eta)$ is irreducible for every non negative integer $n$.
\end{theorem*}

One of our motivations for this PhD thesis was to see whether this is different in the case of wonderful varieties.

\section{Luna limit map}

Now we describe some of the combinatorial tools which we will use to prove our results.
Let $X$ be a smooth projective $G$-variety.
Let $\lambda:\mathbb{C}^* \to G$ be a $1$-parameter subgroup 
of $G$. Denote by $X^\lambda$ the connected component of the 
fixed point set $X^{\lambda(\mathbb{C})^*}$, such that the set 
$X_\lambda = \{x \in X\ |\ \lim_{t \to 0}\lambda(t)x \in X^\lambda\}$ is open in $X$.
The fixed point set $X^\lambda$ is smooth.
We can also define a parabolic subgroup 
$G_\lambda = \{g \in G\ |\ \lim_{t\to 0}\lambda(t)g\lambda(t)^{-1} \in G^\lambda\}$,
where $G^\lambda$ is the centralizer in $G$ of $\lambda(\mathbb{C}^*)$ and 
a Levi of $G_\lambda$.
Then on the open cell $X_\lambda$ of this 
decomposition there exists a limit map 
\begin{align*}
  \varphi_\lambda: X_\lambda &\to X^\lambda \\
  x &\mapsto \lim_{t \to 0}\lambda(t)x,
\end{align*} 
which is a locally trivial fibration in affine spaces.
This gives rise to a rational map
$\varphi_\lambda: X \dashrightarrow X^\lambda \subset X$. For any rational map between smooth projective 
varieties the indeterminacy locus 
is of codimension at least two, so for a positive enough curve $C$ on 
$X$ we can assume that it avoids this indeterminacy locus and we can consider
the curve $\varphi_\lambda(C)$ in $X^\lambda$.
This will later allow us to apply induction on the dimension of wonderful
varieties.


\section{Wonderful symmetric varieties}

Our results are concerned with a special class of wonderful varieties, the wonderful compactifications 
of symmetric spaces, introduced by Corrado De Concini and Claudio Procesi in \cite{deconcini_procesi_csv}.
A symmetric space is a $G$-homogeneous algebraic variety, isomorphic 
to a quotient $G/G^\sigma$ for some involution $\sigma: G \to G$, where $G$ is a reductive group.
Symmetric spaces are spherical. 

Let $G$ be a semisimple algebraic group of adjoint type, then for any symmetric 
space $G/G^\sigma$ there exists a wonderful compactification $X = \overline{G/G^\sigma}$.
Such $X$ is called \textit{wonderful symmetric variety}.

In this case the spherical roots of $X$ are given by $\bar{\alpha}_i = \alpha_i - \sigma(\alpha_i)$
for $1 \leq i \leq \rk(X)$. Note that it may appear that $\rk(X) < \rk(G)$ and 
there exist $1 \leq i \neq j \leq \rk(X)$, such that $\bar{\alpha}_i = \bar{\alpha}_j$.

For a wonderful symmetric variety $X$ there exists an open affine subset $V$, 
such that $V \simeq U \times \mathbb{A}^l$, where $U = R_u(P_{S_X^p})$.
We write $(u, x) \in V$ with $x = (x_1, \dots, x_l)$.
For this $V$ any boundary divisor $X_i$ of $X$ is locally given by
\[
  X_i \cap V = \{(u, x)\ |\ x_i = 0\}.
\]

Let $I \subset \{1, \dots, l\}$ and let $X_I$ be a closed $G$-stable subvariety of $X$.
Then $X_I$ is of the form $G\times^{P^-} X^I$, where $P$ is the parabolic subgroup,
associated to $S_I = \{\alpha \in S\ |\ \bar{\alpha} = \bar{\alpha}_i \text{ for some $i \notin I$}\}$.

The Picard group $\Pic(X) \simeq \mathbb{Z}^{l+s}$, where $s$ is the number of spherical roots $\bar{\alpha}$, 
such that there exist two distinct simple roots $\alpha, \alpha' \in S$, with
either $\alpha \neq -\sigma(\alpha')$ or $\alpha \neq -\sigma(\alpha')$ and $(\alpha, \alpha') \neq 0$.

Also for any closed $G$-stable subvariety $X_I$ there exists a short exact sequence
\[
  \begin{tikzcd}
    0 \arrow[r] & \Pic(G/P_{S_I}^-) \arrow[r, "\pi_I^*"] & \Pic(X_I) \arrow[r,"\iota^*"] & \Pic(X^I) \arrow[r] & 0
  \end{tikzcd}
\]
where $\pi_I : X_I \to G/P_{S_I}^-$ is the $G$-equivariant projection and $\iota: X^I \to X_I$ is the
inclusion of the fiber over the unique $P_{S_I}^-$-fixed point of $G/P_{S_I}^-$.



\section{Reducibility of the moduli space of stable maps}

Let $X$ be a wonderful symmetric variety and let $\eta\in N_1(X)\setminus \{0\}$ be a movable curve class.
For wonderful symmetric varieties this is equivalent to the fact that $\eta$ has nonnegative
intersection with all effective divisors of $X$.
We wish to consider the moduli space $\overline{M}_{0,0}(X, \eta)$ of stable maps of genus $0$ with no marked points.
We define an open subset $M_{0,0}^\circ (X, \eta)$ of irreducible curves meeting the open $G$-orbit of $X$.
It turns out that the set $M_{0,0}^\circ (X, \eta)$ is always nonempty if $\eta = 2\eta'$ for 
some class $\eta' \in N_1(X)$, nonnegative on effective divisors. We also prove the following theorem,
showing that in many cases the space $\overline{M}_{0,0}(X,\eta)$ is reducible.

\begin{customthm}{A}\label{thm_a}
  Let $X$ be a wonderful symmetric variety.
  Let $\eta \in N_1(X)\setminus \{0\}$ be a curve class such that $M_{0,0}^\circ(X, \eta) \neq \emptyset$.
  Assume there exists a decomposition $\eta = \eta_1 + \eta_2$, such that 
  $\eta_1$ and $\eta_2$ are nontrivial effective curve classes on $X$
  and a boundary divisor $D$ with
  \[
    \langle D, \eta_2 \rangle \leq -2.  
  \]

  Then the moduli space $\overline{M}_{0,0}(X, \eta)$ is reducible.
\end{customthm}

In particular we get:
\begin{customthm}{B}\label{thm_b}
  Let $X$ be a wonderful group compactification of a simple group $G$ of adjoint type with $\rk(G) \geq 3$.
  Then $\overline{M}_{0,0}(X, \eta)$ is reducible for any movable class of curves $\eta \in N_1(X)\setminus \{0\}$.
\end{customthm}
This result contrasts with the fact 
that the VMRT is irreducible for wonderful group compactifications.

This result motivates us to consider the set $M_{0,0}^\circ(X,\eta)$ and to prove that it is irreducible 
and can be considered a special component of the whole moduli space. 

\section{Irreducibility of the special component}

Let $G$ be a semisimple algebraic group of adjoint type.
We take $X$ to be the wonderful compactification of the group $G$, viewed as a symmetric space
$G\times G/G$ for the involution $(x,y) \mapsto (y,x)$.
Every fundamental co-weight $\omega_{i_0}$ for some $i_0 \in \{1, \dots, r\}$ defines a $1$-parameter subgroup
of $G$ and also of $G\times G$ with the inclusion $g\mapsto (g,1_G)$.
There exists a limit map $\varphi_{i_0}:X \dashrightarrow X$ with $\varphi_{i_0}(X) \subset X_{i_0}$. 
Let $\eta$ be a numerical equivalence class of mobile curves.
Then $M_{0,0}^\circ(X,\eta)$ is nonempty 
and contains an open set of curves, avoiding the indeterminacy locus of $\varphi_{i_0}$.
So there exists a rational map
\begin{align*}
  \varphi_{i_0}: M_{0,0}^\circ(X, \eta) &\dashrightarrow M_{0,0}^\circ(X_i, \bar{\eta})\\
    C &\mapsto \varphi_{i_0}(C),
\end{align*}
where $\bar{\eta}$ is the class of $\varphi_{i_0}(C)$.
A similar statement is true for closed $G$-stable subvarieties in $X$, 
which we use together with induction to prove the following:

\begin{customthm}{C}\label{thm_c}
  Let $X$ be the wonderful compactification of a semisimple algebraic group $G$ of adjoint type.
  Let $\eta \in N_1(X)$ be a movable curve class.
  Then the space $M_{0,0}^\circ(X,\eta)$ is irreducible.
\end{customthm}

The obvious continuation of this work would be to broaden the class of varieties, for which
the last irreducibility statement is true,
but we could also look towards rational simple connectedness of wonderful varieties. 
A variety $X$ is \textit{rationally connected}, if there exists an irreducible rational curve, 
passing through any two general points of $X$. The obvious 
example is the projective space $\mathbb{P}^N$. \textit{Rational simple connectedness}
requires more, namely that the space of rational curves itself is rationally connected.
Since there are many moduli spaces (for each class of curves) and since these moduli spaces
usually have many irreducible components as Theorem \ref{thm_a} shows,
a good definition of rational simple connectedness requires the existence of a 
"canonical" component of the space of stable maps of some class $\eta$ sufficiently positive which is rationally connected.
The component constructed in Theorem \ref{thm_b} is a good candidate for being this "canonical component".

\chapter{Spherical varieties}

\section{Definitions and preliminaries}\label{section_spherical_varieties}

\begin{definition}
  A normal $G$-variety $X$ is called \textit{spherical} if it there exists an $x\in X$ such that $\overline{Bx} = X$. 
  The stabilizer subgroup $H = \mathrm{Stab}_G(x)$ of the point $x$ is called a \textit{spherical} subgroup of $G$.
  The variety $X$ being spherical is equivalent to $Bx$ being open in $X$ or $BH$ being open in $G$.
  A spherical variety is called \textit{simple} if it contains a unique closed $G$-orbit. 
  The quotient $G/H$ by a spherical subgroup $H$ is called a \textit{spherical homogeneous space}.
\end{definition}

From now on we will assume that $X$ is a spherical variety, unless stated otherwise.

\begin{definition}
  \begin{enumerate}
    \item Define the set of $B$-semi-invariant rational functions on $X$:
    $$\mathbb{C}(X)^{(B)} = \{f \in \mathbb{C}(X)\setminus \{0\}\ |\ bf = 
    \chi(b)f\ \forall b \in B,\text{ where }\chi:B \to \mathbb{C}^*\};$$
    \item For any $f \in \mathbb{C}(X)^{(B)}$, $\chi$ is called the associated 
    $B$\textit{-weight}, it is a character of $B$ and is also denoted $\chi_f$;
    \item The set of all $B$-weights of nonzero rational functions on $X$ is denoted by:
    $$\Lambda(X) = \{\chi_f\ |\ f\in \mathbb{C}(X)^{(B)}\setminus \{0\} \};$$
    \item Define the $\mathbb{Q}$-vector space:
    $$\mathcal{Q}(X) = \mathrm{Hom}_\mathbb{Z}(\Lambda(X), \mathbb{Q})$$
  \end{enumerate}
  These definitions only depend on the open orbit of $X$, which is a spherical homogeneous space.
\end{definition}

\begin{theorem}[{\cite{sumihiro_ec}}]
  Let $X$ be a normal $G$-variety, $Y\subset X$ a $G$-orbit. Then there exists an open set $U$ of $X$ containing $Y$, 
  $G$-stable, and isomorphic to a $G$-stable subvariety of $\mathbb{P}(V)$ for some finite dimensional rational $G$-module $V$.
\end{theorem}

\begin{corollary}
  Any normal $G$-variety $X$ with a unique closed orbit is quasi-projective.
\end{corollary}

\begin{corollary}[{\cite[Corollary 2.2.2]{pezzini_losawv}}]
  Let $X$ be a normal $G$-variety, $Y\subset X$ an irreducible $G$-stable closed subvariety. 
  Then $B$-semi-invariant
  rational functions on $Y$ can be extended to $B$-semi-invariant rational 
  functions on $X$, as a consequence $\Lambda(X) \supset \Lambda(Y)$.
\end{corollary}

\begin{theorem}[{\cite[Corollary 2.1]{brion_vs}}]
  Let $X$ be a spherical $G$-variety. Then both $G$ and $B$ have finite number of orbits in $X$. If $Y$ is a 
  normal $G$-stable closed subvariety of $X$ then $Y$ is spherical.
\end{theorem}


\begin{definition}\label{def_spherical_b_open_subset}
  Let $X$ be a spherical variety.
  We denote by $\mathcal{D}_{X}$ the set of $B$-stable prime divisors.
  The $G$-stable prime divisors in $\mathcal{D}_{X}$ are called \textit{boundary divisors}
  and the $B$-stable but not $G$-stable prime divisors in $\mathcal{D}_{X}$ 
  are called \textit{colors} of $X$. We denote the set of all colors by $\Delta_{X}$.
  Let $Y$ be a $G$-orbit of $X$, define $\Delta_{X,Y} = \{D \in \Delta_{X}\ |\ Y \subset D\}$, 
  $\mathcal{D}_{X,Y} = \{D\in \mathcal{D}_{X}\ |\ Y \subset D\}$
  and $$X_{Y,B} = X \setminus \bigcup_{D\in \Delta_X\setminus\Delta_{X,Y}} D.$$
\end{definition}

Since $X$ has a dense $B$-orbit note that for $\chi \in \Lambda(X)$ there exists a unique up to a scalar
$f \in \mathbb{C}(X)$ with $\chi_f = \chi$. This implies that the following map is well defined.

\begin{definition}\label{def_spherical_var_colors_rho_map}
  Let $X$ be a simple spherical variety.
  Define the map 
  \begin{align*}
    \rho_X : \Delta_X &\to \mathcal{Q}(X)\\
    D &\mapsto \rho_{X}(D)  
  \end{align*} 
  via $\langle \rho_X(D), \chi \rangle = \nu_D(f)$,
  where $f \in \mathbb{C}(X)^{(B)}$ with $\chi_f = \chi$ and $\nu_D(f)$ is the vanishing multiplicity of $f$ on $D$.
\end{definition}

\begin{definition}
  Let $H$ be a spherical subgroup of $G$. A discrete valuation $\nu: \mathbb{C}(G/H) \to \mathbb{Q}\cup \{\infty\}$ is 
  called $G$-\textit{invariant}, if $\nu(g f) = \nu (f)$ for any $g \in G$ and $f \in \mathbb{C}(G/H)$. 
  For $G/H$ a spherical homogeneous space denote by $V_{G/H}$ the set of all $G$-invariant discrete valuations.
\end{definition}

\begin{proposition}[{\cite[Proposition 4.2]{brion_pauer_vehs}}]
  The valuation cone $V_{G/H}$ is a polyhedral convex cone.
\end{proposition}

\begin{definition}
  For a $G$-spherical variety $X$ with dense orbit $G/H$ and $Y \subset X$ a $G$-orbit, 
  define $V_{X,Y} = \{ \nu_D \in V_{G/H}\ |\ D \in \mathcal{D}_{X,Y}\setminus \Delta_{X,Y}\}$.  
\end{definition}

\section{Local structure of spherical varieties}

\begin{theorem}[{\cite[Proposition 2.2 and Theorem 2.3]{brion_vs}}]\label{thm_local_structure_spherical}
  Let $X$ be a spherical $G$-variety and $Y\subset X$ a closed $G$-orbit, then:
  \begin{enumerate}
  \item The set $X_{Y,B}$ is affine $B$-stable and is equal to $\{x \in X\ |\ \overline{Bx}\supset Y\}$;
  \item If $X$ is simple and $Y$ is the unique closed $G$-orbit then any divisor $D\in \Delta_X$ with $Y \not\subset D$
    is Cartier and globally generated;
  \item Let $P\subset G$ be the stabilizer of $X_{Y,B}$ and let $L\subset P$ be a Levi subgroup of $P$. There exists 
  an affine $L$-stable $L$-spherical subvariety $M$ of $X_{Y,B}$ such that the action 
  $$R_u(P) \times M \to X_{Y,B}$$
  is a $P$-equivariant isomorphism.	Moreover $\rk(X) = \rk(M)$.
  \end{enumerate}
\end{theorem}

\begin{corollary}[{\cite[Corollary 1.4.1]{brion_vs}}]
  A spherical variety has rank 0 if and only if it is a complete homogeneous space $G/P$ for some parabolic $P$.
\end{corollary}

\begin{corollary}[{\cite[Corollary 2.3.2]{pezzini_losawv}}]
  The variety $X$ in Theorem \ref{thm_local_structure_spherical} is smooth if and only if $M$ is smooth.
\end{corollary}


\section{Classification of spherical embeddings}

In this section we present a classification of simple spherical varieties via 
combinatorial data called colored cones.

\begin{definition}
  An \textit{embedding} of a spherical homogeneous space $G/H$ is a normal 
  $G$-variety equipped with a point $x \in X$, such that $G x$ is open and $\mathrm{Stab}_G(x) = H$.
\end{definition}


\begin{definition}
  Let $V$ be a $\mathbb{Q}$-vector space and $C \subset V$.
  \begin{enumerate}
    \item $C$ is a \textit{cone} if it is a semi-group with respect to addition in $V$ and stable under
      multiplication with $\mathbb{Q}_{\geq 0}$;
    \item A cone $C$ is \textit{strictly convex} if $C\cap -C = 0$ 
      (equivalently if $C$ contains no line);
    \item A cone $C$ is \textit{polyhedral} if $C=\mathbb{Q}_{\geq 0}v_1+\dots+\mathbb{Q}_{\geq 0}v_n$
      for some $v_1, \dots, v_n$;
    \item $C^\vee = \{\nu \in V^\vee\ |\ \langle \nu, v \rangle \geq 0 \text{ for any } v \in C\}$ is 
      the \textit{dual} cone of $C$;
    \item $C^\circ$ is the relative interior of $C$;
    \item A \textit{face} $F_f$ of a cone $C$ is a subset $F_f = \{v\in C\ |\ f(v)= 0 \}$ for some $f\in C^\vee$.
  \end{enumerate}
\end{definition}

\begin{definition}
  Let $X$ be a simple spherical $G$-variety and $Y\subset X$ a closed $G$-orbit. Define
  the cone $\mathcal{C}_{Y}(X)\subset \mathcal{Q}(X)$ to be the cone 
  generated by $V_{X,Y}$ and $\rho_X(\Delta_{X,Y})$.
\end{definition}




\begin{definition}
  A \textit{colored cone} of $G/H$, is a pair $(\mathcal{C}, \mathcal{F})$
    with $\mathcal{C}\subset \mathcal{Q}(G/H)$ and $\mathcal{F}\subset \mathcal{D}_{G/H}$
      having the following properties:
  \begin{enumerate}
    \item $\mathcal{C}$ is a cone generated by $\rho_{G/H}(\mathcal{F})$ 
      and finitely many elements of $V_{G/H}$
    \item The intersection $\mathcal{C}^\circ\cap V_{G/H}$ is nonempty.
  \end{enumerate}
  The colored cone $\mathcal{C}$ is called \textit{strictly convex} if the cone 
  $\mathcal{C}$ is strictly convex and $0 \notin \mathcal{F}$.
\end{definition}

\begin{proposition}[{\cite[\S 8]{luna_vust_peh}}]
  The map $X\mapsto (\mathcal{C}_Y(X), \Delta_{X,Y})$ is a bijection between the isomorphism classes
  of simple spherical embeddings $X$ of $G/H$ with closed orbit $Y$ and strictly convex colored cones.
  The cone $V_{G/H}$ is a polyhedral cone.

  Furthermore the variety $X$ is complete if $V_{G/H} \subset \mathcal{C}_Y(X)$.
\end{proposition}

\begin{definition}
  A spherical variety $X$ is called \textit{toroidal} if $\Delta_{X,Y} = \emptyset$ 
  for any closed $G$-orbit $Y\subset X$. A spherical embedding of $G/H$
  is called toroidal if it is a toroidal variety.
\end{definition}

\begin{proposition}[{\cite[Corollaire 3.2]{brion_pauer_vehs}}]
  Every toroidal spherical embedding of $G/H$ has an equivariant toroidal completion.
\end{proposition}

\begin{definition}
  Let $ (\mathcal{C}, \mathcal{F})$ be a colored cone. A colored cone 
    $ (\mathcal{C}', \mathcal{F}')$ is called a \textit{colored face} of $ (\mathcal{C}, \mathcal{F})$
      if the following conditions are satisfied:
      \begin{enumerate}
        \item $\mathcal{C}'$ is a face of $\mathcal{C}$
        \item $\mathcal{F}' = \mathcal{F}\cap \rho_X^{-1}(\mathcal{C}')$
      \end{enumerate}
\end{definition}

The closed $G$-stable subvarieties of a spherical variety can be described as follows:

\begin{proposition}[{\cite[Lemma 3.2]{knop_lvtse}}]
  Let $X$ be an embedding of a spherical homogeneous space $G/H$ and let $Y\subset X$ a $G$-orbit.
  There is a bijection between the set of $G$-orbits in $X$
  containing $Y$ in their closure and the set of faces of $(\mathcal{C}_Y(X),\Delta_{X,Y})$.
\end{proposition}

\begin{definition}
  A spherical homogeneous space $G/H$ is called \textit{sober} if it admits a simple
  projective spherical embedding.
\end{definition}

\begin{proposition}[{\cite[Corollary 5.3]{brion_pauer_vehs}}]\label{prop_sober_characterization}
  The following conditions are equivalent:
  \begin{itemize}
    \item $G/H$ is sober;
    \item The cone $V_{G/H}$ is strictly convex;
    \item $N_G(H)/H$ is finite.
  \end{itemize}
\end{proposition}

Recall that simple varieties are quasi-projective.

\begin{proposition}\label{proposition_unique_toroidal_embedding}
  There exists a simple projective toroidal embedding of a sober spherical homogeneous space 
  and it is unique up to isomorphism.
\end{proposition}

\begin{proof}
  The embedding is given by the colored cone $(V_{G/H}, \emptyset)$.
\end{proof}

\section{Curves and divisors on spherical varieties}\label{sec_spherical_curves_divisors}

The following is a well known result on projective rational homogeneous spaces.

\begin{proposition}\label{prop_picard_group_projective_homogeneous_space}
  Let $G$ be a reductive group, $S' \subset S_G$ a subset of simple roots. Let $P = P_{S'}$ be 
  the associated parabolic subgroup of $G$.
  Then 
  \[
    \Pic(G/P^-) = \bigoplus_{\alpha \in S_G\setminus S_P}\mathbb{Z} D_\alpha(G/P^-),
  \]
  where $D_\alpha(G/P^-) = \overline{P s_\alpha P^-/P^-}$ and $s_\alpha$ is a representative of a simple reflection 
  corresponding to $\alpha$ in the Weyl group of $G$. The divisors $D_\alpha(G/P^-)$ are the colors of $G/P^-$.
\end{proposition}

There is another way to describe the Picard group of a homogeneous space:

\begin{proposition}
  Let $\mathcal{L} \in \Pic(G/P^-)$ be a line bundle, $z \in G/P^-$ be the $P^-$-fixed point
  and $\mathcal{L}_z$ be the fiber of $\mathcal{L}$ at $z$.
  Let $\tilde{G}$ be a finite covering of $G$, such that $\mathcal{L}$ is $\tilde{G}$-linearisable
  and $\tilde{P}^- \subset \tilde{G}$ the parabolic subgroup obtained as inverse image of $P^-$.
  Then $\mathcal{L} \simeq \tilde{G}\times^{\tilde{P}^-}\mathcal{L}_z$.
  The $1$-dimensional representation $\mathbb{C}_\chi = \mathcal{L}_z$ of $\tilde{P}^-$, 
  which is given by a character $\chi \in \Lambda(\tilde{G}/\tilde{P}^-)$. We will write $\mathcal{L} = \mathcal{L}_\chi$.
  
  Moreover $\mathcal{O}_{G/P^-}(D_\alpha(G/P^-)) = \mathcal{L}_{\omega_\alpha}$, where
  $\omega_{\alpha}$ is the fundamental weight, corresponding to the simple root $\alpha \in S$.
\end{proposition}

\begin{notation}
  Let $\mathcal{L}$ be a line bundle on $G/B^-$.
  Replace $G$ with a finite covering, such that $\mathcal{L}$ is $G$-linearisable.
  Then it is isomorphic to $G\times^{B^-} \mathbb{C}_\lambda$ 
  for some $B^-$-weight $\lambda$. We will write $\mathcal{L} = \mathcal{L}_\lambda$.

  Let $X$ be a complete spherical variety and let $D$ be a color of $X$. 
  Then we have a line bundle $\mathcal{O}_X(D)$ on $X$. Let $Y$ be a closed $G$-orbit of $X$.
  Let $\iota: Y \to X$ be the inclusion of the closed orbit $Y$ in $X$ and $\pi: G/B^- \to Y$ be the usual map of flag varieties.
  Then $(\iota \circ \pi)^*\mathcal{O}_X(D) = \mathcal{L}_\lambda(G/B^-)$ for some $B^-$-weight $\lambda$.

  If $X$ is simple, then there is exactly one closed $G$-orbit $Y \subset X$.
  If $(\iota \circ \pi)^*\mathcal{O}_X(D) = \mathcal{L}_\lambda(G/B^-)$
  we will write $\mathcal{L}_\lambda = \mathcal{L}_\lambda(G/B^-)$ and
  we will say that $D$ \textit{corresponds} to the $B^-$-weight $\lambda$.
\end{notation}


\begin{notation}
  Denote by $A_k(X)$ the group of $k$-cycles on $X$ modulo rational equivalence.
  Denote by $Eff_k(X) \subset A_k(X)_{\mathbb{Q}}$ the cone generated by all effective $k$-cycles on $X$.
  We will write $Eff^k(X) = Eff_{\dim X - k}(X)$ and $Eff(X)$ for the cone of effective divisors on $X$.
\end{notation}

\begin{proposition}[{\cite[Corollary to Theorem 1]{fulton_macpherson_sottile_sturmfels_itosv}}]\label{cycles_spherical}
  For $X$ spherical the cone of effective cyc\-les $Eff_k(X) \subset {A_k(X)_{\mathbb{Q}}}$ is a polyhedral cone, 
  generated by the classes of the closures of the $B$-orbits.
\end{proposition}
Let $Div(X)$ be the group of Cartier divisors on $X$ and let $Z_1(X)$ be the
group of $1$-cycles on $X$. 

\begin{definition}
  A Cartier divisor is called \textit{effective} if it is a nonnegative linear
  combination of prime divisors. A $1$-cycle is called \textit{effective} if
  it is a nonnegative linear combination of irreducible curves.
\end{definition}
There exists an intersection pairing on
$\langle , \rangle : Div(X) \times Z_1(X)\to \mathbb{Z}$ given by $deg(\restr{\mathcal{O}_{X}(D)}{C})$,
for $\mathcal{O}_{X}(D) \in Div(X)$ and $C$ an irreducible $1$-cycle.

\begin{definition}
  Define $N^1(X)$ to be the quotient of $Div(X)$ by the orthogonal of $Z_1(X)$
  and $N_1(X)$ to be the quotient of $Z_1(X)$ by the orthogonal of $Div(X)$.
\end{definition}
Then we have a perfect pariring between $N_1(X)$ and $N^1(X)$.

\begin{proposition}[{\cite[Proposition 3.2]{brion_candisv}}]
  Let $X$ be a complete spherical variety. The map $\Pic(X) \to N^1(X)$
  is an isomorphism.
\end{proposition}

\begin{definition}
  A Cartier divisor $D$ on a projective variety $X$ is called \textit{nef} if
  $\langle D, C \rangle \geq 0$ for any effective $1$-cycle $C \in Z_1(X)$.
  The class $[D] \in \Pic(X)$ is also called nef.
  We write $Nef(X) \subset N^1(X)_{\mathbb{Q}} = \Pic(X)_{\mathbb{Q}}$ for the cone of nef divisor classes.
\end{definition}


\begin{proposition}[{\cite[Theorem 3.2.1]{kleiman_tnta}}]
  Let $X$ be a projective variety. The nef cone of $X$ is the closure of the cone of 
  ample divisors of $X$.
\end{proposition}

Michel Brion gives a description of the ample divisors of a simple spherical variety in \cite{brion_gpncvs}, which we 
can use to describe $Nef(X)$.

\begin{proposition}[{\cite[Theorem 2.6]{brion_gpncvs}}]\label{simple_spherical_ample_cone}
  Let $X$ be a simple projective spherical variety.
  A Cartier divisor on $X$ is ample is and only if it linearly equivalent to a positive 
  linear combination of divisors in $\Delta_X\setminus \Delta_{X,Y}$.
\end{proposition}

\begin{corollary}\label{prop_picard_wonderful}
  Let $X$ be a simple projective spherical variety.
  The Picard group $\Pic(X)$ is generated by the colors of $X$, which do not contain
  the closed $G$-orbit.
\end{corollary}

\begin{corollary}\label{nef_cone}
  Let $X$ be a simple spherical projective variety. 
  The nef cone is given by 
  \[
    Nef(X) = \bigoplus_{D \in \Delta_X\setminus \Delta_{X,Y}}\mathbb{N}_0 D
  \]
\end{corollary}

In the same paper Michel Brion shows that a divisor on a simple projective spherical variety is generated by its sections 
if and only if it is a nonnegative linear combination of colors of $X$.

\begin{proposition}[{\cite[Proposition 2.6]{brion_gpncvs}}]\label{simple_spherical_generated_by_sections}
  Let $X$ be a simple projective spherical variety. A divisor on $X$ is generated by its sections 
  if and only if it is linearly equivalent to 
  \[
    \sum_{D \in \Delta_X\setminus \Delta_{X,Y}} n_D D,
  \]
  with all $n_D \geq 0$.
\end{proposition}


If the spherical variety is complete, then the space of curve classes is dual to those 
of divisor classes:

\begin{proposition}[{\cite[Prop 3.1 and 3.2]{brion_candisv}}]\label{prop_toroidal_curve_divisor_duality}
  Let $X$ be a complete spherical variety, the map $\Pic(X) \to \Hom(A_1(X), \mathbb{Z})$
  is an isomorphism.
  Moreover the sequence
  \[
    \begin{tikzcd}
      0 \arrow[r] & A_1(X)_{tors} \arrow[r] & A_1(X) \arrow[r] & N_1(X) \arrow[r] & 0
    \end{tikzcd}    
  \]
  is exact. Moreover if $X$ contains a unique closed $G$-orbit, then
  $A_1(X)$ is torsion free and $N_1(X) = A_1(X)$.
\end{proposition}

We defined the nef cone for divisors.
The nef cone is dual to the cone of effective curve classes.
We define a similar cone of curve classes
$Nef_1(X) \subset N_1(X)_{\mathbb{Q}}$ dual to the cone of effective
divisors. 

We now define the movable curve classes.

\begin{definition}
  A curve class $\eta \in N_1(X)_{\mathbb{R}}$ is called \textit{movable} if for any point $x \in X$ there 
  exists a representative $C_x$ of $\eta$, such that $C_x$ is an effective curve passing through $x$.
  Denote by $Mov(X)$ the closure of the cone of movable effective classes of curves.
\end{definition}

It turns out that for smooth projective varieties $\overline{Nef_1(X)} = Mov(X)$:

\begin{proposition}[{\cite[Theorem 2.2]{boucksom_demailly_paun_peternell_pecckmvnkd}}]\label{movable_effective_cones_duality}
  Let $X$ be a smooth projective variety,
  then the cone $Mov(X)$ is dual to the closure of the cone generated by effective divisor classes
  in $N^1(X)_\mathbb{R}$.
\end{proposition}



We finish the section by giving a description of the canonical divisor of a spherical variety:

\begin{proposition}[{\cite[Theorem 6.1.1]{perrin_sl}}]
  Let $X$ be a spherical variety, then the following holds for the anticanonical divisor $-K_X$ of $X$
  \[
    -K_X =  \sum_{D \in \mathcal{D}_{X}\setminus \Delta_X} D  + \sum_{D \in \Delta_X} a_D D,
  \]
  where $a_D$ is a nonnegative integer. 
\end{proposition}

\chapter{Wonderful varieties}

In this chapter we define wonderful varieties. Throughout the rest of the chapter we assume $X$ to be 
a wonderful $G$-variety, unless stated otherwise.

\section{Definitions and preliminaries}

\begin{definition}
  Let $X$ be a $G$ variety. Then $X$ is \textit{wonderful} if:
  \begin{enumerate}
    \item $X$ is smooth and complete
    \item $X$ contains an open $G$-orbit $X^\circ_G$ whose complementary is the union of $G$-stable prime divisors 
    $X_1,\dots, X_l$, which have simple normal crossings
    \item for all $x,y \in X$ we have:
    \[
      Gx = Gy \textrm{ equivalent to }\{i| x\in X_i\} = \{j| y\in X_j\}
    \]
  \end{enumerate}
  The number $l$ is called the \textit{rank} $\rk(X)$ of $X$, 
  and the union of $G$-stable prime divisors is called the 
  \textit{boundary} of $X$, denoted $\partial X$. The wonderful variety $X$ is called \textit{exceptional} if
  $\rk\Pic(X) > \rk(X)$.
\end{definition}

\begin{proposition}
  Any closed $G$-stable subvariety of $X$ is of the form:
  $$
  X_I = \bigcap_{i\in I} X_i
  $$
  for some $I \subset \{1,\dots,r\}$. Moreover, $X_I$ is wonderful and $\rk(X_I) = \rk(X) - |I|$.
\end{proposition}

Wonderful varieties are a special case of spherical varieties, as proved by Luna in \cite{luna_tvmes}.

\begin{theorem}[{\cite{luna_tvmes}}]
  Any wonderful $G$-variety is a simple spherical variety.
\end{theorem}

\begin{definition}
  Let $G/H$ be a spherical homogeneous space and $X$ a wonderful variety,
  such that its open $G$-orbit $X^\circ_G$ is isomorphic to $G/H$. 
  Then $X$ is called a \textit{wonderful embedding of $G/H$}.
\end{definition}

\begin{remark}
  Note that if $G/H$ admits a wonderful embedding, then $G/H$ is sober.
\end{remark}

\begin{proposition}[{\cite[Proposition 3.3.1]{pezzini_losawv}}]
  Let $G/H$ be a sober spherical homogeneous space.
  A spherical embedding of $G/H$ is wonderful if and only if 
  it is the simple projective toroidal embedding of $G/H$ (see Proposition
  \ref{proposition_unique_toroidal_embedding}) and it is smooth. 
  The ranks as wonderful variety and as spherical variety coincide.
\end{proposition}

As a consequence we get a characterization of wonderful varieties in the category of spherical varieties.

\begin{corollary}
  A spherical variety $X$ is wonderful if and only if it smooth simple complete and toroidal.
\end{corollary}

Domingo Luna introduced spherical systems in \cite{luna_vsdta}
and conjectured that they classify wonderful varieties. In the same paper he proved
this conjecture in type $A$. The conjecture was completed
in 2015 by Paolo Bravi and Guido Pezzini, see \cite{bravi_pezzini_pwv}.
This classification is used to classify all spherical homogeneous spaces by combinatorial data. 
We introduce spherical systems in the next section.

\section{Spherical systems}

\begin{definition}\label{def_luna_spherical_roots_and_system}
  Let $H$ be a spherical subgroup of $G$, let $V_{G/H}$ be the cone of $G$-invariant discrete valuations 
  on the spherical homogeneous space $G/H$. We define the set $\Sigma_{G/H}$ of \textit{spherical roots}
  as the unique set verifying:
  \begin{itemize}
    \item $\Sigma_{G/H}$ consists of linearly independent primitive elements of $\Lambda(G/H)$
    \item $V_{G/H} = \{v \in \mathcal{Q}(G/H)\ |\ \langle v, \gamma\rangle \leq 0, \gamma\in \Sigma_{G/H}\}$
  \end{itemize}
  Since the spherical roots lie in a half-plane, we can find a Borel subgroup $B$ of $G$, such that
  any spherical root is a nonnegative linear combination of simple roots of $G$.
  Define the \textit{support of a spherical root $\gamma$} denoted by $\textrm{Supp}(\gamma)$ 
  as the set of simple roots of $G$, such that the 
  coefficients of $\gamma$ are nonzero. The support of a subset of $\Sigma_{G/H}$ is the union 
  of supports of its elements.
\end{definition}

\begin{proposition}[{\cite[Section 1.2]{luna_vsdta}}]\label{proposition_t_weights_spherical_roots}
  Let $X$ be a wonderful $G$-variety with $Y \subset X$ its unique closed orbit. 
  Let $z \in Y \subset X$ be the unique fixed point of $B^-$.
  Then the set of weights of the $T$-action on the normal space
  \[
    \frac{T_zX}{T_z(G z)}
  \]
  coincides with $\Sigma_X$ seen as elements of $\Lambda(X)$.
\end{proposition}

\begin{definition}\label{def_roots_move_colors}
  Let $X$ be a wonderful $G$-variety. Let $D\in \Delta_X$ be a color of $X$. 
  A simple root $\alpha \in S_G$ \textit{moves} $D$ if
  $D$ is not stable under the action of the minimal parabolic $P_{\alpha}$. 
  Denote by $\Delta_X(\alpha)$ the set of colors moved by $\alpha$.
  With that we define the following objects:
  \begin{itemize}
    \item $S_X^p\subset S$ is the set of simple roots of $G$ moving no color;
    \item $A_{G,X}$ the set of colors which are moved by simple roots $\alpha \in G$ that
      are also spherical roots, i.e. $\alpha \in S \cap \Sigma_X$.
  \end{itemize}
  The triple $(S_X^p, \Sigma_X, A_{G,X})$ is called the \textit{spherical $G$-system} of $X$.
\end{definition}

\begin{remark}
  The notion of spherical $G$-system can be defined in a purely combinatorial way, not depending on
  $X$, see Definition 1.2.2 in \cite{bravi_pezzini_pwv}.
\end{remark}

\begin{proposition}[{\cite[Section 2.7]{luna_gcplvs}}]\label{prop_color_types}
  Let $X$ be a wonderful $G$-variety with open $G$-orbit $G/H$. Then for any simple root of $G$ 
  exactly one of the four cases is true:
  \begin{itemize}
    \item[$(p)$] $\Delta_X(\alpha) = \emptyset$, if $\alpha \in S_X^p$;
    \item[$(a)$] $\Delta_X(\alpha)$ contains two elements and $\alpha \in \Sigma_X$;
    \item[$(a')$] $\Delta_X(\alpha)$ contains exactly one element and $2\alpha \in \Sigma_X$;
    \item[$(b)$] $\Delta_X(\alpha)$ contains exactly one element, but $2\alpha \notin \Sigma_X$;
  \end{itemize}
  We will write that $\alpha \in S_X^p$ (resp. $S_X^a$, $S_X^{a'}$ or $S_X^b$) if the case (p) (resp. (a), (a'), (b)) holds.
  The union of sets 
  \[
    \bigcup_{\alpha \in S^a} \Delta_X(\alpha) \cup \bigcup_{\alpha \in S^{a'}} \Delta_X(\alpha) \cup \bigcup_{\alpha \in S^b} \Delta_X(\alpha)
  \]
  is disjoint.
\end{proposition}


\section{Classification of wonderful varieties}

Now that we have defined the spherical system we can give the Theorem in \cite{bravi_pezzini_pwv}
classifying the wonderful varieties.

\begin{theorem}[{\cite[Theorem 1.2.3]{bravi_pezzini_pwv}}]\label{thm_wond_var_classification}
  There is a bijection between isomorphism classes of wonderful $G$-varieties and 
  spherical $G$-systems. 
\end{theorem}

\begin{lemma}\label{lemma_closed_orbit_form}
  Let $X$ be a wonderful $G$-variety. Then the closed orbit $Y$ of $G$ is isomorphic to $G/P_{S_X^p}^-$.
\end{lemma}

\begin{proof}
  The closed orbit $Y$ is projective homogeneous variety, so it is of the form $G/P^-$
  for some parabolic subgroup $P \subset G$. As any wonderful variety is toroidal, we have that
  $X_{Y,B} \cap Y$ is nonempty and is the open $B$-orbit in $Y = G/P^-$
  by Theorem \ref{thm_local_structure_spherical}.
  We therefore have $P = \mathrm{Stab}(X_{Y,B}\cap Y)$.
  We want to show that $P_{S_X^p} = P$.
  
  We have that $P_{S_X^p}$ stabilizes $X_{Y,B}$, thus $P_{S_X^p}$
  stabilizes $X_{Y,B} \cap Y$ and $P_{S_X^p} \subset P$.
  
  Let $P_{\alpha}$ be a minimal parabolic. Suppose there exists a color $D \in \Delta_X$, 
  such that $P_\alpha D = X$. Then $P_\alpha (D\cap Y) = P_\alpha D \cap Y = Y$, 
  so $P_\alpha \not \subset P$.
  We have $P \subset \langle P_\alpha \ |\ \Delta_X(\alpha) = \emptyset\rangle$,
  i.e. $P \subset P_{S_X^p}$.
  $\qed$
\end{proof}

\begin{proposition}[{\cite[page 3]{pezzini_aowv}}]
  Let $X$ be a wonderful variety.
  Let $z$ be the unique fixed point of $B^-$ in the closed orbit $Y\subset X$.
  For any boundary divisor $X_i$ there exists a unique spherical root $\gamma_i$,
  such that the $T$-weight of the quotient $T_z X/T_z X_i$ is precisely $\{\gamma_i\}$. 
\end{proposition}

\begin{proposition}
  Let $X_I \subset X$ be a closed $G$-stable subvariety for some $I \subset \{1, \dots, \rk(X)\}$. 
  Then the spherical system of ${X_I}$ is given by $(S_{X_I}^p, \Sigma_{X_I}, A_{G, {X_I}})$, where
  \begin{itemize}
    \item $S_{X_I}^p = S_X^p$
    \item $\Sigma_{X_I} = \{\gamma_i\in \Sigma_X\ |\ i \notin I\}$
    \item $A_{G,{X_I}} = \bigcup_{\alpha \in S_G \cap \Sigma_{X_I}}\Delta(\alpha)$
  \end{itemize}
\end{proposition}

\begin{proof}
  As seen in Lemma \ref{lemma_closed_orbit_form}, $S_X^p$ is also defined by the closed orbit $Y$ of $X$.
  Since $Y \subset X_I \subset X$, we have that $S_{X_I}^p = S_X^p$.
  Next we consider $\Sigma_{X_I}$. We have that $X_i \cap X_I$ is a boundary divisor of $X_I$ for every $i \notin I$ 
  and the corresponding spherical root is the $T$-weight of the quotient 
  $T_z X/T_z X_i = T_z X_I/T_z (D_i\cap X_I)$, because $\partial X$ is a smooth normal crossings divisor.
  So $\Sigma_{X_I} = \{\gamma_i \in \Sigma_X\ |\ i \notin I\}$.
  Lastly $A_{G, {X_I}} = \bigcup_{\alpha \in S_G \cap \Sigma_{X_I}}\Delta(\alpha)$ by definition of $A_{G,X_I}$.
  $\qed$
\end{proof}

\section{Curves and divisors on wonderful varieties}

Wonderful varieties are spherical, so we can use everything
already mentioned in Section \ref{sec_spherical_curves_divisors}.
In particular, the colors form a basis of $\Pic(X)$ and we now consider the dual basis
of $N_1(X)$.

\begin{lemma}[{\cite[Lemma 3.2.2]{luna_gcplvs}}]\label{lemma_wond_color_dual_curve_description}
  For any divisor $D \in \Delta_X$ there exists a curve $C_D$ on $X$, such that 
  for any $D' \in \Delta_X$:
  \[
    \langle D', [C_D] \rangle = 
    \begin{cases}
      1 &D = D',\\
      0 &D\neq D'.
    \end{cases}  
  \]
\end{lemma}

Let $Y$ be as usual the closed orbit of $X$, we have seen that $Y = G/P^-$.
Let $\alpha \in S\setminus S_X^p$. The Schubert curve in $Y$, associated to $\alpha$ is
$C_\alpha = \overline{B^- s_\alpha P^-}/P^-$. 
The classes $([C_\alpha])_{\alpha \in S\setminus S_X^p}$ form a basis of $N_1(Y)$,
which is dual to the basis $(D_\alpha(G/P^-))_{\alpha \in S\setminus S_X^p}$ of $\Pic(Y)$.
Denote by $\iota: Y \to X$ the inclusion and consider $\iota_*[C_\alpha]$.

\begin{lemma}[{\cite[Lemma 3.2.1]{luna_gcplvs}}]\label{lemma_wond_color_closed_orbit_curves}
  Let $\alpha \in S\setminus S_X^p$ be a simple root of $G$. We have the following for $\iota_*[C_\alpha]$:
  \begin{itemize}
    \item if $\alpha \in S^a_X$, we have $\iota_*[C_\alpha] = [C_{D^+}] + [C_{D^-}]$, where $\Delta_X(\alpha) = \{D^+, D^-\}$;
    \item if $\alpha \in S^{a'}_X$, we have $\iota_*[C_\alpha] = 2[C_D]$, where $\Delta_X(\alpha) = \{D\}$;
    \item if $\alpha \in S^{b}_X$, we have $\iota_*[C_\alpha] = [C_D]$, where $\Delta_X(\alpha) = \{D\}$.
  \end{itemize}
\end{lemma}


Wonderful varieties are complete and spherical, so the colors generate the Picard group of $X$,
and there is an explicit way to write the boundary divisors as linear combinations 
of colors:

\begin{proposition}[{\cite[Proposition 3.1]{luna_gcplvs}}]\label{prop_wond_var_boundary_divisors_as_colors}
  Let $X$ be a wonderful $G$-variety for a connected reductive group $G$. 
  Let $X_i$ be its boundary divisors with corresponding spherical roots
  $\gamma_i$ and let $\Delta_X$ be the set of colors of $X$. Then 
  \[
    X_i = \sum_{D \in \Delta_X}\langle \rho_X(D), \gamma_i \rangle D,  
  \] 
  where $\rho_X: \Delta_X \to \mathcal{Q}$ is defined in \ref{def_spherical_var_colors_rho_map}.
\end{proposition}

\begin{corollary}\label{corollary_wond_var_boundary_divisors_intersection_with_curves}
  Let $X$ be a wonderful $G$-variety, $X_i$ its boundary divisors
  with spherical roots $\gamma_i$ and $\Delta_X$ its set of colors. 
  Let $\eta \in N_1(X)$ be the class of a curve on $X$ with $\eta = \sum_{D\in \Delta_X} c_D [C_D]$.
  Then 
  \[
    \langle X_i, \eta\rangle = \sum_{D\in \Delta_X} c_D \langle\rho_X (D), \gamma_i\rangle
  \]
\end{corollary}

\begin{proposition}[{\cite[Proposition 3.4]{luna_gcplvs}}]\label{prop_wond_var_rho_value}
  Let $X$ be a wonderful $G$-variety, $\gamma$ be a spherical root. 
  Let $\alpha$ be a simple root of $G$, then
  \begin{itemize}
    \item if $\Delta_X(\alpha) = \{D\}$ and $2\alpha\in \Sigma_X$, then 
      $\langle \rho_X(D), \gamma\rangle = \frac{1}{2}\alpha^\vee(\gamma)$;
    \item if $\Delta_X(\alpha) = \{D^+, D^-\}$, then 
      $\langle \rho_X(D^+) + \rho_X(D^-), \gamma \rangle = \alpha^\vee(\gamma)$;
    \item if $\Delta_X(\alpha) = \{D\}$ and $2\alpha \notin \Sigma_X$, $\langle \rho_X(D), \gamma\rangle = \alpha^\vee(\gamma)$.
  \end{itemize}
\end{proposition}

\section{Luna cells and limit map}\label{sec_limitstuff}

In this section we assume $X$ to be smooth simple projective variety, unless stated otherwise
and $Y \subset X$ to be the unique closed orbit of $X$.

\begin{definition}\label{def_limit_stuff}
  Fix a subset $S' \subset S$ of the set of simple roots. We introduce the following notation: 
  \begin{enumerate}
    \item $G^{S'} = P_{S'}\cap P_{S'}^-$, the Levi subgroup of $P_{S'}$, containing $T$;
    \item $X_{S'} = \{x \in X \ |\ \overline{P_{S'}x}\supset Y\}$, an open subset in $X$;
    \item $X^{S'} = \{x \in X\ |\ Z(G^{S'})^\circ x = x\}$, 
      the fixed points under the action of the connected centre of $G^{S'}$;
    \item $\Delta_X(S') = \cup_{\alpha\in S'}\Delta_X(\alpha)$, 
      the set of colors not stabilized by the parabolic $P_{S'}$.
  \end{enumerate}
\end{definition}

\begin{definition}
  Let $\lambda: \mathbb{C}^* \to T$ be a $1$-parameter subgroup of $G$, 
  write $\lambda \in \Lambda(T)^*$. 
  We will call $\lambda$ \textit{adapted} to $S'$ if:
  \begin{equation*}
  \langle \lambda,\alpha\rangle 
    \begin{cases}
    =0 & \text{if }\alpha \in S',\\
    > 0 & \text{if }\alpha \in S\setminus S'.
    \end{cases}
  \end{equation*}
\end{definition}\

\begin{definition}
  Let $\lambda$ be a $1$-parameter subgroup of $G$, we define the following:
  \begin{enumerate}
    \item $X_Z = \{x \in X\ |\ \lim_{t \to 0}\lambda(t) x \in Z\}$ for $Z \subset X$ 
      a connected component of the smooth fixed points set $X^{\lambda(\mathbb{C}^*)}$;
    \item $X^\lambda$ the connected component of $X^{\lambda(\mathbb{C}^*)}$ 
      such that $X_{X^\lambda}=: X_\lambda$ is open in $X$;
    \item $G^\lambda = C_G(\lambda(\mathbb{C^*}))$ a reductive subgroup of $G$;
    \item $G_\lambda = \{g \in G\ |\ \lim_{t \to 0} \lambda(t)g\lambda(t)^{-1} \in G^\lambda\}$ a parabolic subgroup of $G$,
      whose Levi is $G^\lambda$.
    \item $\pi_\lambda: X_\lambda \to X^\lambda$ defined by 
      $\pi_\lambda(x)= \lim_{t\to 0}\lambda(t)x$ and by abuse of notation
      $\pi_\lambda: G_\lambda \to G^\lambda$ defined by $\pi_\lambda(g)= \lim_{t\to 0}\lambda(t)g\lambda(t)^{-1}$.
      Clearly $\pi_\lambda(gx) = \pi_\lambda(g)\pi_\lambda(x)$.
  \end{enumerate}
\end{definition}

\begin{proposition}[{\cite[Propositions 1.3 and 2.6]{luna_gcplvs}}]\label{prop_luna_limit_stuff}
  Let $S'\subset S$. Then
  \begin{enumerate}
    \item for every $\lambda \in \Lambda(T)^*$ adapted to $S'$ we have $X_\lambda = X_{S'}$ and $X^\lambda = X^{S'}$;
    \item $R_u(P_{S'})$ fixes $X^{S'}$;
    \item the $P_{S'}$-variety $X^{S'}$ is simple spherical;
    \item if $X$ is a wonderful $G$-variety, then $X^{S'}$ is a wonderful $G^{S'}$-variety.
  \end{enumerate}
  For any $D \in \Delta_X(S')$ we have $D \cap X^{S'} \in \Delta_{X^{S'}}$.
  Moreover the assigment
  \begin{align*}
    \Delta_X(S') &\to \Delta_{X^{S'}}\\
    D &\mapsto D\cap X'
  \end{align*}
  is bijective and its inverse is given by $D' \mapsto \overline{\pi_\lambda^{-1}(D')} = \overline{P_{S'}D}$.
  
  Furthermore the set of spherical roots of $X^\lambda$ is given by
  \[
    \Sigma_{X^\lambda} = \Sigma_X \cap \Lambda_{G^{S'}},
  \]
  where $\Lambda_{G^{S'}}$ is the lattice, generated by the roots of $G^{S'}$.
\end{proposition}

\begin{lemma}\label{lemma_parabolic_induction_and_limit_stuff}
  Let $S' \subset S$ be a subset of simple roots of $G$ and let $P = P_{S'}$ 
  and let $L = G^{S'}$ be its Levi subgroup.
  If $X'$ is a wonderful $L$-variety, then $X = G \times^{P^-} X'$, where $P^-$ acts via its
  Levi subgroup $L$, is a wonderful $G$-variety and $X^{S'} = X'$.
\end{lemma}

\begin{proof}
  The projection $G\times^{P^-} X' \to G/P^-$ is a locally trivial fibration.
  Any $G$-orbit $O \subset X$ is of the form $G O'$ for a unique 
  $L$-orbit $O' \subset X'$ and the assigment $O \mapsto O'$
  preserves codimensions.
  It follows that $X$ is wonderful.

  Let $\lambda$ be a $1$-parameter subgroup, adapted to $S'$.
  Let $L/H$ be the 
  open $L$-orbit of $X'$ for some spherical subgroup $H \subset L$.
  The variety $X'$ is wonderful, so in particular $X'$ is simple.
  By Proposition \ref{prop_sober_characterization} we have that $N_L(H)/ H$ is finite.
  The $1$-parameter subgroup $\lambda(\mathbb{C}^*)$ lies in $Z(L)^\circ \subset N_L(H)$, as
  $\lambda$ is adapted to $S'$ and $L$ is the Levy subgroup of $P_{S'}$.
  By the isomorphism theorem
  the quotient $\lambda(\mathbb{C}^*) / H \cap \lambda(\mathbb{C}^*)$ is 
  isomorphic to the subgroup $\lambda(\mathbb{C}^*)H/H$
  of $N_L(H)/H$ and therefore finite. The $1$-parameter subgroup $\lambda(\mathbb{C}^*)$
  is nontrivial and connected, so $\lambda(\mathbb{C}^*) \subset H$ and the variety $X'$ is fixed by $\lambda$.
  This means that $X^{\prime \lambda} = X'$.

  Now we consider the action of $\lambda$ on $G/P^-$. 
  We have that $(G/P^-)_{S'} = \{x \in G/P^- \ |\ \overline{P_{S'}x}\supset G/P^-\}$,
  as $G/P^-$ is homogeneous and projective. 
  In other words $(G/P^-)_{S'}$ is the open $P$-orbit of $G/P^-$. By Proposition \ref{prop_luna_limit_stuff}
  $(G/P^-)_{S'} = \{ x \in G/P^-\ |\ \lim_{t\to 0} \lambda(t)x \in (G/P^-)^{S'}\}$, where 
  $(G/P^-)^\lambda = (G/P^-)^{S'}$ is the set of $P^-$-fixed points in $(G/P^-)_{S'}$.
  Then $(G/P^-)^{S'}$ is just the single $P^-$-fixed point $[1_G] \in (G/P^-)_{S'}$.
  As the projection $G\times^{P^-} X' \to G/P^-$ is $G$-equivariant,
  $\pi(X^{S'}) = \pi(X)^{S'} = (G/P^-)^{S'}$ is the point $[1_G]$ and $X^{S'}$ is contained
  in the fiber $\pi^{-1}([1_G])$, therefore ${X}^{S'}$ lies in the fiber over $[1_G]$.
  The fiber $\pi^{-1}([1_G])$ is precisely $X'$, thus $X^{S'} \subset X'$.
  
  Since $X'$ is fixed by $\lambda$,
  we have $X' \subset X^{S'}$, proving the result.  
  $\qed$
\end{proof}


\begin{lemma}[{\cite[Proposition 2.6]{luna_gcplvs}} for (a) and (b)]\label{lemma_color_intersection_grosses_cellules}
  Let $X$ be a wonderful $G$-variety, $S' \subset S$, $\lambda$ adapted to $S'$.
  For $D\in \Delta_X$ in $X$ there are two possibilities:
  \begin{enumerate}
    \item $\langle \rho_X(D), \lambda \rangle = 0$, $D \cap X_{S'} \neq \emptyset$ and $D$ is not $P_{S'}$-stable,
    \item $\langle \rho_X(D), \lambda \rangle \neq 0$, $D \cap X_{S'} = \emptyset$ and $D$ is $P_{S'}$-stable.
  \end{enumerate}
  Furthermore any prime $G$-stable divisor intersects $X_{S'}$ and 
  $$X = X_{S'} \cup \bigcup_{D \in \Delta_X\setminus\Delta_X(S')} D.$$
\end{lemma}

\begin{proof}
  First, note that $X_{S'}$ contains the closed orbit, which is the intersection of all the $G$-stable divisors.   
  The points (a) and (b) can be found in \cite{luna_gcplvs} proof of Proposition 2.6.

  We consider $X_\emptyset$ as in Definition \ref{def_limit_stuff} with $\emptyset \subset S$
  then $X_\emptyset \cap Y \neq \emptyset$, where $Y \subset X$ is the closed $G$-orbit. 
  For any closed $G$-stable subvariety $Z \subset X$, we have 
  $Y \subset Z$, so $Z \cap X_\emptyset \neq \emptyset$. For any $S' \subset S$ we have $X_\emptyset \subset X_{S'}$,
  so $Z\cap X_\emptyset \subset Z \cap X_{S'} \neq \emptyset$. 
    
  Next, we have that $X_\emptyset \cap G/H$ is the open $B$-orbit of $G/H$ and its complement is the union of all colors.
  For any $S' \subset S$ we have $X_\emptyset \subset X_{S'}$, so $X \setminus X_{S'}$ is a union of some colors, 
  therefore 
  \[
    X \setminus X_{S'} = \bigcup_{D \in \Delta_X, D\cap X_{S'}=\emptyset}D.
  \]
  Part (a) tells us that the set $\Delta_X(S')$ of colors that are not $P_{S'}$-stable is given by 
  $\{ D \in \Delta_X \ |\ D \cap X_{S'} \neq \emptyset\}$.
  $\qed$
\end{proof}

\section{Wonderful symmetric varieties}\label{sec_wonderful_compactifications_sym_space}

Here and later $G$ stands for a semisimple algebraic group of adjoint type.
Let $\sigma: G \to G$ be an involution and $H = G^\sigma = \{g \in G\ |\ \sigma(g) = g\}$
be the fixed point set. Then the space $G/H$ is called a \textit{symmetric} space.
We will explain the construction of the wonderful embedding of $G/H$ 
and give a precise description of $\Pic(X)$.

\begin{lemma}[{\cite[Proposition 1.1]{deconcini_procesi_csv}}]\label{lemma_dcp_torus_choice}
  Every $\sigma$-stable torus in $G$ is contained in a maximal torus of $G$ which is $\sigma$-stable.
\end{lemma}

Let us fix a $\sigma$-stable maximal torus $T$, such that the subtorus $T_1 = \{ t \in T \ |\ \sigma(t) = t^{-1}\}^\circ$ 
(the connected component, containing $id_G$) is of maximal dimension.
The involution $\sigma$ acts on the root system $R_G$ of $G$ via composition:
$\sigma(\alpha) = \alpha\circ\sigma$. With that we split $R$ into
$R_0 = \{\alpha \in R\ |\ \sigma(\alpha) = \alpha\}$ and 
$R_1 = R \setminus R_0$.

\begin{lemma}[{\cite[Lemma 1.2]{deconcini_procesi_csv}}]\label{lemma_basis_of_roots_for_sigma}
  There is a basis of simple roots 
  \[
    S= \{\beta_1,\dots,\beta_m,\alpha_1,\dots,\alpha_n\}
  \]
  of $R$, such that 
  $\{\beta_1,\dots,\beta_m\}\subset  R_0$ and $\{\alpha_1,\dots,\alpha_n\} \subset R_1$.
  Denote by $R^+$ the positive roots defined by $S$ and $R^- = -R^+$.
  With this definition $\sigma(\alpha) \in R^-_1$ for any $\alpha \in R^+_1$.

\end{lemma}

We fix a choice of a Borel subgroup $B\subset G$, such that $R^+$ are its roots.
Note that $BH$ is open in $G$.

\begin{lemma}[{\cite[pages 5-6]{deconcini_procesi_csv}}]\label{lemma_dcp_spherical_roots_not_more_that_two}
  There exists a permutation $\tilde{\sigma} : S \cap R_1 \to S \cap R_1$ of order at most $2$, such that 
  $\sigma(\alpha) = - \tilde{\sigma}(\alpha) - \beta$, where $\beta$ is a nonnegative linear combination
  of simple roots of $R_0$. 
\end{lemma}

\begin{definition}\label{def_dcp_spherical_roots}
  The weights $\bar{\alpha}_i = \alpha_i - \sigma(\alpha_i)$ are called
  the \textit{spherical roots} of $G/H$. The number of spherical roots is called 
  the \textit{rank} of $G/H$ denoted $l = \rk(G/H)$, this is also the dimension of $T_1$.
  A simple root $\alpha$ of $G$ is called \textit{exceptional}
  if $\alpha \neq \tilde{\sigma}(\alpha)$ and $(\alpha, \tilde{\sigma}(\alpha)) \neq 0$. 
  We will
  also call $\bar{\alpha}$ an \textit{exceptional} spherical root of $X$.
  We reorder the simple roots of $G$, so that $\bar{\alpha}_1, \dots, \bar{\alpha}_l$ 
  are all distinct.
\end{definition}

\begin{remark}\label{remark_at_most_two_simple_per_spherical}
  By Lemma \ref{lemma_basis_of_roots_for_sigma}
  there exist at most two simple roots $\alpha, \alpha' \in S$,
  such that $\bar{\alpha} = \bar{\alpha}'$ and if $\bar{\alpha} = \bar{\alpha}'$,
  then $\alpha = \tilde{\sigma}(\alpha')$.

  There is also an induced permutation $\tilde{\sigma} : \{1, \dots, n\} \to \{1, \dots, n\}$, defined by 
  $\tilde{\sigma}(\alpha_i) = \alpha_{\tilde{\sigma}(i)}$.
  For any spherical root $\bar{\alpha}_i$ with $1 \leq i \leq l$ 
  there are at most two simple roots $\alpha_i$ and $\alpha_{\tilde{\sigma}(i)}$ with $\tilde{\sigma}(i) > l$, such that 
  $\bar{\alpha}_i = \bar{\alpha}_{\tilde{\sigma}(i)}$.
\end{remark}


\begin{lemma}[{\cite[pages 5-6]{deconcini_procesi_csv}}]\label{lemma_fundamental_weights_under_sigma}
  Let $\alpha \in S$ be a simple root and $\omega_\alpha$ be the corresponding fundamental weight,
  then 
  \[
    \sigma(\omega_\alpha) = \begin{cases}
      \omega_{\tilde{\sigma}(\alpha)} &\text{if }\alpha \in R_0,\\
      -\omega_{\tilde{\sigma}(\alpha)} &\text{if }\alpha \in R_1.
    \end{cases}
  \]
  We will denote $\omega_{\alpha_i}$ by $\omega_i$.
  Additionally if $1 \leq i \leq l$ and $\tilde{\sigma}(i)\neq i$, then we have $\tilde{\sigma}(i) > l$. 
\end{lemma}

\begin{theorem}[{\cite[Proposition 4.7]{richardson_oiratiorg}}]
  The set $\bar{R} = \{\alpha - \sigma(\alpha)\ |\ \alpha \in R\}$ is a (possibly nonreduced)
  root system. The corresponding Weyl group is given by $N_G(T_1)/Z_G(T_1)$.
\end{theorem}

\begin{lemma}[{\cite[Lemma 1.5 and Proposition 1.8]{deconcini_procesi_csv}}]
  Let $\lambda$ be a dominant weight of $G$ and $V_\lambda$ be the irreducible representation 
  of $\tilde{G}$ with highest weight $\lambda$, where $\tilde{G}$ is the simply
  connected covering of $G$. Let $\tilde{H}$ be the preimage of $H$ in $\tilde{G}$.
  Let $V_\lambda^{\tilde{H}}$ be the subspace of $\tilde{H}$-invariant vectors in the irreducible representation $V_\lambda$.
  Then $\dim V_\lambda^{\tilde{H}} \leq 1$.
  If $V_\lambda^{\tilde{H}} \neq 0$, then $\sigma(\lambda) = -\lambda$.
  If $\sigma(\lambda) = -\lambda$, then $\dim V_{2\lambda}^{\tilde{H}} = 1$.
\end{lemma}

\begin{definition}\label{def_spherical_weights}
  A weight is called \textit{spherical} if $V_\lambda^{\tilde{H}} \neq 0$.
  A fundamental weight $\omega_\alpha$ of $G$ is called \textit{exceptional} if the simple root 
  $\alpha$ is exceptional. 

\end{definition}

Let $\lambda$ be a regular dominant spherical weight.
Let $V_\lambda$ be an irreducible representation of the simply connected covering $\tilde{G}$
of $G$ of highest weight $\lambda$
and $V_\lambda^\sigma$ be the $\sigma$-twisted representation, i.e.
$g  v = \sigma(g)  v$.
It is shown in Section 1 of \cite{deconcini_procesi_csv},
that $V_\lambda^\sigma \simeq V_\lambda^*$.
Let $v_\lambda \in V_\lambda$ be the highest weight vector. Pick $v^\lambda \in V_\lambda^*$,
such that $v^\lambda(v_\lambda) = 1$ and $v^\lambda(W) = 0$,
where $W$ is the unique $T$-stable complement of the line $\mathbb{C}v_\lambda$.
Let $\tilde{h}: V_\lambda^* \to V_\lambda$ be the unique $\sigma$-linear
isomorphism, such that $\tilde{h}(v^\lambda) = v_\lambda$.
Then $\tilde{h} \in V_\lambda \otimes V_\lambda$ and we consider the projection
$p: V_\lambda \otimes V_\lambda \to V_{2\lambda}$. Finally
the element $[p(\tilde{h})] \in \mathbb{P}(V_{2\lambda})$ is $H$-invariant.
We denote this element by $h$. Note that the action of $\tilde{G}$ on the
projective space factors through $G$.
We have the following result.
 
\begin{theorem}[{\cite[Theorem 3.1 and Section 4]{deconcini_procesi_csv}}]
  Let $X = \overline{G h} \subset \mathbb{P}(V_{2\lambda})$. Then $X$ is the wonderful embedding of $G/H$.
  As an abstract $G$-variety, $X$ is independent of the choice of $\lambda$.
\end{theorem}

Let us fix $X$ to be the wonderful embedding of $G/H$, 
with $X = \overline{G h} \subset \mathbb{P}(V_\lambda)$ for some $\lambda$ and $h$.
The next proposition gives us 
an affine chart of $X$, that allows us to study $X$ locally. 
Let $v_\lambda$ denote the highest weight vector in $V_\lambda$.
Let $W$ be the $T$-stable complement of $\mathbb{C}v_\lambda$, such that $V_\lambda = \mathbb{C}v_\lambda + W$.
Define $A \subset \mathbb{P}(V_\lambda)$ as the open subset where the coordinate on $v_\lambda$ is nonzero.
By construction in \cite{deconcini_procesi_csv} we have $h \in A$.

\begin{lemma}[{\cite[Lemma 2.2]{deconcini_procesi_csv}}]\label{lemma_dcp_affine_torus_orbit}
  The closure in $A$ of the $T_1$ orbit $T_1 h$ is isomorphic to $\mathbb{A}^l$.
  Let $x = (x_1, \dots, x_l)$ be a point in $\mathbb{A}^l$, then 
  \[
    t  (x_1, \dots, x_l) = (-\bar{\alpha}_1(t)x_1, \dots, -\bar{\alpha}_l(t)x_l),
  \]
  and the orbit $T_1 h$ corresponds to the set $\{x \in \mathbb{A}^l\ |\ x_i \neq 0\ \forall i\}$.
\end{lemma}

The variety $X$ is wonderful, so we can define the rank of $X$, as well as the boundary divisors $X_i$ 
for $i \in \{1, \dots, \rk(X)\}$ and with that the closed $G$-stable subvarieties 
$X_I = \cap_{i \in I}X_i$ for $I \subset \{1, \dots, \rk(X)\}$.
The closed $G$-orbit of $X$ is $Y = X_{\{1, \dots, \rk(X)\}}$ and is isomorphic to 
$G/P_{S_X^p}^-$.

There exists an explicit version of local structure Theorem
\ref{thm_local_structure_spherical} for wonderful varieties.

\begin{proposition}[{\cite[Proposition 2.3]{deconcini_procesi_csv}}]\label{prop_dcp_affine_subset}
  Let $U^- = \prod_{\alpha \in R_1^-} U_\alpha \subset G$, then the map 
  \begin{align*}
    U^- \times \mathbb{A}^l &\to X \\
    (u,x) &\mapsto u  x
  \end{align*}
  is an open immersion, where $\mathbb{A}^l \simeq \overline{T_1 h}\cap A$. 
  Let $V^-$ denote its image, then $V^-$ is $B^-$ stable and $X = \cup_{g \in G}gV^-$.

  Moreover $\rk(X) = l$ and for any closed $G$-stable subvarieties $X_I$ we have 
  \[
    X_I \cap V^- = \{(u,x)\in U \times \mathbb{A}^l\ |\ x_i = 0 \text{ for all }i \in I\}.  
  \]
\end{proposition}

\begin{remark}\label{remark_symmtric_roots_not_moving_colors}
  The translation of $V^-$ by the longest element of the Weyl group of $G$ is 
  an open $B$-stable subset $V$. This $V$ is isomorphic to $U \times \mathbb{A}^l$, 
  where the weights of the $T$-action on $\mathbb{A}^l$ are $\bar{\alpha}_i$ for $1 \leq i \leq l$ 
  and $U = \prod_{\alpha \in R_1^+} U_\alpha$.

  The subset $V \subset X$ is the subset $X_{Y,B}$ introduced in Definition \ref{def_spherical_b_open_subset}.
  It coincides with $X_\emptyset$ introduced in Definition \ref{def_limit_stuff}.
  Note that also $S_X^p = S \cap R_0$, since the parabolic $P_{S \cap R_0}$ stabilizes $V$.
\end{remark}

\begin{corollary}\label{corollary_symmetric_roots_move_one_color}
  The spherical roots $\bar{\alpha}_i$, defined in \ref{def_dcp_spherical_roots}
  coincide with the ones defined in \ref{def_luna_spherical_roots_and_system}.
  Any color $\alpha \in S \setminus S_X^p$ is of type $(a')$ or $(b)$ in the sense of Proposition \ref{prop_color_types}, 
  in particular any simple root $\alpha \in S$ moves at most one color.
\end{corollary}
  
\begin{proof}
  The first part is a direct consequence of Proposition \ref{proposition_t_weights_spherical_roots} and \ref{prop_dcp_affine_subset}
  and the previous Remark.
  Let $z = (id, 0) \in U \times \mathbb{A}^l$, it is a point in the closed orbit $Y$ of $X$.
  We have that $V \cap Y \simeq U = R_u(P)$, where $R_u$ denotes the unipotent radical.
  Then $z$ is the $B^-$-fixed point of $G/P^-$, and the weights of the $T$-action on the normal space $T_z X/T_z Y$
  are given by 
  \[
    t  (x_1, \dots, x_l) = (\bar{\alpha}_1(t)x_1, \dots, \bar{\alpha}_l(t)x_l)
  \]

  If $D$ is a color of type $(a)$, that would mean that there exists a simple root $\alpha \in S$, such 
  that $\alpha$ is a spherical root. Then $\alpha = \alpha_i - \sigma(\alpha_i)$ for some $1 \leq i \leq l$, a contradiction to the fact 
  that $\alpha$ is simple.
  $\qed$
\end{proof}

\begin{proposition}[{\cite[Theorem 5]{deconcini_procesi_csv}}]\label{prop_boundary_divisor_projection}
  Let $X_I \subset X$ be a closed $G$-stable subvariety. Let 
  $S_I = \{\alpha_i, \tilde{\sigma}(\alpha_i) \in S_G\ |\ i\notin I\}\cup S_X^p$.
  
  There exists an equivariant fibration $X_I \simeq G\times^{P_{S_I}^-} X^I$
  with fibers isomorphic to the wonderful compactification $X^I$ of the 
  symmetric space $\overline{L_I}/\overline{(L_I \cap H)}$, where 
  $P_{S_I}^-$ acts on $X^I$ via its adjoint Levi quotient $\overline{L_I}$.
\end{proposition}

\begin{remark}
  Note that
  $S_I = \{\alpha \in S\ |\ \bar{\alpha}_i = \alpha - \sigma(\alpha) \text{ for some }i \notin I\} \cup S_X^p$,
  this follows directly from Remark \ref{remark_at_most_two_simple_per_spherical}.
\end{remark}




Next we wish to describe the Picard group of $X$ and of any closed $G$-stable subvariety $X_I \subset X$.
This is the first step to understanding the space of curves, as $Pic(X) = A_1(X)^\vee$ by Proposition
\ref{prop_toroidal_curve_divisor_duality}.

\begin{proposition}[{\cite[Proposition 8.1]{deconcini_procesi_csv}}]\label{prop_wond_injection_pullback_injective_on_pic}
  Let $\iota: Y \to X$ denote the inclusion of the closed orbit $Y \subset X$. Then the pullback 
  \[
    \iota^*: \Pic(X) \to \Pic(Y)
  \]
  is injective.
\end{proposition}


\begin{proposition}[{\cite[Theorems 4.2 and 4.8]{deconcini_springer_csv}}]\label{prop_picard_grp_as_weight_lattice}
  The Picard group $\Pic(X)$ of $X$ is isomorphic to the sublattice $\Lambda_X$ of $X(T)^*$, generated by 
  spherical weights and exceptional fundamental weights. The rank of the Picard group $\Pic(X)$ is 
  equal to $l+s$, where $l$ is the rank of $\bar{R}$ and $s$ is the number of spherical roots $\bar{\alpha}$,
  such that $\alpha$ is exceptional.
\end{proposition}

\begin{definition}
  The variety $X$ is called \textit{exceptional} if it has exceptional spherical roots, 
  or equivalently $\rk(X) < |\Delta_X|$. Otherwise $X$ is called \textit{nonexceptional}.
\end{definition}

We can also describe the Picard group for any closed $G$-stable subvariety $X_I$ in terms of Picard groups
of wonderful symmetric varieties and projective homogeneous varieties.

\begin{proposition}[{\cite[page 7]{chirivi_maffei_pnocsv}}]\label{prop_boundary_picard_exact_sequence} 
  Let $X_I$, $X^I$ and $P_{S_I}$ be as in Proposition \ref{prop_boundary_divisor_projection}, 
  with map $\pi_I:X_I \to G/P_{S_I}^-$.
  Let $\iota_I: X^I \to X_I$ be the inclusion of the fiber $\pi_I^{-1}([1_G])$ for $[1_G] \in G/P_{S_I}^-$.
  The following sequence is exact:
  \[
  \begin{tikzcd}
    0 \arrow[r] & \Pic(G/P_{S_I}^-) \arrow[r, "\pi_I^*"] & \Pic(X_I) \arrow[r,"\iota_I^*"] & \Pic(X^I) \arrow[r] & 0
  \end{tikzcd}
  \]
\end{proposition}

\begin{corollary}\label{corollary_explicit_exact_seq_morphism}
  The map
  \begin{align*}
    \Pic(G/P_{S_I}^-) \oplus \Pic(X^I) &\to \Pic(X_I)\\
    (Z_{base}, Z_{fiber}) &\mapsto \pi_I^*Z_{base} + \overline{BZ_{fiber}}
  \end{align*}
  is an isomorphism of abelian groups.
\end{corollary}

\begin{proof}
  Let $S_I \subset S$ be the set of simple roots of $G$, such that $P_{S_I}$ is 
  the parabolic associated to $S_I$.
  Recall that $X_I = G \times^{P_{S_I}^-} X^I$ by Proposition \ref{prop_boundary_divisor_projection}.
  By Lemma \ref{lemma_parabolic_induction_and_limit_stuff}
  we have $X^I \simeq (X_I)^{S_I}$.
  
  By Proposition \ref{prop_luna_limit_stuff} we have that $\Delta_{X_I}(S_I) = \Delta_{X^I}$,
  where the bijection is given by $Z \mapsto \overline{B Z}$ for $Z \in \Delta_{X^I}$
  and $D \mapsto D \cap X^I$ for $D \in \Delta_{X_I}(S_I)$.
  This bijection also tells us that $Z = \iota^*(D) \neq 0$ in $\Pic(X^I)$. 
  
  Suppose now that $\iota^*(D) = 0$ for some $D \in \Delta_{X_I}$. Then 
  $D \in \Ima(\pi_I^*)$, so there exists some divisor $Z \in \Pic(G/P_{S_I}^-)$, 
  such that $D = \pi_I^{-1}(Z)$. Then $\pi_I(D) = \pi_I(\pi_I^{-1}(Z)) = Z$
  and since $D$ is $B$-stable and irreducible and $\pi_I$ is a $G$-equivariant projection,
  we get that $Z$ is also $B$-stable and irreducible.
  $\qed$ 
\end{proof}

Recall that any line bundle on $X$ can be identified with a $B^-$-weight as follows:
first replace $G$ by its simply connected covering, then
let $\psi: G/B^- \to X$ be the composition of the map between flag varieties $G/B^- \to Y$ and 
the inclusion $Y \subset X$. Then for any divisor $D$ on $X$ we have 
\[
  \psi^*\mathcal{O}_X(D) \simeq G/B^- \times^{B^-} \mathbb{C}_\chi = \mathcal{L}_\chi
\]
for some $B^-$-weight $\chi$. 
The fact that $\iota^*: \Pic(X) \to \Pic(Y)$ is injective for wonderful symmetric varieties, 
means that the weight $\chi$ determines the color $D$.

First we need to prove that any color $D \in \Delta_X$ is moved by at most two simple roots.

\begin{lemma}\label{lemma_each_color_is_moved_by_at_most_2_roots}
  Let $D$ be a color in $\Delta_X$. If $\{D\} = \Delta_X(\alpha_i) = \Delta_X(\alpha_j)$
  for distinct $\alpha_i, \alpha_j \in S \setminus S_X^p$, then $\alpha_i = -\tilde{\sigma}(\alpha_j)$
  and $(\alpha_i, \alpha_j) = 0$. Otherwise there exists unique $\alpha_i \in S$,
  such that $\{D\} = \Delta_X(\alpha_i)$.
\end{lemma}

\begin{proof}
  Let $D$ be a color, suppose that there exist two distinct simple roots $\alpha_i, \alpha_j \in S \setminus S_X^p$,
  such that $\{D\} = \Delta_X(\alpha_i) = \Delta_X(\alpha_j)$.
  Then by Proposition \ref{prop_wond_var_rho_value}
  and Corollary \ref{corollary_symmetric_roots_move_one_color} for any spherical root
  $\alpha - \sigma(\alpha) \in \Sigma_X$ we have
  \[
    \alpha_i^\vee (\alpha - \sigma(\alpha)) = \alpha_j^\vee (\alpha - \sigma(\alpha))
  \]
  or
  \[
    (\alpha_i^\vee - \alpha_j^\vee) (\alpha - \sigma(\alpha)) = 0.
  \]
  So $\alpha_i^\vee - \alpha_j^\vee$ is orthogonal to $\bar{R}$
  and therefore $\alpha_i^\vee - \alpha_j^\vee \in \langle R_0^\vee \rangle$.
  So we have
  \[
    \alpha_i - \alpha_j = \sigma(\alpha_i - \alpha_j)
  \]
  or
  \[
    \alpha_i - \sigma(\alpha_i) = \alpha_j - \sigma(\alpha_j)
  \]
  and by Remark \ref{remark_at_most_two_simple_per_spherical} we have
  \[
    \alpha_i = \tilde{\sigma}(\alpha_j).  
  \]
  The order of the involution $\tilde{\sigma}$ is two, so
  any color is moved by at most two simple roots, namely $\alpha_i$ and $\tilde{\sigma}(\alpha_i)$.

  By Remark \ref{remark_at_most_two_simple_per_spherical} we have that $\alpha_i$ and $\alpha_j$
  are the only simple roots, associated to the spherical root $\bar{\alpha}_i = \bar{\alpha}_j$.
  Without loss of generality we can assume that $i < j$,
  so $1 \leq i \leq l$ as in Definition \ref{def_dcp_spherical_roots}.
  For $I = \{1, \dots l\}\setminus \{i\}$ we introduce
  $S_I = \{\alpha_i, \tilde{\sigma}(\alpha_i)\ |\ i \notin I\} \cup S_X^p$,
  as in Proposition \ref{prop_boundary_divisor_projection}.
  Let $P = P_{S_I}$ be the parabolic subgroup, associated to $S_I$ and $L = G^{S_I}$ be the Levi
  quotient of $P$.
  We have $X_I = G \times^{P^-}X^I$, where $X^I$ is an wonderful
  symmetric $L$-variety with $\rk(X^I) = \rk(X_I) = 1$.
  We also have that
  \[
    \Delta_X(S_I) = \bigcup_{\alpha\in S_I} \Delta_X(\alpha) = \Delta_X(\alpha_i) \cup \Delta_X(\alpha_j),
  \]
  as $\Delta_X(\alpha) = \emptyset$ for any $\alpha \in S_X^p$.

  The root system of $L$ is generated by the roots $\alpha_i$ and $\alpha_j$ and $S_X^p$.
  By Proposition \ref{prop_picard_grp_as_weight_lattice} we have 
  \[
    \rk(\Pic(X^I)) = \begin{cases}
      1 &\text{if } (\alpha_i, \alpha_j) = 0,\\
      2 &\text{else}. 
    \end{cases}
  \]
  On the other hand the Picard group $\Pic(X^I)$ is generated by $\Delta_{X^I}$,
  so $|\Delta_{X^I}| = \rk(\Pic(X^I))$. 
  By Lemma \ref{lemma_parabolic_induction_and_limit_stuff} we have
  $X^I = X^{S_I}$. We apply Proposition \ref{prop_luna_limit_stuff} to get that
  $\Delta_X(S_I)$ is in bijection with $\Delta_{X^I}$, so
  \[
    |\Delta_X(S_I)| = \begin{cases}
      1 &\text{if } (\alpha_i, \alpha_j) = 0,\\
      2 &\text{else}.
    \end{cases}
  \]
  By Corollary \ref{corollary_symmetric_roots_move_one_color} we have
  $|\Delta_X(\alpha)| = 1$ for any $\alpha \in S\setminus S_X^p$, so
  \[
    \Delta_X(S_I) = \begin{cases}
      \Delta_X(\alpha_i) = \Delta_X(\alpha_j)   &\text{if } (\alpha_i, \alpha_j) = 0,\\
      \Delta_X(\alpha_i) \mathbin{\mathaccent\cdot\cup} \Delta_X(\alpha_j) &\text{else}.
    \end{cases}
  \]
  The result follows.
  $\qed$
\end{proof}

Now that we have studied $\{\alpha \in S\setminus S_X^p\ |\ D \in \Delta_X(\alpha)\}$
for any color $D \in \Delta_X$ we can state the following proposition,
describing the relationship between colors and $B^-$-weights.

\begin{proposition}\label{prop_explicit_basis_picard_lattice}
  Let $X$ be a wonderful symmetric variety of rank $l$.
  Let $D \in \Delta_X$ be a color.
  Suppose that $D$ is moved by nonexceptional simple roots. Then one of the following holds:
  \begin{itemize}
    \item $\exists!\alpha \in S: D \in \Delta_X(\alpha), \alpha = -\sigma(\alpha)$,
      then $\iota^*(\mathcal{O}_X(D))\simeq \mathcal{L}_{2\omega_\alpha}$;
    \item $\exists!\alpha \in S: D \in \Delta_X(\alpha), \alpha \neq -\sigma(\alpha)$,
      then $\iota^*(\mathcal{O}_X(D))\simeq \mathcal{L}_{\omega_\alpha}$;
    \item $\exists \alpha,\alpha' \in S: \alpha \neq \alpha', D \in \Delta_X(\alpha) = \Delta_X(\alpha')$,
      then $\iota^*(\mathcal{O}_X(D))\simeq \mathcal{L}_{\omega_\alpha + \omega_{\alpha'}}$;
  \end{itemize}
\end{proposition}

\begin{proof}
  Let $D$ be a color of $X$ and let $\iota: Y \to X$ be the inclusion of the closed orbit.
  We have $Y \simeq G/P^-$, where $P = P_{S_X^p}$. The Picard group of $Y$ 
  is generated by the Schubert divisors $D_\alpha(G/P^-)$ for $\alpha \in S\setminus S_X^p$.
  We know that $N_1(Y)$ is dual to $\Pic(X)$ and any divisor on $Y$ is defined by
  the intersections with Schubert curves.
  In particular for $\iota^*D$ we have
  \[
    \iota^*D = \sum_{\alpha \in S \setminus S_X^p} \langle \iota^*D, [C_\alpha]\rangle D_\alpha(G/P^-),
  \]
  where $[C_\alpha]$ denote a Schubert curve class.

  In order to describe $\iota^*D$ we need to compute $\langle \iota^*D, [C_\alpha]\rangle$
  for any Schubert curve class $[C_\alpha]$.
  By projection formula we have
  \[
    \langle \iota^*D, [C_\alpha] \rangle = \langle D, \iota_*[C_\alpha] \rangle.
  \]
  By Corollary \ref{corollary_symmetric_roots_move_one_color} any simple root $\alpha$
  is either in $S_X^{a'}$ or in $S_X^b$ and so by
  Lemma \ref{lemma_wond_color_closed_orbit_curves} we can write $\iota_*[C_\alpha]$
  as a multiple of a curve class $[C_{D'}]$ for some $D' \in \Delta_X$.
  Here $[C_{D'}]$ is a curve class on $X$, such that
  \[
    \langle D, [C_{D'}] \rangle = \begin{cases}
      1 &\text{if }D = D',\\
      0 &\text{else},
    \end{cases}
  \]
  see Lemma \ref{lemma_wond_color_dual_curve_description}.

  Let $S_D = \{\alpha\in S\ |\ D \in \Delta_X(\alpha)\}$ be the set of simple roots moving $D$.
  By Corollary \ref{corollary_symmetric_roots_move_one_color} each root moves at most one color,
  so $\Delta_X(S_D) = \{D\}$.
  For any simple root $\alpha \in S \setminus S_D$ we have $D \notin \Delta_X(\alpha)$ and
  $[C_\alpha]$ is a multiple of $[C_{D'}]$ for some $D' \neq D$, so
  $\langle \iota^*D, [C_{\alpha}] \rangle = 0$.
  With that we can write
  \[
    \iota^*D = \sum_{\alpha \in S_D} \langle \iota^*D, [C_\alpha]\rangle D_\alpha(G/P^-).
  \]
  By Lemma \ref{lemma_each_color_is_moved_by_at_most_2_roots} every color is moved either by one root $\alpha_D$
  or by two orthogonal simple roots $\alpha_D$ and $\tilde{\sigma}(\alpha_D')$.
  This means that $S_D$ has either one or two elements.

  Suppose that $S_D = \{\alpha_D\}$ contains a single simple root.  
  If $\alpha_D = -\sigma(\alpha_D)$, we have that
  $\alpha_D - \sigma(\alpha_D) = 2\alpha_D \in \Sigma_X$.
  By Lemma \ref{lemma_wond_color_closed_orbit_curves} we have $\iota_*[C_{\alpha_D}] = 2 [C_D]$,
  so $\langle \iota^*D, [C_{\alpha_D}] \rangle = 2$ and
  \[
    \iota^*D = 2 D_{\alpha_D}(G/P^-).
  \]
  If $\alpha_D \neq -\sigma(\alpha_D)$, then $2\alpha_D \notin \Sigma_X$ and $\alpha \in S^b_X$, 
  so by Lemma \ref{lemma_wond_color_closed_orbit_curves} we have $\iota_*[C_{\alpha_D}] = [C_D]$ and similarly
  \[
    \iota^*D = D_{\alpha_D}(G/P^-).
  \]
  Suppose now that $S_D = \{\alpha_D, \tilde{\sigma}(\alpha_D)\}$ and $(\alpha, \tilde{\sigma}(\alpha)) = 0$.
  Then neither $2\alpha_D$ not $2\tilde{\sigma}(\alpha_D)$ are spherical roots,
  so by Lemma \ref{lemma_wond_color_closed_orbit_curves} we have $\iota_*[C_{\alpha_D}] = [C_D]$ and
  $\iota_*[C_{\tilde{\sigma}(\alpha_D)}] = [C_D]$, therefore
  \[
    \iota^*D = D_{\alpha_D}(G/P^-) + D_{\tilde{\sigma}(\alpha_D)}(G/P^-).
  \]
  Lastly we note that $\mathcal{O}_Y(D_\alpha(Y)) = \mathcal{L}_{\omega_\alpha}$ and the result follows.
  $\qed$ 
\end{proof}

\chapter{Moduli space of stable maps of a wonderful symmetric variety}

\section{Moduli space of stable maps}

\begin{definition}
  Let $X$ be a smooth projective complex variety and $\eta \in N_1(X)$. Let $M_{0,n}(X, \eta)$ 
  be the set of isomorphism classes of pointed maps $(C, p_1, \dots, p_n, \mu)$,
  where $C$ is a projective smooth rational curve, $p_1, \dots, p_n$ are distinct points on $C$
  and $\mu : C \to X$ is a morphism, such that $\mu_*([C]) = \eta$. When $X = \{pt\}$ is a point 
  we will write $M_{0,n} = M_{0,n}(\{pt\}, 0)$.

  A $n$-pointed rational complex \textit{quasi-stable} curve 
  \[
    (C,p_1, \dots, p_n)  
  \]
  is a projective connected reduced (at worst) nodal rational curve with $n$ distinct
  smooth marked points. Let $S$ be an algebraic scheme over $\mathbb{C}$. A \textit{family}
  of $n$-pointed rational quasi-stable curves over $S$ is a flat projective map $\pi: \mathcal{C} \to S$
  with $n$ sections $p_1, \dots, p_n$, such that each geometric fiber $(\mathcal{C}_s, p_1(s), \dots, p_n(s))$
  is an $n$-pointed rational quasi-stable curve. 
  A \textit{family of maps} over $S$ from $n$-pointed rational quasi-stable
  curves to $X$ consists of data $(\pi: \mathcal{C} \to S, \{p_i\}_{1 \leq i \leq n}, \mu : \mathcal{C} \to X)$:
  \begin{itemize}
    \item A family of $n$-pointed rational quasi stable curves $\pi: \mathcal{C} \to S$ with $n$ sections $\{p_1, \dots, p_n\}$.
    \item A morphism $\mu : \mathcal{C} \to X$.
  \end{itemize}
  Two families of maps over $S$,
  \[
    (\pi: \mathcal{C} \to S, \{p_i\}_{1 \leq i \leq n}, \mu : \mathcal{C} \to X), 
      (\pi': \mathcal{C}' \to S, \{p_i'\}_{1 \leq i \leq n}, \mu' : \mathcal{C}' \to X)
  \]
  are \textit{isomorphic} if there exists a scheme isomorphism $\tau : \mathcal{C} \to \mathcal{C}'$
  satisfying $\pi = \pi' \circ \tau, p_i' = \tau \circ p_i, \mu = \mu' \circ \tau$.
  When $\pi: \mathcal{C} \to \Spec(\mathbb{C})$ is the structure map, 
  $(\pi: \mathcal{C} \to S, \{p_i\}_{1 \leq i \leq n}, \mu : \mathcal{C} \to X)$ is written as $(C, \{p_i\}, \mu)$.

  Let $(C, \{p_i\}, \mu)$ be a map from an $n$-pointed rational quasi-stable curve to $X$. The \textit{special points}
  of an irreducible component $E \subset C$ are the marked points and the component intersections of $C$, that 
  lie on $E$. The map $(C, \{p_i\}, \mu)$ is \textit{stable} if every component $E \simeq \mathbb{P}^1$, 
  that is mapped to a point by $\mu$ contains at least three special points.
  A family of pointed maps $(\pi: \mathcal{C} \to S, \{p_i\}, \mu)$ is \textit{stable} if the pointed map on each
  geometric fiber of $\pi$ is stable.

  Let $\eta \in N_1(X)$. A map $\mu: C\to X$ \textit{represents} 
  $\eta$ if the $\mu$-push-forward of the fundamental class $[C]$ equals $\eta$.

  We define a functor $\overline{\mathcal{M}}_{0,n}(X, \eta)$ from the category of complex 
  algebraic schemes to sets:
  \[
    \overline{\mathcal{M}}_{0,n}(X, \eta)(S) = \left\{\begin{array}{l}
      \text{isomorphism classes of stable families over $S$}\\
        \text{of maps from $n$-pointed rational curves}\\
          \text{to $X$ representing the class $\eta$}
    \end{array}\right\}
  \]
\end{definition}

\begin{theorem}[{\cite[Theorem 1]{fulton_pandharipande_nosmaqc}}]\label{thm_fulton_pandharipande_thm_1}
  There exists a projective coarse moduli space $\overline{M}_{0,n}(X, \eta)$, representing 
  the functor $\overline{\mathcal{M}}_{0,n}(X, \eta)$ and containing $M_{0,n}(X, \eta)$ as an open subset.
\end{theorem}

\begin{definition}
  Let $X$ be a smooth projective variety, then $X$ is called \textit{convex} if 
  $H^1(\mathbb{P}^1, \mu^*(T_X)) = 0$ for every map $\mu: \mathbb{P}^1 \to X$.
\end{definition}

\begin{theorem}[{\cite[Theorem 2]{fulton_pandharipande_nosmaqc}}]\label{thm_fulton_pandharipande_thm_2}
  Let $X$ be a projective, smooth, convex variety. 
  \begin{itemize}
    \item $\overline{M}_{0,n}(X, \eta)$ is a normal, projective variety of pure dimension
    \[
      \dim(X) - \langle K_X, \eta \rangle + n - 3.
    \]
    \item $\overline{M}_{0,n}(X, \eta)$ is locally a quotient of a smooth variety by a finite group.
  \end{itemize}
\end{theorem}

Projective homogeneous varieties are convex, but there is a stronger result by Thomsen,
that states that $\overline{M}_{0,n}(X,\eta)$ is irreducible.

\begin{proposition}[{\cite[Theorem 1]{thomsen_im0ngp}}]\label{prop_homogeneous_variety_mod_space_irreducible}
  Let $X$ be a flag variety and $\eta \in N_1(X)$ an effective class. The space
  $\overline{M}_{0,n}(X,\eta)$ is irreducible for every non negative integer $n$.
\end{proposition}

\begin{remark}
  Clearly $\{pt\}$ is projective homogeneous,
  so in particular we get that $\dim(\overline{M}_{0,n}) = n - 3$ for $n \geq 3$
  and $M_{0,n}$ is open dense in $\overline{M}_{0,n}$.
  For $n < 3$ there is a single isomorphism class of maps $\mathbb{P}^1 \to \Spec(\mathbb{C})$,
  so $M_{0,n} = \{pt\}$.
\end{remark}

\begin{definition}
  Let $X$ be a $G$-variety with a unique open $G$-orbit $X^\circ$.
  Denote by $M_{0,n}^\circ(X, \eta) \subset M_{0,n}(X,\eta)$ the subset of stable maps $(\mathbb{P}^1, p_1, \dots, p_n, \mu)$, 
  such that the image of $\mu$ meets $X^\circ$.
\end{definition}

\begin{definition}
  A proper complex curve $C$ over $\mathbb{C}$ is called \textit{smoothable} if there is an irreducible
  pointed complex variety $(0,T)$ and a proper flat family of curves $g : Z \to T$ such that 
  $C \simeq g^{-1}(0)$ and the generic fiber of $g$ is smooth.
  
  Any such $g : Z \to T$ is called a \textit{smoothing} of $C$.
\end{definition}

\begin{proposition}[{\cite[Example II.1.12]{kollar_rcoav}}]
  Let $C = \sum_{i = 0}^n C_i$ be a curve, which is obtained from the smooth and irreducible
  curve $C_0$ by attaching trees of smooth rational curves to it. 

  Then $C$ is smoothable.
\end{proposition}

\begin{remark}
  Note that the curves in $\overline{M}_{0,0}(X, \eta)$ are smoothable.
\end{remark}

Recall the definition of $\Hom(C,X)$ from \cite{kollar_rcoav}.
It is a scheme, locally of finite type.

\begin{theorem}[{\cite[Theorem II.1.7]{kollar_rcoav}}]\label{thm_kollar_ii.1.7}
  Let $X$ be a smooth projective variety, let $\mu: C\to X$ be a curve on $X$.
  The Zariski tangent space of $\Hom(C, X)$ at $\mu$ is isomorphic to 
  \[
    H^0(C, \mu^*T_X).  
  \]
  The dimension of every irreducible component of $\Hom(C, X)$ at $\mu$ 
  is at least 
  \[
    h^0(C, \mu^*T_X) - h^1(C, \mu^*T_X).
  \]
\end{theorem}

\begin{proposition}[{\cite[Theorem II.1.11]{kollar_rcoav}}]\label{thm_kollar_ii.1.11}
  Let $C$ be a proper smoothable curve over $\mathbb{C}$. Then there 
  is a normal pointed complex variety $(0, S)$ and 
  a flat and proper family of curves $f: \mathcal{C} \to S$ such that:
  \begin{enumerate}
      \item[$(1)$] $f^{-1}(0) \simeq C$;
      \item[$(2)$] if $V$ is an irreducible pointed complex variety and $t: v \in V \to 0 \in S$ 
          a morphism such that 
          \[
              [V\times_S \mathcal{C} \to V] \simeq [V\times_{\Spec(\mathbb{C})} C \to V],
          \]
          then $t(V) = {0}$;
      \item[$(3)$] $\dim S = -3\chi(\mathcal{O}_C) + \dim \Aut(C)$.
  \end{enumerate}
\end{proposition}

Recall the definition of $\mathrm{Hom}(\mathcal{C}/S, X \times S/S)$ from \cite{kollar_rcoav}.
It is a scheme, locally of finite type.

\begin{proposition}[{\cite[Proposition II.1.13]{kollar_rcoav}}]\label{prop_kollar_ii.1.13}
  Let $C$ be a proper algebraic curve. Let $h : C \to X$ be a morphism
  to a quasi projective variety $X$ of pure dimension $\dim X$. Assume that $X$ is a local 
  complete intersection and the image of every component of $C$ intersects the smooth locus of $X$.
  Assume that $C$ is smoothable and let $f: \mathcal{C} \to S$ be as in Proposition above. Then 
  \[
      \dim_{[h]} \mathrm{Hom}(\mathcal{C}/S, X \times S/S) \geq 
          - \langle K_X, C \rangle + (\dim X - 3)\chi(\mathcal{O}_C) + \dim \mathrm{Aut}(C)
  \]
\end{proposition}

\begin{remark}
  Note that \cite{kollar_rcoav} actually proves that any irreducible component of 
  $\mathrm{Hom}(\mathcal{C}/S, X \times S/S)$ containing $[h]$ has dimension at least 
  $- \langle K_X, C \rangle + (\dim X - 3)\chi(\mathcal{O}_C) + \dim \mathrm{Aut}(C)$.
\end{remark}

\begin{proposition}\label{prop_moduli_space_dimension}
  Let $X$ be a smooth projective variety, let $\eta \in N_1(X)$ be an effective curve class.
  Then any irreducible component of $\overline{M}_{0,n}(X, \eta)$ has dimension at least
  \[
      - \langle K_X, \eta \rangle + \dim X - 3 + n
  \]
  Moreover a point $(\mu: C \to X, p_1, \dots, p_n)$ is smooth in $\overline{M}_{0,n}(X, \eta)$
  if $H^1(C, \mu^*T_X) = 0$.
\end{proposition}

The fact that any irreducible component of $\overline{M}_{0,n}(X, \eta)$ is of dimension at least
$- \langle K_X, \eta \rangle + \dim X - 3 + n$ can be found in \cite[Lemma 4.2]{harris_roth_starr_rcohold}.
Also in the first paragraph of the proof of \cite[Proposition 7.4]{behrend_manin_sosmagwi}
it is shown that a point in the moduli space is smooth if $H^1(C, \mu^*T_X) = 0$.
However both articles use different notations to the ones presented here.

\begin{proof}
  The map $\overline{M}_{0, n + 1}(X, \eta) \to \overline{M}_{0, n}(X, \eta)$ is flat 
  of relative dimension one, see \cite{knudsen_pomsosc}, so we may assume $n = 0$
  for the first statement.
  
  Let $\mu : C \to X$ be a map in an irreducible component $M \subset M_{0,0}(X, \eta)$.
  Let $n_i$ denote the number of irreducible components of $C$ with $3-i$ nodes for $i \in \{1, 2, 3\}$. 
  Note that, because $C$ is connected and stable,
  $n_3 > 0$ if and only if $C$ is irreducible and in that case $n_1 = n_2 = 0$.
  
  The automorphisms group $\Aut(C)$ consists of automorphisms exchanging the irreducible components and 
  automorphisms of each irreducible component of $C$. As there is a finite number of irreducible components, 
  the image of $\Aut(C)$ in the permutation group of the set of components is finite,
  so $\dim \Aut(C)$ is equal to that of the subgroup, fixing the components.
  Let $C_0$ be an irreducible component of $C$, then $C_0 \simeq \mathbb{P}^1$ and the dimension of 
  the group of automorphisms of $C_0$, fixing the nodes 
  is equal to $3 - |\{x \in C_0\ |\ x \text{ is a node in }C\}|$. 
  Set $n = n_1 + 2 n_2 + 3 n_3$, then $n = \dim(\Aut(C))$. This allows us to consider the space $\overline{M}_{0,n}$,
  which will now have $3$ marked points on every component, so for any $(C, p_1, \dots, p_n)$ we have that 
  $\dim \Aut(C, p_1, \dots, p_n) = 0$.
  Write $(C, p_1, \dots, p_n)$ for points in $\overline{M}_{0,n}$ with the obvious morphism $C \to \Spec(\mathbb{C})$.
  By \cite{knudsen_pomsosc} $\overline{M}_{0,n+1}$ is the universal curve over $\overline{M}_{0,n}$, so
  the fiber over some point $(C, p_1,\dots,p_n) \in \overline{M}_{0,n}$ is isomorphic to $C$.
  We write $\mathcal{C} = \overline{M}_{0,n + 1}$ and $S = \overline{M}_{0,n}$ for brevity.
  Note that the pointed normal complex variety $(C, p_1, \dots, p_n) \in S$ together 
  with the map $\mathcal{C} \to S$ give an explicit example of the family in Proposition \ref{thm_kollar_ii.1.11}.
  
  We get a map $p: \Hom(\mathcal{C}/S, X \times S/S) \to \overline{M}_{0,0}(X,\eta)$ given by the composition 
  of maps $\Hom(\mathcal{C}/S, X \times S/S) \to \overline{M}_{0,n}(X,\eta)$ and
  $\overline{M}_{0,n}(X,\eta) \to \overline{M}_{0,0}(X,\eta)$.
  The map $\Hom(\mathcal{C}/S, X \times S/S) \to \overline{M}_{0,n}(X,\eta)$ is dominant by definition 
  and injective onto an open subset, that contains $(C, p_1, \dots, p_n, \mu)$,
  because $(C, p_1, \dots, p_n, \mu)$ has no automorphisms as a stable map. The second map
  $\overline{M}_{0,n}(X,\eta) \to \overline{M}_{0,0}(X,\eta)$ is 
  flat and proper of relative dimension $n$, so $p$ is dominant.

  Let $M \subset \overline{M}_{0,0}(X, \eta)$ be an irreducible component, that contains $(C, \mu)$, 
  denote by $H$ an irreducible component of $\Hom(\mathcal{C}/S, X \times S / S)$ 
  such that $p(H) \subset M$ is an open subset, that again contains $(C, \mu)$. 
  
  The fiber $F$ of $\overline{M}_{0,n}(X, \eta) \to \overline{M}_{0,0}(X, \eta)$ is of dimension $n$
  and by Proposition \ref{prop_kollar_ii.1.13} we have that
  \[
      \dim_{\mu} H \geq 
          -\langle K_X, \eta \rangle + \dim X - 3 + n
  \]

  Altogether
  \begin{align*}
      \dim_{[\mu]} M & \geq \dim \mathrm{Im}(H \to \overline{M}_{0,0}(X, \eta)) \\
          & \geq \dim_{\mu} \mathrm{Hom}(\mathcal{C}/S, X \times S /S) - \dim F \\
          & = -\langle K_X, \eta \rangle + \dim X - 3,
  \end{align*}
  proving the first part of the result.

  Now we prove that if $H^1(C, \mu^*T_X) = 0$, then the point
  $(\mu: C \to X)$ is a smooth point of $\overline{M}_{0,n}(X, \eta)$.
  We have the open neighbourhood $\Hom(\mathcal{C}/S, X \times S/S)$
  of $(\mu: C \to X, p_1, \dots, p_n)$ in $\overline{M}_{0,n}(X, \eta)$,
  so we only need to prove that
  $(\mu: C \to X, p_1, \dots, p_n)$ is smooth in $\Hom(\mathcal{C}/S, X \times S/S)$.

  By Proposition \ref{thm_kollar_ii.1.7} we have that $\Hom(C, X)$ is smooth
  and of minimal dimension at $\mu$, as $H^1(C, \mu^*T_X) = 0$.
  This is the fiber of the morphism
  $\Hom(\mathcal{C}/S, X \times S/S) \to S$ at $0 \in S$.

  By \cite[Theorem I.2.15.4]{kollar_rcoav} we have that if
  $\Hom(C,X)$ is of minimal dimension at $\mu$, then the map
  $\Hom(\mathcal{C}/S, X \times S/S) \to S$ is flat at $\mu$.
  Since the fiber over $0 \in S$ is $\Hom(C,X)$ and is smooth, the variety
  $\Hom(\mathcal{C}/S, X \times S/S) \to S$ is therefore smooth at $\mu$ if and only if
  $0 \in S$ is smooth. This is true, since $S = \overline{M}_{0,n}$ is smooth.

  Therefore $\mu$ is smooth in $\overline{M}_{0,n}(X, \eta)$.
  Then the map $\overline{M}_{0,n}(X, \eta) \to \overline{M}_{0,n-1}(X, \eta)$
  is flat for any $n$ and smoothness descends via flat morphisms.

  Lastly for arbitrary $n \geq 0$ we have that $\overline{M}_{0,n+1}(X, \eta) \to \overline{M}_{0,n}(X, \eta)$
  is flat. For any point $(\mu: C \to X, p_1, \dots, p_n, p_{n+1})$ the image is
  $(\mu: C \to X, p_1, \dots, p_n)$, which is smooth. The fiber over $(\mu: C \to X, p_1, \dots, p_n)$
  is $\{(\mu: C \to X, p_1, \dots, p_n, x)\ |\ x \in C\}$.
  Since marked points are smooth for stable maps we have that the point
  $(\mu: C \to X, p_1, \dots, p_n, p_{n+1})$ is smooth in the fiber
  and in the base, so by a property of flat morphisms
  for locally Noetherian schemes it is smooth in $\overline{M}_{0,n+1}(X, \eta)$.

  $\qed$
\end{proof}

The open set $M^\circ_{0,n}(X,\eta)$ is particularly interesting, because the tangent bundle 
is globally generated on the open $G$-orbit which gives a result similar to the 
one in the case of a convex variety $X$.

\begin{corollary}\label{corollary_moduli_circ_exact_dimension}
  For $X$ a $G$-variety with a dense $G$-orbit,
  the points of $M^\circ_{0,n}(X, \eta)$ are smooth in $M_{0,n}(X, \eta)$ and 
  \[
    \dim M^\circ_{0,n}(X, \eta) = \dim X  + \langle - K_X, \eta \rangle + n - 3.
  \]
\end{corollary}

\begin{proof}
  Any $(\mu: C \to X, p_1, \dots, p_n) \in M_{0.n}^\circ(X,\eta)$ meets the open orbit of $X$,
  therefore $\mu^* T_X$ is globally generated and $H^1(C, \mu^*T_X) = 0$.
  $\qed$
\end{proof}

\begin{definition}
  For a moduli space of rational curves $\overline{M}_{0,n}(X, \eta)$ with $n > 0$ marked points $p_1, \dots p_n$ 
  define $n$ maps
  \begin{align*}
    ev_{X,i}: \overline{M}_{0,n}(X, \eta) &\to X\\
    (f:C\to X) &\mapsto f(p_i)
  \end{align*}
  Sometimes indices can be dropped if there is no danger of confusion.
\end{definition}

\section{General reducibility of $M_{0,0}(X,\eta)$}\label{sec_mod_space_reducible}

Let $G$ be a semisimple algebraic group of adjoint type and
let $X$ be the wonderful embedding of a symmetric space $G/H$
of rank $\rk(X) = l$.

\begin{lemma}\label{lemma_positivity_of_curve_meeting_oo}
  Let $\eta \in N_1(X)\setminus \{0\}$ be the class of an irreducible curve meeting the open $G$-orbit. 
  Then $\langle D, \eta \rangle \geq 0$ for any $B$-stable effective divisor $D$ of $X$.
\end{lemma}

\begin{proof}
  By Proposition \ref{cycles_spherical} the cone of effective divisors on $X$ 
  is generated by the closures of codimension one $B$-orbits. 
  The codimension one $B$-orbit closure in $X$ are the colors and the boundary divisors.
  By Proposition \ref{nef_cone} the colors of $X$ generate the nef cone, so $\langle D, \eta \rangle \geq 0$
  for any color $D \in \Delta_X$, because $\eta$ is effective.

  Let $C$ be an irreducible representative of $\eta$, meeting the open orbit.
  Then $C \cap X_i$ is of dimension $0$ for any boundary divisor $X_i$, 
  and we get that $\langle X_i, \eta \rangle \geq |C \cap X_i| \geq 0$.
  $\qed$
\end{proof}

\begin{lemma}\label{lemma_positivity_of_closed_orbit_class}
  Let $\eta \in N_1(X)$ be an effective curve class, let $G/P^- = \cap_{i = 1}^l X_i$ be the closed orbit of $X$. 
  Then there exists an effective curve class $\eta_{G/P^-} \in N_1(G/P^-)$, such that 
  $\iota_* \eta_{G/P^-} = 2 \eta$. 

  If $X$ is a wonderful group compactification, then there exists 
  $\eta_{G/P^-}$ effective such that $\iota_* \eta_{G/P^-} = \eta$.
\end{lemma}

\begin{proof}
  Let $\eta_{G/P^-}$ be a curve class on $G/P^-$. Then $\eta_{G/P^-} = \sum_{\alpha \in S\setminus S_X^p}c_\alpha [C_\alpha]$,
  where $[C_\alpha]$ is the curve class, dual to the color $D_\alpha(G/P^-)$.
  Let $\iota: G/P^- \to X$ be the inclusion of the closed orbit in $X$.
  Then by Lemma \ref{lemma_wond_color_closed_orbit_curves} we have that $\iota_*[C_\alpha] = k [C_D]$, 
  where $D \in \Delta_X(\alpha)$ and $k = 1$, if $\alpha \in S_X^b$ or $k = 2$, if $\alpha \in S_X^a$.
  We have that $2A_1(X) \subset \Ima(\iota_*)$.

  Let $\eta = \sum_{D \in \Delta_X} c_D [C_D]$ be an effective curve class on $X$.
  Suppose that $D$ is moved by a single simple root $\alpha \in S_X^{b}$,
  then we set $c_\alpha = 2c_D$.
  Otherwise for any color $D \in \Delta_X(\alpha)$ for some $\alpha \in S\setminus S_X^p$
  we set $c_\alpha = c_D$. 
  Then $\eta_{G/P^-} = \sum_{\alpha \in S\setminus S_X^p}c_\alpha [C_\alpha]$
  is effective as $c_\alpha \geq 0$ for every $\alpha \in S \setminus S_X^p$ and $\iota_*(\eta_{G/P^-}) = 2\eta$.

  Suppose now that $X$ is a wonderful group compactification. Then for any simple root $\alpha \in S$ 
  we have $\sigma(\alpha) \neq \alpha$, so by Proposition \ref{prop_color_types} $\alpha \in S_X^b \cup S_X^p$.
  With that there is no curve class $[C_\alpha] \in N_1(G/P^-)$, such that $\iota_*[C_\alpha] = 2[C_D]$ for 
  some $D$ and therefore $N_1(X) \subset \Ima(\iota_*)$.
  For any color $D$ there are exatcly two simple roots $\alpha_D$ and $\beta_D$ moving $D$.
  Then we can pick $\eta_{G/P^-}$ effective, such that $c_{\alpha_D} = c_D$ and $c_{\beta_D} = 0$.
  $\qed$
\end{proof}

\begin{lemma}\label{lemma_irreducible_rep_meeting_open_orbit}
  Let $\eta \in N_1(X)$ be a movable curve class on $X$. Let 
  \[
    M_{0,0}^{irr}(X, \eta) = \{(f:C \to X)\ |\ C\text{ irreducible}\}
  \]
  be the subset of $\overline{M}_{0,0}(X,\eta)$ of stable maps, having an irreducible source.
  Suppose that $M_{0,0}^{irr}(X, \eta)$ is nonempty, then 
  $M_{0,0}^{\circ}(X, \eta)$ is nonempty as well.
\end{lemma}

\begin{proof}
  First note that for any morphism $f : C \to X$ with $C$ irreducible
  there exists a minimal $G$-stable variety $X_I = \cap_{i \in I}X_i$ with $I \subset \{1, \dots, \rk(X)\}$,
  such that $f(C) \subset X_I$. Let $\iota_I: X_I \to X$ be the inclusion, let $\eta_I \in N_1(X_I)$ be 
  any curve class, such that $(\iota_I)_* \eta_I = \eta$
  then we get an inclusion
  \begin{align*}
    M_{0,0}^\circ(X_I, \eta_I) &\to M_{0,0}^{irr}(X, \eta) \\
    (f: C \to X) &\mapsto \iota_I \circ f
  \end{align*}
  and therefore 
  \[
    M_{0,0}^{irr}(X, \eta) = \bigcup_{I \subset \{1, \dots, \rk(X)\}} \bigcup_{(\iota_I)_*\eta_I = \eta} M_{0,0}^\circ(X_I, \eta_I),
  \]
  so we have 
  \[
    \dim M_{0,0}^{irr}(X, \eta) = \max_{I \subset \{1, \dots, \rk(X)\}} \max_{(\iota_I)_*\eta_I = \eta}\dim M_{0,0}^{\circ}(X_I, \eta_I).
  \]

  For any $I \subset \{1, \dots, \rk(X)\}$ by Corollary \ref{corollary_moduli_circ_exact_dimension} we have 
  \[
    \dim M_{0,0}^\circ(X_I, \eta_I) = - \langle K_{X_I}, \eta_I \rangle + \dim X_I - 3.
  \]
  We apply adjunction and projection formulas to $\langle K_{X_I}, \eta_I \rangle$:
  \[
    - \langle K_{X_I}, \eta_I \rangle = - \langle \restr{(K_X + \sum_{i \in I}X_i)}{X_I}, \eta_I \rangle 
      = - \langle K_X + \sum_{i \in I}X_i, \eta \rangle.
  \]
  Also $X_I$ is an intersection of transversal divisors, so $\dim X_I = \dim X - |I|$ and we get
  \[
    \dim M_{0,0}^\circ(X_I, \eta_I) = - \langle K_X, \eta \rangle + \dim X - 3 - |I| - \sum_{i \in I}\langle X_i, \eta \rangle.
  \]
  
  On the other hand by Proposition \ref{prop_moduli_space_dimension} we have 
  \[
    \dim M_{0,0}^{irr}(X, \eta) \geq - \langle K_X, \eta \rangle + \dim X - 3.  
  \]
  
  Let us compare the estimates for $\dim M_{0,0}^{irr}(X, \eta)$:
  \begin{align*}
    \dim M_{0,0}^{irr}(X, \eta) - \dim M_{0,0}^\circ (X_I, \eta_I)\\
      = |I| + \sum_{i \in I}\langle X_i, \eta \rangle
  \end{align*}
  As $\eta$ is movable, by Proposition \ref{movable_effective_cones_duality}
  it lies in the dual cone to the cone of effective divisors, so $\langle X_i, \eta \rangle \geq 0$.
  With that we have that 
  \[
    \dim M_{0,0}^{irr}(X, \eta) - \dim M_{0,0}^\circ (X_I, \eta_I) = 0  
  \]
  if and only if $|I| = 0$ and so $M_{0,0}^\circ (X, \eta)$ is nonempty.
  $\qed$
\end{proof}

\begin{corollary}\label{corollary_wond_g_stable_var_curve_rep}
  Let $\eta \in N_1(X)\setminus \{0\}$ be a movable class on $X$.
  Then $M_{0,0}^\circ(X, 2\eta)$ is nonempty.
  If $X$ a wonderful group compactification, then $M_{0,0}^\circ(X, \eta)$ is nonempty as well.
\end{corollary}

\begin{proof}
  By Lemma \ref{lemma_irreducible_rep_meeting_open_orbit} we only have to show 
  that $M_{0,0}^{irr}(X, 2\eta)$ is nonempty.

  A movable class is nonnegative on colors of $X$, so
  by Lemma \ref{lemma_positivity_of_closed_orbit_class} there exists an effective curve class $\eta_{G/P^-} \in N_1(G/P^-)$,
  such that $\iota_* \eta_{G/P^-} = 2\eta$, where $\iota: G/P^- \to X$ is the inclusion of the closed $G$-orbit.
  Note that if $X$ is a wonderful group compactification, then there exists
  $\eta_{G/P^-}$ effective such that $\iota_* \eta_{G/P^-} = \eta$.
  Let $\eta' = \iota_* \eta_{G/P^-}$.

  By Proposition \ref{prop_homogeneous_variety_mod_space_irreducible} we have that $M_{0,0}^\circ(G/P^-, \eta_{G/P^-})$
  is nonempty. It embeds into $M_{0,0}^{irr}(X, \eta')$ via composition with $\iota$, therefore
  $M_{0,0}^{irr}(X, \eta') \neq \emptyset$ and the result follows.
  $\qed$
\end{proof}

\begin{theorem}\label{thm_wond_mod_space_reducible}
  Let $\eta \in N_1(X)\setminus \{0\}$ be a curve class such that $M_{0,0}^\circ(X, \eta) \neq \emptyset$.
  Assume there exists a decomposition $\eta = \eta_1 + \eta_2$, such that 
  $\eta_1$ and $\eta_2$ are nontrivial effective curve classes on $X$
  and a boundary divisor $D$ with
  \[
    \langle D, \eta_2 \rangle \leq -2.
  \]

  Then the moduli space $\overline{M}_{0,0}(X, \eta)$ is reducible.
\end{theorem}

\begin{proof}
  Let $I_1 = \{ i \in \{1, \dots, l\}\ |\ \langle X_i, \eta_1 \rangle < 0\}$ 
  and $I_2 = \{ i \in \{1, \dots, l\}\ |\ \langle X_i, \eta_2 \rangle < 0\}$. 
  By Lemma \ref{lemma_positivity_of_curve_meeting_oo} the class $\eta$ is nonnegative on boundary divisors, 
  so $I_1 \cap I_2$ is empty. Write $X_1 = X_{I_1}$ and $X_2 = X_{I_2}$. 
  Then $\eta_1$ (resp. $\eta_2$) is nonnegative on boundary divisors and colors of $X_1$ (resp. $X_2$), 
  so in particular by Proposition \ref{movable_effective_cones_duality} $\eta_1$ (resp. $\eta_2$) is 
  movable on $X_1$ (resp. $X_2$). With that there exist representatives of both $\eta_1$ and $\eta_2$, 
  that meet in at least one point of $X_1 \cap X_2$, so the moduli space 
  $\overline{M}_{0,1}(X_1, \eta_1) \times_{X} \overline{M}_{0,1}(X_2,\eta_2)$ is nonempty. 
  The class of the union of these representatives is $\eta_1 + \eta_2 = \eta$, 
  and in general we get a map 
  \[
    \overline{M}_{0,1}(X_1, \eta_1) \times_{X} \overline{M}_{0,1}(X_2,\eta_2) \to \overline{M}_{0,0}(X, \eta),
  \]
  unmarking the point on both $\overline{M}_{0,1}(X_1, \eta_1)$ and $\overline{M}_{0,1}(X_2,\eta_2)$.

  By our assumption the space $M_{0,0}^\circ(X, \eta)$ is nonempty and by Corollary \ref{corollary_moduli_circ_exact_dimension}
  it is of dimension 
  \[
    -\langle K_X, \eta \rangle + \dim X - 3.
  \]
  We wish to estimate the dimension of $\overline{M}_{0,1}(X_1, \eta_1) \times_{X} \overline{M}_{0,1}(X_2,\eta_2)$ from below 
  to show that $\overline{M}_{0,1}(X_1, \eta_1) \times_{X} \overline{M}_{0,1}(X_2,\eta_2)$ cannot lie in the same irreducible
  components as $M_{0,0}^\circ(X, \eta)$.

  To do so we consider the projection:
  \[
    p_1: \overline{M}_{0,1}(X_1, \eta_1) \times_{X} \overline{M}_{0,1}(X_2,\eta_2) \to \overline{M}_{0,1}(X_1, \eta_1).
  \]
  We have that
  \[
    \Ima(p_1) = \{(f, x_1) \in \overline{M}_{0,1}(X_1, \eta_1)\ |\ f(x_1) \in X_1 \cap X_2\},
  \]
  because $\eta_2$ is movable, so there is a curve passing through any point of $X_2$.
  Let us consider the map 
  \begin{align*}
    ev_{X_1}:  \overline{M}_{0,1}(X_1, \eta_1) &\to X_1\\
    (f, x_1) &\mapsto f(x_1),
  \end{align*}
  then $\Ima(p_1) = ev_{X_1}^{-1}(X_2\cap X_1)$. 
  From that we get that 
  \[
    \mathrm{codim}_{\overline{M}_{0,1}(X_1, \eta_1)}\Ima(p_1) \leq \mathrm{codim}_{X_1} X_1\cap X_2,  
  \]
  or equivalently:
  \[
    \dim \Ima(p_1) \geq \dim \overline{M}_{0,1}(X_1, \eta_1) - \dim X_1 + \dim X_1\cap X_2.
  \]
  By Proposition \ref{prop_moduli_space_dimension} we get that 
  \[
    \dim \overline{M}_{0,1}(X_1, \eta_1) \geq -\langle K_{X_1}, \eta_1 \rangle + \dim X_1 - 3 + 1,
  \]
  so 
  \[
    \dim \Ima(p_1) \geq -\langle K_{X_1}, \eta_1 \rangle + \dim X_1\cap X_2 - 2.
  \]

  Let $(f, x_1) \in \Ima(p_1)$, we wish to estimate the dimension of $p_1^{-1}(f, x_1)$.
  We have that
  \begin{align*}
    p_1^{-1}(f, x_1) = &\{(f, x_1) \times_X \overline{M}_{0,1}(X_2, \eta_2)\} \\
      = \{\big((f, x_1), &(g, x_2)\big)\ |\  (g, x_2)\in \overline{M}_{0,1}(X_2, \eta_2) \text{ with } f(x_1) = g(x_2) \}.
  \end{align*}
  We consider the second projection
  \[
    p_2: \overline{M}_{0,1}(X_1, \eta_1) \times_{X} \overline{M}_{0,1}(X_2,\eta_2) \to \overline{M}_{0,1}(X_2, \eta_2),
  \]
  then 
  \[
    p_2(p_1^{-1}(f, x_1)) = \{(g, x_2) \in \overline{M}_{0,1}(X_2, \eta_2)\ |\ f(x_1) = g(x_2) \}.
  \]
  On the other hand we have the map 
  \begin{align*}
    ev_{X_2}:  \overline{M}_{0,1}(X_2, \eta_2) &\to X_2\\
    (g, x_2) &\mapsto g(x_2),
  \end{align*}
  and 
  \[
    ev_{X_2}^{-1}(f(x_1)) = \{(g, x_2) \in \overline{M}_{0,1}(X_2, \eta_2)\ |\ g(x_2) = f(x_1)\} = p_2(p_1^{-1}(f, x_1)).
  \]
  With that we get an estimate for the dimension of the fiber $p_2(p_1^{-1}(f, x_1))$: 
  \[
    \mathrm{codim}_{\overline{M}_{0,1}(X_2, \eta_2)} p_2(p_1^{-1}(f, x_1)) \leq \mathrm{codim}_{X_2} \{x_1\},
  \]
  or equivalently:
  \[
    \dim p_2(p_1^{-1}(f, x_1)) \geq \dim \overline{M}_{0,1}(X_2, \eta_2) - \dim X_2,
  \]
  and again by Proposition \ref{prop_moduli_space_dimension} we have:
  \[
    \dim p_2(p_1^{-1}(f, x_1)) \geq -\langle K_{X_2}, \eta_2 \rangle - 2.
  \]

  To combine these two estimates to get a lower bound for 
  $\dim \overline{M}_{0,1}(X_1, \eta_1) \times_{X} \overline{M}_{0,1}(X_2,\eta_2)$ we pick an irreducible component 
  $M_{im}$ of $\Ima(p_1)$ of maximal dimension and an irreducible component $M_s$ of 
  $p_1^{-1}(M_{im})$, dominating $M_{im}$ also of maximal dimension. 
  We consider the restriction $\restr{p_1}{M_s}$ to get
  \[
    \dim M_s = \dim M_{im} + \dim F,
  \]
  where $F$ is a general fiber of $\restr{p_1}{M_s}$. We have an estimate for the dimension of any fiber of $p_1$, so 
  in total we have
  \begin{align*}
    \dim \overline{M}_{0,1}&(X_1, \eta_1) \times_{X} \overline{M}_{0,1}(X_2,\eta_2) \geq \dim M_{im} + \dim F\\
      \geq \dim &\Ima(p_1) + \dim M_{0,1}(X_2, \eta_2) - \dim X_2 \\
        &\geq -\langle K_{X_1}, \eta_1 \rangle + \dim X_1\cap X_2 - 2 -\langle K_{X_2}, \eta_2 \rangle - 2 \\
        &= \sum_{j = 1,2}(-\langle K_{X_j}, \eta_j \rangle) + \dim X_1\cap X_2 - 4.
  \end{align*}
  By adjunction formula $K_{X_j} = K_X - \sum_{i \in I_j} X_i$
  for $j = 1,2$ and we get
  \begin{align*}
    \dim& \overline{M}_{0,1}(X_1, \eta_1) \times_{X} \overline{M}_{0,1}(X_2,\eta_2) \\
    \geq & \sum_{j = 1,2}\big(-\langle K_{X}, \eta_j \rangle - \sum_{i \in I_j} \langle X_i, \eta_j \rangle\big)
      + \dim X_1 \cap X_2 - 4 \\
    = & -\langle K_X, \eta_1 + \eta_2\rangle - \sum_{j = 1,2}\sum_{i \in I_j}\langle X_i, \eta_j \rangle
      + \dim X_1 \cap X_2 - 4 \\
    = & -\langle K_X, \eta\rangle - \sum_{j = 1,2}\sum_{i \in I_j}\langle X_i, \eta_j \rangle
    + \dim X_1 \cap X_2 - 4 
  \end{align*}
  As $X_1 \cap X_2 = \cap_{i \in I_1 \cup I_2}X_i$ and all $X_i$'s intersect transversally,
  we have $\dim X_1\cap X_2 = \dim X - \mathrm{codim} X_1 \cap X_2 = \dim X - |I_1 \cup I_2|$.
  Applying that $I_1$ and $I_2$ are disjoint yields the final form of the lower bound:
  \begin{align*}
    \dim \overline{M}_{0,1}&(X_1, \eta_1) \times_{X} \overline{M}_{0,1}(X_2,\eta_2) \geq \\
    -\langle K_X,& \eta\rangle  - \sum_{j = 1,2}\sum_{i \in I_j}\langle X_i, \eta_j \rangle
      + \dim X - (|I_1| + |I_2|) - 4
  \end{align*}

  Any curve, that is a union of curves of classes $\eta_1$ and $\eta_2$ will not lie 
  in $M_{0,0}^\circ (X, \eta)$, because it is reducible. We compare 
  the dimension of $M_{0,0}^\circ (X, \eta)$ given by Corollary \ref{corollary_moduli_circ_exact_dimension}
  and the lower bound on the dimension of $\overline{M}_{0,1}(X_1, \eta_1) \times_{X} \overline{M}_{0,1}(X_2,\eta_2)$:
  \begin{align*}
    \dim M_{0,0}^\circ&(X, \eta) - \dim \overline{M}_{0,1}(X_1, \eta_1) \times_{X} \overline{M}_{0,1}(X_2,\eta_2) \\
    \leq& -\langle K_X,\eta \rangle + \dim X - 3 \\
    &- \big(-\langle K_X, \eta\rangle  - \sum_{j = 1,2}\sum_{i \in I_j}\langle X_i, \eta_j \rangle
      + \dim X - |I_1| - |I_2| - 4\big) \\
    =& 1 + |I_1| + |I_2| + \sum_{i \in I_1}\langle X_i, \eta_1 \rangle 
      + \sum_{i \in I_2}\langle X_i, \eta_2 \rangle
  \end{align*}
  By definition of $I_j$ we have that $\langle X_i, \eta_j \rangle < 0$ for all $i \in I_j$, so
  \[
    |I_j| + \sum_{i \in I_j}\langle X_i, \eta_j \rangle \leq 0, 
  \]
  and so
  \[
    \dim M_{0,0}^\circ(X, \eta) - \dim \overline{M}_{0,1}(X_1, \eta_1) \times_{X} \overline{M}_{0,1}(X_2,\eta_2) \leq 1.
  \]

  We apply the assumption on $\eta$ and its decomposition $\eta = \eta_1 + \eta_2$,
  that there exists $i_0 \in \{1, \dots, l\}$ such that $D = X_{i_0}$ is a boundary divisor with
  \[
    \langle D, \eta_2 \rangle < -2,
  \]
  then 
  \[
    |I_2| + \sum_{i \in I_2}\langle X_i, \eta_2 \rangle \leq -1,
  \]
  so 
  \[
    \dim M_{0,0}^\circ(X, \eta) - \dim \overline{M}_{0,1}(X_1, \eta_1) \times_{X} \overline{M}_{0,1}(X_2,\eta_2) \leq 0.
  \]
  The general curve in $\overline{M}_{0,1}(X_1, \eta_1) \times_{X} \overline{M}_{0,1}(X_2,\eta_2)$ is reducible by construction, 
  whereas the general curve in $M_{0,0}^\circ(X, \eta)$ is irreducible, therefore
  $M_{0,0}^\circ(X, \eta)$ cannot lie in $\overline{M}_{0,1}(X_1, \eta_1) \times_{X} \overline{M}_{0,1}(X_2,\eta_2)$.
  Altogether the whole space $\overline{M}_{0,0}(X,\eta)$ has at least two irreducible components 
  and is reducible.
  $\qed$
\end{proof}

The condition of the theorem is satisfied in a very big class of wonderful varieties,
namely the wonderful group compactifications.

\begin{corollary}\label{corollary_wond_group_comp_mod_space_reducible}
  Let $X$ be a wonderful group compactification of a simple group $G$ of adjoint type with $\rk(G) \geq 3$.
  Then $\overline{M}_{0,0}(X, \eta)$ is reducible for any movable class of curves $\eta \in N_1(X)\setminus \{0\}$.
\end{corollary}

\begin{proof}
  By Corollary \ref{corollary_wond_g_stable_var_curve_rep} $M_{0,0}^\circ(X, \eta) \neq \emptyset$.
  We are looking for a boundary divisor $D$ and a decomposition $\eta = \eta_1 + \eta_2$, 
  such that $\langle D, \eta_2 \rangle \leq -2$.

  The variety $X$ is the wonderful compactification of the symmetric space $G\times G/ (G\times G)^\sigma$,
  where $\sigma(x, y) = (y,x)$ is an involution on $G \times G$.
  The root system of $G\times G$ is the product of two copies of the root systems of $G$.
  Let $\Phi_G$ be the first copy and $\Phi_G'$ the second one.
  Then the spherical roots of $X$ are of the form $\alpha_i + \alpha_i'$, where 
  $\alpha_i \in \Phi_G$ and $\alpha_i' \in \Phi_G'$ are simple roots, such that $\sigma(\alpha_i) = - \alpha_i'$.
  Denote by $S_G$ (resp. $S_G'$) the sets of simple roots of $\Phi_G$ (resp. $\Phi_G'$).
  Note that $\alpha_i$ and $\alpha_j'$ are orthogonal for any $i, j \in \{1, \dots, \rk(G)\}$,
  and with that we can conclude that $X$ is nonexceptional and by Proposition
  \ref{prop_picard_grp_as_weight_lattice} we have that $\rk \Pic(X) = \rk (G)$.
  
  For any simple root $\alpha \in S$ we have that neither $\alpha$, nor $2\alpha$ is a spherical root of $X$, 
  so by Proposition \ref{prop_color_types}, we have that $\alpha \in S_X^b \cup S_X^p$.
  By Lemma \ref{lemma_wond_color_dual_curve_description} 
  any curve class $\eta$ is of the form $\sum_{D \in \Delta_X} c_D [C_D]$.
  For any simple root $\alpha \in S_G$, there exists a divisor $D_\alpha \in \Delta_X(\alpha)$,
  moved by $\alpha$.
  Note that in the case of group compactification $D_\alpha$ is the closure of the Bruhat cell of $G$,
  so in particular for $\alpha \neq \beta$, we have $D_\alpha \neq D_\beta$.
  This means that the colors $\{ D_\alpha \ |\ \alpha \in S_G\}$ are linearly independent and therefore 
  generate $\Pic(X)$, since 
  \[
    \Pic(X) \simeq \mathbb{Z}^{\rk(G)} \simeq \bigoplus_{\alpha \in S_G}\mathbb{Z}D_\alpha
  \]
  are all of the same rank and $\Pic(X)$ is generated by colors of $X$.
  With that we can write 
  \[
    \eta = \sum_{i = 1}^{\rk(G)} c_i [C_{D_{\alpha_i}}].
  \]
  
  As $G$ is simple and of rank $\geq 3$, there exist simple roots of $G$,
  which are connected to more that one root in the Dynking diagram of $G$.
  We will call such roots nonextremal in this proof.

  We wish to show that there exists an nonextremal simple root $\alpha_{i_0} \in S$, such that $c_{i_0} \neq 0$.
  As $\eta \neq 0$, there exists at least one $1 \leq j \leq \rk(G)$, such that $c_j \neq 0$.
  If $\alpha_j$ is nonextremal, we are done and $i_0 = j$.
  Suppose that $\alpha_j$ is extremal, we claim that $c_{i_0} \neq 0$ for $\alpha_{i_0}$ nonextremal, adjacent to $\alpha_j$.
  By Corollary \ref{corollary_wond_var_boundary_divisors_intersection_with_curves} and Proposition \ref{prop_wond_var_rho_value} we have
  for any boundary divisor $X_{i'}$
  \[
    \langle X_{i'}, \eta \rangle = \sum_{i = 1}^{\rk(G)}  c_i\langle \rho_X(D_{\alpha_i}), \alpha_{i'} - \sigma({\alpha_{i'})}\rangle
      = \sum_{i = 1}^{\rk(G)}  c_i \alpha_i^\vee (\alpha_{i'} - \sigma({\alpha_{i'}}))
  \]
  We have that $\alpha_i \in S_X^b$ for every $1 \leq i \leq \rk(G)$ and that $\alpha_i$ is orthogonal to $\sigma(\alpha_{i'})$, as
  $\sigma(\alpha_{i'})$ lies in $S_{G'}$. 
  With that and the assumption that $\eta$ is movable, we get for the boundary divisor $X_{i_0}$:
  \[
    0 \leq \langle X_{i_0}, \eta \rangle = \sum_{i = 1}^{\rk(G)}  c_i\alpha_i^\vee(\alpha_{i_0}) 
      \leq c_j\alpha_j^\vee(\alpha_{i_0}) + c_{i_0}\alpha_{i_0}^\vee(\alpha_{i_0})
  \]
  and $c_j\alpha_j^\vee(\alpha_{i_0}) < 0$, so $c_{i_0}$ has to be positive.

  Let $\eta_1 = c_{i_0}[D_{\alpha_{i_0}}]$ and $\eta_2 = \eta - c_{i_0}[D_{\alpha_{i_0}}]$, 
  then $\eta = \eta_1 + \eta_2$.
  For all $k$ with $\alpha_k$ adjacent to $\alpha_{i_0}$ we have
  $\langle X_{k}, \eta_1 \rangle = c_{i_0}\alpha_{i_0}^\vee(\alpha_{k}) < 0$, thus
  we have that $\langle X_{k}, \eta_2 \rangle > 0$, because 
  $\eta$ is nonnegative on effective divisors. This means that $c_k$
  is positive for all $k$, such that $\alpha_k$ is adjacent to $\alpha_{i_0}$

  Since $i_0$ is nonextremal, we have
  \[
    \langle X_{i_0}, \eta_2 \rangle = \sum_{i \neq i_0} c_i\alpha_i^\vee(\alpha_{i_0}) 
      \leq \sum_{\alpha_k \text{ adj. }\alpha_{i_0}}c_{k} \alpha_{k}^\vee(\alpha_{i_0}) \leq -2,
  \]
  and
  both $\eta_1$ and $\eta_2$ are effective, as they are nonnegative linear combinations of simple coroots.
  The divisor $X_{i_0}$ is the $D$ we are looking for.
  The result follows by Theorem \ref{thm_wond_mod_space_reducible}.
  $\qed$
\end{proof}

\begin{corollary}
  For $\eta$ the class of the VMRT of a wonderful compactification of a simple group of adjoint type 
  the moduli space $\overline{M}_{0,0}(X, \eta)$ is reducible.
\end{corollary}

\begin{example}\label{example_pgl4}
  Let $X = \overline{\mathrm{PGL}_4(\mathbb{C})}$ be the wonderful group compactification. In other 
  words it is the wonderful embedding of the symmetric 
  space $G\times G/G$ for $G = \mathrm{PGL}_4(\mathbb{C})$.
  Let 
  \[
    \{\alpha_1, \alpha_2, \alpha_3 ,\beta_1, \beta_2, \beta_3\}
  \]
  be the usual simple roots of $\mathrm{PGL}_4(\mathbb{C})\times \mathrm{PGL}_4(\mathbb{C})$ with $\beta_1, \beta_2, \beta_3$
  denoting the simple roots of the second term of the product.
  The involution $\sigma: G\to G$ given by $(x,y) \mapsto (y,x)$ acts on the simple roots as follows:
  \[
    \sigma(\alpha_i) = -\beta_i.
  \]
  For $i \in \{1, 2, 3\}$ there is a color $D_i$ moved by roots $\alpha_i, \beta_i$.
  A boundary divisor $X_i$ corresponds to the spherical root 
  \[
    \bar{\alpha}_i = \alpha_i - \sigma(\alpha_i) = \alpha_i + \beta_i
  \]
  by calculating its presentation as a linear combination of colors using \ref{prop_wond_var_boundary_divisors_as_colors}.

  Take the class $\eta$ to be the unique class of minimal rational curves,
  described by Michel Brion and Baohua Fu in \cite{brion_fu_mrc}. This class in type $A$ 
  is given by the coroot of the highest root.
  In case of $\mathrm{PGL}_4(\mathbb{C})$ it is given by $\alpha_1^\vee + \alpha_2^\vee + \alpha_3^\vee$.
  
  Then for any color $D_i$ of $X$ we have that $\langle D_i, \eta \rangle = 1$, whereas for
  the boundary divisors we have the following values: $\langle X_1, \eta \rangle = \langle X_3, \eta \rangle = 1$
  and $\langle X_2, \eta \rangle = 0$.

  We apply our result for $\bar{\alpha}_{i_0} = \bar{\alpha}_2 = \alpha_2 + \beta_2$ and we get a decomposition
  $\eta_1 = \alpha_2^\vee$ and $\eta_2 = \alpha_1^\vee + \alpha_3^\vee$.
  Now we calculate the respective intersections of these classes with boundary divisors and we get
  \begin{equation*}
    \begin{split}
      \langle X_1, \eta_1 \rangle &= -1\\
      \langle X_2, \eta_1 \rangle &= 2\\
      \langle X_3, \eta_1 \rangle &= -1
    \end{split}
    \qquad\qquad
    \begin{split}
      \langle X_1, \eta_2 \rangle &= 2\\
      \langle X_2, \eta_2 \rangle &= -2\\
      \langle X_3, \eta_2 \rangle &= 2
    \end{split}
  \end{equation*}

  As required in the statement of the theorem we have that $\langle X_2, \eta_2 \rangle \leq -2$.
  A similar approach works for any wonderful compactification of a group $G$, such that 
  the entries in at least one row of the Cartan matrix of $G$ add up to a nonpositive value.
\end{example}

\section{Irreducibility of the special component}\label{sec_spec_comp_irred}

Let $X$ be the wonderful compactification of a semisimple algebraic group $G$ of adjoint type.
In terms of notation of Section \ref{sec_wonderful_compactifications_sym_space}, $X$ is the compactification of
the symmetric space $G\times G/G$, where the involution is defined as
\begin{align*}
  \sigma: G\times G & \to G\times G \\  
  (g,g') & \mapsto (g',g)
\end{align*}
We write $\mathcal{G} = G\times G$.
Denote by $R$ the root system of $\mathcal{G}$ and by $S$ its basis of simple roots 
\[
  S = \{\alpha_1, \dots, \alpha_l, \beta_1, \dots \beta_l\},
\] 
such that $\sigma(\alpha_i) = -\beta_i$. The corresponding Borel subgroup is 
of the form $B\times B^-$, where $B \subset G$ is a Borel subgroup associated with the 
basis of simple roots $\alpha_1, \dots, \alpha_l$.
Let $\omega_i$ (resp. $\zeta_i$) be the fundamental weight, corresponding to $\alpha_i$ (resp. $\beta_i$)
for $1 \leq i \leq l$.
Denote by 
\[
  \Sigma_X = \{\bar{\alpha}_1,\dots, \bar{\alpha}_l\}
\]
the set of spherical roots of $X$ 
where $\bar{\alpha}_i = \alpha_i - \sigma(\alpha_i) = \alpha_i + \beta_i$.
To every $\bar{\alpha}_i$ there exists a corresponding boundary divisor $X_i$,
such that $\bar{\alpha}_i$ is the weight of the $T$-action on $T_zX/T_zX_i$, 
where $z$ is the unique $B^-$-fixed point of the closed orbit $Y \subset X$.

\begin{lemma}\label{lemma_wond_group_comp_colors}
  The closed orbit $Y$ of $X$ is isomorphic to $G\times G/B\times B^-$.
  The simple roots of $\mathcal{G}$ are all of type $(b)$, i.e. $S_X^b = S_\mathcal{G}$, and
  $\Delta_X(\alpha_i) = \Delta_X(\beta_i) = \{D_i\}$ for every $1 \leq i \leq l$.
  Moreover $\Delta_X = \{D_1, \dots, D_l\}$.
\end{lemma}

\begin{proof}
  For any simple root $\alpha \in S$ we have that $\alpha, 2\alpha \notin \Sigma_X$.
  Moreover, there are no simple roots stable under $\sigma$, so $S_X^p = \emptyset$, 
  so by Lemma \ref{lemma_closed_orbit_form} $Y = \mathcal{G}/B\times B^-$.
  By Proposition \ref{prop_color_types} $S = S_X^b$, in other words for any spherical root $\bar{\alpha}_i$ 
  we have $\langle \rho_X(D), \bar{\alpha}_i \rangle = \alpha^\vee(\bar{\alpha}_i)$, where $\alpha$
  is some simple root of $G$ with $\{D\} = \Delta_X(\alpha)$.

  On the other hand any color $D$ of $X$ is the closure of a codimension one Bruhat cell of $G$.
  This means that $D = \overline{B s_{\alpha_i} B^-}$ for some $\alpha_i$. We have that
  $P_{\alpha_i} B s_{\alpha_i} B^-$ and $B s_{\alpha_i} B^- P_{\alpha_i}^-$ are both open in $G$.
  Let $\mathcal{P}_{\alpha_i}$ denote the minimal parabolic of $G\times G$ given by $P_{\alpha_i} \times \{id\}$
  and $\mathcal{P}_{\beta_i}$ denote the minimal parabolic of $G\times G$ given by $\{id\} \times P_{\alpha_i}^-$.
  Then we have that 
  \[
    P_{\alpha}B s_{\alpha_i} B^- = \mathcal{P}_{\alpha_i}  B s_{\alpha_i} B^-
  \]
  and 
  \[
    B s_{\alpha_i} B^- P_{\alpha_i}^- = \mathcal{P}_{\beta_i}  B s_{\alpha_i} B^-.
  \]
  This means that every $D_i$ is moved by exactly two simple roots of $\mathcal{G}$, namely $\alpha_i$ and $\beta_i$.

  In total we get that $\Delta_X = \bigcup_{i = 1}^l \Delta_X(\alpha_i) \cup \Delta_X(\beta_i) = \{D_1, \dots , D_l\}$.
  $\qed$
\end{proof}


Let $V \simeq U \times \mathbb{A}^l$ be the affine subset introduced in Remark \ref{remark_symmtric_roots_not_moving_colors}.
As $X$ is the wonderful embedding of $\mathcal{G}/G$, $U$ is the unipotent radical of the 
Borel subgroup $B \times B^-$, 
where $B$ is the Borel of $G$ corresponding to the choice of the base of simple roots and $\mathbb{A}^l$
is the affine space, on which the maximal torus chosen in \ref{lemma_dcp_torus_choice} acts via
$$t(x_1, \dots x_l) = (\bar{\alpha}_1(t)x_1, \dots, \bar{\alpha}_l(t)x_l)$$
for $(x_1, \dots x_l) \in \mathbb{A}^l$.

Let $\lambda$ be the $1$-parameter subgroup given by the fundamental coweight 
$\omega_{i_0}$ of the simple root $\alpha_{i_0}$.
Then for any simple root $\alpha$ of $\mathcal{G}$ we have that 
\[
  \langle \lambda, \alpha \rangle = \begin{cases}
    1 &\text{if }\alpha = \alpha_{i_0},\\
    0 &\text{else}.
  \end{cases}
\]

Define the family of automorphisms on $X$ by $\varphi_t: x \to \lambda(t)x$ for $t \in \mathbb{C}^*$ 
and consider the limit map as $t \to 0$:
\begin{align*}
  \varphi: X &\dashrightarrow X \\ 
  x &\mapsto \lim_{t \to 0} \lambda(t)x
\end{align*}

We will show in the next Proposition 
that this map is defined on $V$, which is an open subset of $X$ intersecting every $G$-orbit of $X$.
In particular this means that for any $G$-stable subvariety $X_I$ for some $I \subset \{1, \dots, l\}$
there exists a restriction $\restr{\varphi}{X_I}: X_I\dashrightarrow X_I$.

\begin{proposition}\label{prop_limit_map_image_properties}
  Let $X_I \subset X$ be a closed $\mathcal{G}$-stable subvariety. 
  The rational map $\varphi$ restricts to a rational map
  \[
    \restr{\varphi}{X_I}: X_I \dashrightarrow X_I \cap X_{i_0},
  \]
  which has indeterminacy locus of codimension at least two and the open 
  $\mathcal{G}$-orbit of $X_I$ is mapped inside the open $\mathcal{G}$-orbit of $X_{i_0}\cap X_I$.
\end{proposition}

\begin{proof}
  First we note that the indeterminacy locus of any rational morphism 
  between smooth projective varieties is of codimension 
  at least $2$. 
  
  The map $\varphi$ is defined everywhere on $V \simeq U\times \mathbb{A}^l$. 
  Indeed, for $(u,x)$, where $u \in U$ and $x = (x_1, \dots x_l)\in \mathbb{A}^l$ we have 
  \[
    \lambda(t)(u,x) = (\lambda(t)u\lambda(t)^{-1}, (t^{\langle \lambda, \bar{\alpha}_i \rangle}x_i)_{1\leq i\leq l}).
  \]
  Compute $\langle \lambda, \bar{\alpha}_i \rangle = \langle \omega_{i_0}, \alpha_i - \beta_i \rangle = \delta_{i_0 i}$ and 
  take the limit on the $\mathbb{A}^l$ part to get 
  \[
  \lim_{t\to 0 }(t^{\langle \lambda, \bar{\alpha}_i \rangle}x_i)_{1 \leq i \leq l} = (x_1, \dots, x_{i_0-1}, 0, x_{i_0+1}, \dots  x_l).
  \]
  As for the part $U$, use that $U = \prod_{\alpha \in R^+}U_\alpha$ with every $U_\alpha \simeq \mathbb{C}$ and 
  for some $u_\alpha \in U_\alpha$ identified with the complex numbers
  \[
    \lambda(t)u_\alpha \lambda(t)^{-1} = t^{\langle\lambda,\alpha \rangle}u_\alpha,
  \]
  so for any $\alpha \in R^+$ we have
  \[
    \lim_{t\to 0} \lambda(t)u_\alpha \lambda(t)^{-1} = \begin{cases}
      1_G &\text{if }\alpha \geq \alpha_{i_0},\\
      u_\alpha &\text{else}.
    \end{cases}
  \]  


  Recall that $V\cap X_{i_0} = \{ (u,x) \in V\ |\ x_{i_0} = 0\}$ and $\lambda(t)$ for $t \in \mathbb{C}^*$ acts via 
  \[
    \lambda(t) (u,x) = (\lambda(t)u \lambda(t)^{-1}, x_1, \dots, tx_{i_0}, \dots, x_l)
  \]
  so it follows immediately that the image of $\varphi$ lies in $X_{i_0}$.  
  Next, $X_I$ is $\mathcal{G}$-stable and $\varphi$ is defined via the $\mathcal{G}$-action, so it 
  restricts to $X_I$:
  \[
    \restr{\varphi}{X_I}: X_I \dashrightarrow X_I\cap X_{i_0}.
  \]
  
  The intersection of $V$ with the open $\mathcal{G}$-orbit of $X_{i_0}$ is given by 
  \[
    \{ (u,x) \in V \cap X_{i_0}\ |\ x_i \neq 0,\ i\notin \{i_0\}\}, 
  \]
  and the intersection of $V$ with the open $\mathcal{G}$-orbit of $X_I$ is given by 
  \[
    \{ (u,x) \in V \cap X_I\ |\ x_i \neq 0,\ i\notin I\}, 
  \]
  since $X_I = \cap_{i \in I}X_i$ for some $I \subset \{1, \dots l\}$.
  Therefore any point in the intersection of $V$ with the open $\mathcal{G}$-orbit of $X_I$ 
  is mapped via $\varphi$ to the open $\mathcal{G}$-orbit of $X_{i_0}\cap X_I$.
  $\qed$
\end{proof}

By Proposition \ref{prop_boundary_picard_exact_sequence} we have
\[
  \rk\Pic(X_{i_0}) = \rk \Pic(\mathcal{G}/\mathcal{P}_{S_{\{i_0\}}}^-) + \rk \Pic(X^{\{i_0\}}).
\]
Recall that $S_{\{i_0\}} = S \setminus \{\alpha, \sigma(\alpha)\}$, so $rk \Pic(\mathcal{G}/\mathcal{P}_{S_{\{i_0\}}}^-) = 2$
and that $X^{\{i_0\}}$ is the wonderful group compactification of the Levi quotient of $L_{S_{\{i_0\}}}$,
so $\rk \Pic(X^{\{i_0\}}) = \rk\Pic(X) - 1$. Alltogether $\rk \Pic(X_{i_0}) = \rk \Pic(X) + 1$.
Similarly for any $X_I \subset X$ closed $\mathcal{G}$-stable subvariety, not contained in $X_{i_0}$,
$\rk\Pic(X_I\cap X_{i_0}) = \rk\Pic(X_I)+1$.

We will recall a classical result from algebraic geometry due to Steven L. Kleiman.
This will be used to invoke a general position argument for a certain curve to avoid 
the indeterminacy locus of $\restr{\varphi}{X_I}$, so that the limit map is 
defined everywhere on that curve.

\begin{proposition}[{\cite[Kleiman Bertini Theorem]{kleiman_tgt}}]\label{KBT}
  Let $G$ be an algebraic group. If $X$ is a homogeneous $G$-variety and $Y, Z \subset X$ subvarieties, 
  then there exists an open dense subset $G_0 \subset G$, such that for any $g \in G_0$
  the subvarieties $gY$ and $Z$ intersect properly and every component of $gY \cap Z$ has 
  dimension $\dim Y + \dim Z - \dim X$.
  In particular if $\dim Y + \dim Z < \dim X$, then $gY \cap Z = \emptyset$.
\end{proposition}

We will use the following definitions and notations for the rest of the section.

\begin{notation}\label{def_settings_for_irreducibility_proof}
  Let $X_I\subset X$, where as before $X_I = \cap_{i \in I} X_i$, $i_0 \notin I$ 
  be a closed $\mathcal{G}$-stable subvariety, not contained in $X_{i_0}$,
  write $I'$ for $I \cup \{i_0\}$ and $X_{I'}$ for the divisor $X_I\cap X_{i_0}$ of $X_I$.

\end{notation}

From now on let $\eta$ be a movable curve class on $X_I$.

\begin{lemma}\label{lemma_m_circ_non_empty}
  Any component of $M_{0,0}^\circ(X_I, \eta)$ contains 
  an irreducible representative of class $\eta$, that meets the open
  $\mathcal{G}$-orbit of $X_I$ and avoids the indeterminacy locus of $\restr{\varphi}{X_I}$.
  Furthermore, $M_{0,0}^\circ(X_I, \eta)$ is nonempty.
\end{lemma}

\begin{proof}
  By Corollary \ref{corollary_wond_g_stable_var_curve_rep} 
  there exists an irreducible curve $C$ of class $\eta$ passing through the open 
  $\mathcal{G}$-orbit of $X_I$, so $M_{0,0}^\circ(X_I, \eta) \neq \emptyset$.

  By definition every component of $M_{0,0}^\circ(X_I, \eta)$ must contain an 
  irreducible curve. Let $M \subset M_{0,0}^\circ(X_I, \eta)$ be some 
  component and let $C_M$ be such a curve.

  We have that $X_I$ is a union of homogeneous $\mathcal{G}$-varieties, namely the orbits of $\mathcal{G}$
  in $X_I$. The curve $C_M$ meets the open orbit and is irreducible, 
  so its intersection with every nonopen orbit is of dimension $0$.
  Since $\varphi$ is restricts to every orbit closure in $X_I$, we get that the 
  indeterminacy locus $I(\restr{\varphi}{X_I})$
  is of codimension at least $2$ for the open $\mathcal{G}$-orbit $X_I^\circ \subset X_I$
  by Proposition \ref{prop_limit_map_image_properties}.
  On the other hand $\restr{\varphi}{X_I}$ is defined on the open set $V$, which intersects
  every $\mathcal{G}$-orbit closure in $X_I$, so for any $\mathcal{G}$-orbit $O$
  the set $I(\restr{\varphi}{X_I})\cap O$ is of codimension at least $1$.
  This gives us that $\dim I(\restr{\varphi}{X_I})\cap O + \dim C_M \cap O < \dim O$.

  Next we apply Proposition \ref{KBT} to every $\mathcal{G}$ orbit $O$, therefore
  there exist open subsets $\mathcal{G}_O \subset \mathcal{G}$,
  such that for any $g \in \mathcal{G}_O$ we have $gC_M \cap I(\varphi)\cap O = \emptyset$.
  There is a finite number of orbits, let $\mathcal{G}^\circ = \cap_{O \subset X_I} \mathcal{G}_O$. 
  It is still open in $\mathcal{G}$ and for any $g \in \mathcal{G}^\circ$ the curve 
  $gC_M$ is a representative of class $\eta$ with required properties.
  $\qed$
\end{proof}

Let $C \in M_{0,0}^\circ(X_I, \eta)$ be an irreducible curve, avoiding the indeterminacy locus of $\varphi$.

\begin{lemma}\label{lemma_varphi_irreducible_passes_though_op_or}
  The image of $C$ under $\varphi$ is an irreducible curve $\varphi(C)$, 
  meeting the open $\mathcal{G}$-orbit of $X_{I'}$.
\end{lemma}

\begin{proof}
  The curve $C$ avoids the indeterminacy locus of $\varphi$, so $\varphi$ is defined everywhere on $C$ and
  as $C$ is irreducible, its image $\varphi(C)$ is irreducible as well.

  By Proposition \ref{prop_limit_map_image_properties} the open $\mathcal{G}$-orbit of $X_I$ is mapped via $\varphi$ 
  inside the open $\mathcal{G}$-orbit of $X_{I'}$. As $C$ meets the open $\mathcal{G}$-orbit of $X_I$,
  its image $\varphi(C)$ meets the open $\mathcal{G}$-orbit of $X_{I'}$.
  $\qed$
\end{proof}

\begin{lemma}\label{lemma_rat_equiv_limit_curve_groups}
  Let $\iota_{I'}^I: X_{I\cup \{i_0\}} \to X_I$ be the inclusion.
  Then $(\iota_{I'}^I)_*[\varphi(C)] = \eta$.
  In other words $\varphi(C)$ is rationally equivalent to $C$ in $X_I$.
\end{lemma}

\begin{proof}
  Let $V'$ be the open set, on which $\varphi$ is defined.
  For every $t \in \mathbb{C}^*$ the rational map
  $\varphi_t(x): x \mapsto \lambda(t) x$ is a true morphism.
  Consider the morphism
  \begin{align*}
    F: V' \times \mathbb{A}^1& \to X\\
      (x,t) &\mapsto \varphi_t(x),
  \end{align*}
  where $\varphi_0 = \varphi$. Let $p: V' \times \mathbb{A}^1$ denote the projection.
  Then $F(p^{-1}(1)) = C$ and $F(p^{-1}(0)) = \varphi(C)$, so
  $C$ and $\varphi(C)$ are members of a family of cycles,
  parametrized by a rational variety $\mathbb{A}^1$.
  By \cite[Example 10.1.7]{fulton_it} we have that $\varphi(C)$ is rationally equivalent to $C$ and
  $(\iota_{I'}^I)_*[\varphi(C)] = \eta$.
  $\qed$
\end{proof}


  

Let $\iota_I: X_I \to X$ be the inclusion of a closed $G$-stable subvariety
$X_I = \mathcal{G} \times^{\mathcal{P}_{S_I}^-}X^I$ as in Proposition \ref{prop_boundary_divisor_projection},
where $\mathcal{P}_{S_I}$ is the parabolic subgroup of $\mathcal{G}$ associated to the 
set of simple roots $S_I$. In case of wonderful group compactification
$S_I = \{\alpha_i, \beta_i\ |\ i \notin I\}$.
By Proposition \ref{prop_boundary_picard_exact_sequence}
$\Pic(X_I) = \Pic(X^I) \oplus \Pic(\mathcal{G}/{\mathcal{P}_{S_I}^-})$,
By Proposition \ref{prop_luna_limit_stuff} we have 
\[
  \Delta_{X^I} = \Delta_X(S_I),
\]
so 
\[
  \Pic(X^I) = \sum_{i \notin I} \mathbb{Z}D_i \cap X^I.
\]
For the homogeneous part we have 
\[
  \Pic(\mathcal{G}/{\mathcal{P}_{S_I}^-}) = \sum_{\alpha \in S \setminus S_I} \mathbb{Z}D_{\alpha}(\mathcal{G}/{\mathcal{P}_{S_I}^-}).
\]
Therefore in total we have 
\[
  \Pic(X_I) = \big(\sum_{i \notin I} \mathbb{Z}D_i \cap X_I \big)\oplus
    \big(\sum_{\alpha \in S \setminus S_I} \mathbb{Z}\pi_I^*D_{\alpha}(\mathcal{G}/{\mathcal{P}_{S_I}^-}) \big),
\]
where $\pi_I : X_I \to \mathcal{G}/{\mathcal{P}_{S_I}^-}$ is the $G$-equivariant projection.
In particular we have that 
\[
  \Delta_{X_I} = \{D_i \cap X_I\ |\ i \notin I\}\cup 
    \{\pi_I^*D_{\alpha_i}(\mathcal{G}/{\mathcal{P}_{S_I}^-}), \pi_I^*D_{\beta_i}(\mathcal{G}/{\mathcal{P}_{S_I}^-})\ |\ i \in I\}.
\]
We will write $D_i(X_I) = D_i \cap X_I$ for $i \notin I$ and $D_i^+(X_I) = \pi_I^*D_{\alpha_i}(\mathcal{G}/{\mathcal{P}_{S_I}^-})$,
as well as $D_i^-(X_I) = \pi_I^*D_{\beta_i}(\mathcal{G}/{\mathcal{P}_{S_I}^-})$ for $i \in I$.
Similar results hold for $X_{I'}$ and we want to understand the map $\Pic(X_I) \to \Pic(X_{I'})$ 
induced by the inclusion $\iota_{I'}^I : X_{I'} \to X_I$.

\begin{lemma}\label{lemma_colors_pullback}
  Let $\iota_{I'}^I : X_{I'} \to X_I$ be the inclusion of a boundary divisor
  and $I' = I \cup \{i_0\}$. Then
  \begin{itemize}
    \item $(\iota_{I'}^I)^* D_{i_0}(X_I) = D_{i_0}^+(X_I') + D_{i_0}^-(X_I')$;
    \item $(\iota_{I'}^I)^* D_i(X_I) = D_i(X_{I'})$, if $i \notin I'$;
    \item $(\iota_{I'}^I)^* D_i^+(X_I)  = D_{i}^+(X_I')$, if $i \in I$;
    \item $(\iota_{I'}^I)^* D_i^-(X_I)  = D_{i}^-(X_I')$, if $i \in I$;
  \end{itemize}
\end{lemma}

\begin{proof}
  For $i \notin I$ any color $D_i(X_I)$ is moved by exactly two simple roots $\alpha_i$ and $\beta_i$ of $G$,
  so the pullback $(\iota_{I'}^I)^* D_i(X_I)$ is moved by $\alpha_i$ and $\beta_i$ as well.
  
  If $i \notin I'$, then $(\iota_{I'}^I)^* D_i(X_I)$ is multiple of $D_i(X_{I'})$.
  
  If $i = i_0$, then $(\iota_{I'}^I)^* D_{i_0}(X_I)$ is a positive linear combination of $D_{i_0}^+(X_I')$ and $D_{i_0}^-(X_I')$,
  otherwise it would not be moved by both roots.

  For $i \in I$ any color $D_i^+(X_I)$ or $D_i^-(X_I)$ is moved by precisely one simple root of $G$,
  namely $\alpha_i$ or $\beta_i$ respectively. Therefore the pullback
  $(\iota_{I'}^I)^* D_i^+(X_I)$ or $(\iota_{I'}^I)^* D_i^-(X_I)$
  is a multiple of $D_{i}^+(X_I')$ or $D_{i}^-(X_I')$ respectively.
  
  Let $[C_D]$ be a curve class, dual to some color $D \in \Delta_{X_I}$
  as in Lemma \ref{lemma_wond_color_dual_curve_description}.
  Then we have that
  \[
    \langle D', [C_D] \rangle = \begin{cases}
      1 &\text{if }D' = D,\\
      0 &\text{else },
    \end{cases}
  \]
  for all $D' \in \Delta_{X_I}$, therefore the same also holds
  for all the curve classes $[C_{\alpha_i}]$ or $[C_{\beta_i}]$ of the closed orbit, as
  by Lemma \ref{lemma_wond_color_closed_orbit_curves} we have that $\iota_Y^I [C_{\alpha_i}] = [C_D]$ for
  some color $D$, moved by $\alpha_i$.
  So for any $D \in \Delta_{X_I}$ any coefficient is at most $1$ in the presentation of $(\iota_{I'}^I)^*D$
  as a linear combination of colors of $X_{I'}$.
  The result follows.
  $\qed$
\end{proof}

\begin{proposition}\label{prop_class_of_limit_curve_unique_groups}
  Let $\eta$ be a movable curve class.
  There exists a rational map 
  \begin{align*}
    \varphi_*: M^\circ_{0,0}(X_I, \eta) &\dashrightarrow M^\circ_{0,0}(X_{I\cup \{i_0\}}, \bar{\eta})\\
    C &\mapsto \varphi(C)
  \end{align*}
  where the class $\bar{\eta}$ is uniquely defined.
  Suppose that 
  \[
    \eta = \sum_{i \notin I} c_i [C_i(X_I)] + \sum_{i \in I} (a_i [C_{\alpha_i}(X_I)] + b_i [C_{\beta_i}(X_I)]),
  \]
  where the curve class $[C_i(X_I)]$ is dual to $D_i(X_I)$ and $[C_{\alpha_i}(X_I)]$ (resp. $[C_{\beta_i}(X_I)]$) 
  is dual to $D_i(X_I)^+$ (resp. $D_i(X_I)^-$), then with similar notation for $X_{I'}$
  \[
    \bar{\eta} = \sum_{i \notin I'} c_i [C_i(X_{I'})] + c_{i_0} [C_{\beta_{i_0}}(X_{I'})] +
      \sum_{i \in I} (a_i [C_{\alpha_i}(X_{I'})] + b_i [C_{\beta_i}(X_{I'})]).
  \]
\end{proposition}

\begin{proof}
  Proposition \ref{prop_limit_map_image_properties} tells us that $\varphi(C) \subset X_I \cap X_{i_0} = X_{I'}$
  and Lemma \ref{lemma_varphi_irreducible_passes_though_op_or} tell us that it lies in
  $M^\circ_{0,0}(X_{I'}, \bar{\eta})$ for some $\bar{\eta}$.
  We calculate the class $\bar{\eta}$ independently of the representative $C$.
  
  Write $I' = I \cup \{i_0\}$ for brevity, then we have 
  \[
    \Delta_{X_I} = \{D_i(X_I)\ |\ i \notin I\} \cup \{D_i(X_I)^+, D_i(X_I)^-\ |\ i \in I\}
  \]
  and 
  \[
    \Delta_{X_{I'}} = \{D_i(X_{I'})\ |\ i \notin I'\} \cup \{D_i(X_{I'})^+, D_i(X_{I'})^-\ |\ i \in I'\} 
  \]
  By Lemma \ref{lemma_colors_pullback} 
  we have $(\iota_{I'}^I)^*D_{i_0}(X_I) = D_{i_0}(X_{I'})^+ + D_{i_0}(X_{I'})^-$.
  For $i \neq i_0$ we have $(\iota_{I'}^I)^*D_i(X_I) = D_i(X_{I'})$ 
  (resp. $(\iota_{I'}^I)^*D_i^+(X_I) = D_i^+(X_{I'})$ and $(\iota_{I'}^I)^*D_i^-(X_I) = D_i^-(X_{I'})$).
  By Lemma \ref{lemma_wond_color_dual_curve_description}
  \[
    \bar{\eta} = \sum_{i \notin I'} \bar{c}_i [C_i(X_{I'})] + \sum_{i \in I'} (\bar{a}_i [C_{\alpha_i}(X_{I'})] + \bar{b}_i [C_{\beta_i}(X_{I'})]),
  \]
  where the curve class $[C_i(X_{I'})]$ is dual to $D_i(X_{I'})$ and $[C_{\alpha_i}(X_{I'})]$ (resp. $[C_{\beta_i}(X_{I'})]$) 
  is dual to $D_i(X_{I'})^+$ (resp. $D_i(X_{I'})^-$).
  Similarly we write
  \[
    \eta = \sum_{i \notin I} c_i [C_i(X_I)] + \sum_{i \in I} (a_i [C_{\alpha_i}(X_I)] + b_i [C_{\beta_i}(X_I)]).
  \]

  Now we can apply the projection formula 
  to the inclusion $(\iota_{I'}^I) : X_{I'} \to X_I$:
  \[
    \langle (\iota_{I'}^I)^*D, \bar{\eta} \rangle = \langle D, \iota_*\bar{\eta} \rangle = \langle D, \eta \rangle,
  \]
  for any color $D \in \Delta_{X_I}$, where the last equality follows from Lemma \ref{lemma_rat_equiv_limit_curve_groups}.
  This yields the following values for $\langle (\iota_{I'}^I)^*D, \bar{\eta} \rangle$:
  \[
    \langle (\iota_{I'}^I)^*D, \bar{\eta} \rangle = \langle D, \eta \rangle = \begin{cases}
      c_i &D = D_i(X_I), i \notin I, \\
      a_i &D = D_i^+(X_I), i \in I, \\
      b_i &D = D_i^-(X_I), i \in I. \\
    \end{cases}
  \]
  On the other hand we have that $(\iota_{I'}^I)^*D \in \Delta_{X_{I'}}$ for $D \neq D_{i_0}(X_I)$, 
  so 
  \[
    \langle (\iota_{I'}^I)^*D, \bar{\eta} \rangle = \begin{cases}
      \bar{c}_i &D = D_i(X_I), i \notin I \cup \{i_0\}, \\
      \bar{a}_i &D = D_i^+(X_I), i \in I, \\
      \bar{b}_i &D = D_i^-(X_I), i \in I. \\
    \end{cases}
  \]
  Therefore for any $i \neq i_0$ we have $\bar{a}_i = a_i$, $\bar{b}_i = b_i$ and $\bar{c}_i = c_i$.

  Lastly we wish to compute $\bar{a}_{i_0}$ and $\bar{b}_{i_0}$.
  We already know that 
  \[
    \langle D_{i_0}^+(X_{I'}) + D_{i_0}^-(X_{I'}), \bar{\eta} \rangle = \bar{a}_{i_0} + \bar{b}_{i_0} = 
      c_{i_0} = \langle D_{i_0}(X_I), \eta \rangle  
  \]
  
  We wish to compute either $\bar{a}_{i_0}$ or $\bar{b}_{i_0}$ to get the final result.
  Let us look at the following commutative diagram:
  \[
    \begin{tikzcd}
      X_I \arrow[d, "p_I"] \arrow[r, dashed, "\varphi"] & X_{I'} \arrow[d, "p_{I'}"] \\
      \mathcal{G}/\mathcal{P}_I \arrow[rd, dashed] \arrow[r, dashed, "\tilde{\varphi}"] 
        & \mathcal{G}/\mathcal{P}_{I'} \arrow[d, "\theta"] \\
      & \mathcal{G}/\mathcal{P}_{i_0}
    \end{tikzcd}
  \]
  Here $\mathcal{P}_I$, $\mathcal{P}_{I'}$ and $\mathcal{P}_{i_0}$ are
  parabolic subgroups of $\mathcal{G}$ associated to $\{\alpha_i, \beta_i \ |\ i \notin I\} \subset S$,
  $\{\alpha_i, \beta_i \in S\ |\ i \notin I'\}$ and $\{\alpha_i, \beta_i \ |\ i \neq i_0\}$ respectively.
  Note that $\mathcal{P}_{i_0} = P_{i_0} \times P_{i_0}^-$, where $P_{i_0}$ is the maximal parabolic, 
  associated to $\{\alpha_i \ |\ i \neq i_0\}$.
  The map $\theta:\mathcal{G}/\mathcal{P}_{I'} \to \mathcal{G}/\mathcal{P}_{i_0}$ is the usual map 
  of flag varieties, since $\mathcal{P}_{I'} \subset \mathcal{P}_I \subset \mathcal{P}_{i_0}$.
  So we can write 
  \[
    \mathcal{G}/\mathcal{P}_{i_0} = G/P_{i_0} \times G/P_{i_0}^-
  \]
  On $\mathcal{G}/\mathcal{P}_{i_0}$ we have two divisors, $D_{\alpha_{i_0}}(\mathcal{G}/\mathcal{P}_{i_0})$
  and $D_{\beta_{i_0}}(\mathcal{G}/\mathcal{P}_{i_0})$. We have that 
  $D_{i_0}^+(X_{I'}) = \theta^*D_{\alpha_{i_0}}(\mathcal{G}/\mathcal{P}_{i_0})$
  and $D_{i_0}^-(X_{I'}) = \theta^*D_{\beta_{i_0}}(\mathcal{G}/\mathcal{P}_{i_0})$.
  
  The limit map $\varphi$ is defined on an open affine subset of $\mathcal{G}/\mathcal{P}_{i_0}$.
  This open affine subset is isomorphic to
  \[
    \prod_{\alpha \geq \alpha_{i_0}} U_\alpha \times \prod_{\beta \geq \beta_{i_0}} U_\beta,
  \]
  where $\alpha$ is a linear combination of roots $\alpha_1, \dots, \alpha_l$ and 
  $\beta$ is a  linear combination of roots $\beta_1, \dots, \beta_l$.
  As the weight $\lambda$, defining the limit map $\varphi$, is trivial on the roots $\beta_{i_0}$, 
  the limit map $\varphi$ itself is constant on the part $\prod_{\beta \geq \beta_{i_0}} U_\beta$
  and therefore it is trivial on the second factor of $\mathcal{G}/\mathcal{P}_{i_0}$.
  On the first factor, since $\langle \lambda, \alpha \rangle > 0$ for 
  any root $\alpha \geq \alpha_{i_0}$ of $G$, we get that 
  \[
    \varphi (U_\alpha) = \lim_{t \to 0} \lambda(t) U_\alpha
      = \lim_{t \to 0} t^{\langle \lambda, \alpha \rangle}U_\alpha = 0.
  \]
  In other words, $\varphi$ contracts an open subset to a point. 

  Then we have 
  \[
    \langle D_{i_0}^+(X_{I'}), \bar{\eta} \rangle
      = \langle \theta^*D_{\alpha_{i_0}}(\mathcal{G}/\mathcal{P}_{i_0}), \bar{\eta} \rangle.
  \]
  We can again apply the projection formula to $\theta$:
  \[
    \langle \theta^*D_{\alpha_{i_0}}(\mathcal{G}/\mathcal{P}_{i_0}), \bar{\eta} \rangle
      = \langle D_{\alpha_{i_0}}(\mathcal{G}/\mathcal{P}_{i_0}), \theta_*\bar{\eta} \rangle.
  \]

  By Lemma \ref{lemma_m_circ_non_empty} there exists an open dense subset $M \subset M_{0,0}^\circ(X_I, \eta)$ 
  of curves, that avoid the indeterminacy locus of $\varphi$. This allows us to pick 
  a representative $C \in M$ and to look at its image directly.
  The class $\theta_*\bar{\eta}$ is represented by the curve $\theta(\varphi(C))$. 
  As $\varphi$ is defined via the $\mathcal{G}$-action, it commmutes with $\theta$,
  which is $\mathcal{G}$-equivariant. Then we get 
  \[
    \langle D_{i_0}^+(X_{I'}), \bar{\eta} \rangle
      = \langle D_{\alpha_{i_0}}(\mathcal{G}/\mathcal{P}_{i_0}), [\theta(\varphi(C))] \rangle
        = \langle D_{\alpha_{i_0}}(\mathcal{G}/\mathcal{P}_{i_0}), [\varphi(\theta(C))] \rangle.
  \]
  note that 
  \[
    \langle D_{\alpha_{i_0}}(\mathcal{G}/\mathcal{P}_{i_0}), [\varphi(\theta(C))] \rangle = 0,
  \]
  because $\varphi(\theta(C))$ is of dimension $0$ on 
  the first factor of $\mathcal{G}/\mathcal{P}_{i_0}$.

  Altogether we have that 
  \[
    \bar{a}_{i_0} = \langle D_{i_0}^+(X_{I'}), \bar{\eta} \rangle = 0,
  \] 
  which means that 
  \[
    \bar{b}_{i_0} = \langle D_{i_0}^-(X_{I'}), \bar{\eta} \rangle = \langle D_{i_0}, \eta \rangle.
  \]

  So in the end we have that 
  \[
    \bar{\eta} = \sum_{i \notin I'} c_i [C_i(X_{I'})] + c_{i_0} [C_{\beta_{i_0}}(X_{I'})] +
      \sum_{i \in I} (a_i [C_{\alpha_i}(X_{I'})] + b_i [C_{\beta_i}(X_{I'})]),
  \]
  proving the result.
  $\qed$
\end{proof}

\begin{proposition}\label{prop_curve_point_smooth_boundary_div}
  The point $\varphi(C)$ is smooth in $\overline{M}_{0,0}(X_{I'},\bar{\eta})$ and the 
  point $(\iota_{I'}^I)\circ\varphi(C)$ is smooth in $\overline{M}_{0,0}(X_I, \eta)$.
\end{proposition}

\begin{proof}
  Let $(\iota_{I'}^I): X_{I'} \to X_I$ be the inclusion, 
  $f: \mathbb{P}^1 \to X_I$ be the curve $C$ and $\bar{f} = \varphi\circ f$ be
  the curve $\varphi(C)$ on $X_{I'}$.
  Then $(\iota_{I'}^I) \circ \bar{f}: \mathbb{P}^1 \to X_I$ is a curve on $X_I$.
  As $\varphi(C)$ meets the open $\mathcal{G}$-orbit of $X_{I'}$, 
  we have that $H^1(\bar{f}^*T_{X_{I'}}) = 0$ and the point $\varphi(C)$ is smooth 
  in $M_{0,0}(X_{I'}, \bar{\eta})$ by Proposition \ref{prop_moduli_space_dimension}.
  
  Now we want to show that $H^1(((\iota_{I'}^I) \circ \bar{f})^*T_{X_I}) = 0$.
  Let us consider the following exact sequence:
  \[
  \begin{tikzcd}
    0 \arrow[r] & T_{X_{I'}} \arrow[r] & \restr{T_{X_I}}{X_{I'}} \arrow[r] & O_{X_{I'}}(X_{I'})
    \arrow[r] & 0
  \end{tikzcd}
  \]
  and pull it back onto $\mathbb{P}^1$ to get:
  \[
  \begin{tikzcd}
    0 \arrow[r] & \bar{f}^*T_{X_{I'}} \arrow[r] & \bar{f}^*\restr{T_{X_I}}{X_{I'}} \arrow[r] & \bar{f}^*O_{X_{I'}}(X_{I'})
    \arrow[r] & 0
  \end{tikzcd}
  \]
  We have $\restr{T_{X_I}}{X_I'} = (\iota_{I'}^I)^* T_{X_I}$. 
  Replace $\bar{f}^*\restr{T_X}{X_{I'}}$ by $(\iota_{I'}^I \circ \bar{f})^*T_{X_I}$:
  \[
  \begin{tikzcd}
    0 \arrow[r] & \bar{f}^*T_{X_{I'}} \arrow[r] & (\iota_{I'}^I \circ \bar{f})^*T_{X_I} \arrow[r] & \bar{f}^*O_{X_{I'}}(X_{I'})
    \arrow[r] & 0
  \end{tikzcd}
  \]
  Here we notice that $\bar{f}^*O_{X_{I'}}(X_{I'}) = O_{\mathbb{P}^1}(\langle X_{I'}, \eta \rangle)$, 
  so when we take $H^1$ of the previous exact sequence we get:
  \[
  \begin{tikzcd}
    H^1(\bar{f}^*T_{X_{I'}}) \arrow[r] & H^1((\iota_{I'}^I \circ \bar{f})^*T_{X_I}) \arrow[r] & H^1(O_{\mathbb{P}^1}(\langle X_{I'}, \eta \rangle)) \arrow[r] & 0
  \end{tikzcd}
  \]
  In this exact sequence $H^1(O_{\mathbb{P}^1}(\langle X_{I'}, \eta \rangle))$ vanishes because
  $\langle X_{I'}, \eta \rangle \geq 0 \geq -1$ by definition of $\eta$ 
  and $H^1(\bar{f}^*T_{X_{I'}})$ vanishes since the curve 
  $f^*C$ meets the open $\mathcal{G}$-orbit of $X_{I'}$.
  With that $H^1((\iota_{I'}^I \circ \bar{f})^*T_{X_I}) = 0$ and the point $\varphi(C)$ is smooth in $\overline{M}_{0,0}(X_I, \eta)$
  by Proposition \ref{prop_moduli_space_dimension}.
  $\qed$
\end{proof}

\begin{lemma}\label{lemma_irreducibility_induction_step}
  Let $X_{I'} = X_I \cap X_{i_0} \subset X_I$ be a boundary divisor of a $G$-stable subvariety $X_I \subset X$.
  Suppose that the space $M_{0,0}^\circ(X_{I'}, \bar{\eta})$ is irreducible,
  then the space $M_{0,0}^\circ({X_I, \eta})$ is irreducible as well.
\end{lemma}

\begin{proof}
  Write $M$ for the open subset of $M_{0,0}^\circ(X,\eta)$, on which $\varphi_*$ is defined.
  By Lemma \ref{lemma_m_circ_non_empty} any component of $M_{0,0}^\circ(X_I, \eta)$
  meets $M$ and by Proposition \ref{prop_moduli_space_dimension}
  $M_{0,0}^\circ(X_I, \eta)$ only contain smooth points, so $M$ is dense in $M_{0,0}^\circ(X_I, \eta)$.
  By Lemma \ref{lemma_rat_equiv_limit_curve_groups}
  $(\iota_{I'}^I)_*\bar{\eta} = \eta$, so
  $\varphi_*(M) \subset M_{0,0}^\circ(X_{I'}, \bar{\eta})$.
  We can also apply the limit map $\varphi_*$ to the set $M_{0,0}^\circ(X_{I'}, \bar{\eta})$,
  as it is endowed with a $G$-action as well, so we get a following diagram:
  \[
    \begin{tikzcd}
      M_{0,0}^\circ(X_{I'}, \bar{\eta}) \arrow[r, "\iota_*"] \arrow[rd, "\varphi_*", dashed] 
        & \overline{M} \arrow[d, "\varphi_*", dashed] \\
      & M_{0,0}^\circ(X_{I'}, \bar{\eta}) 
    \end{tikzcd}
  \]

  For a curve $C \subset X_I$, such that $\varphi$ is defined on $C$,
  the curve $\varphi(C)$ is stable under the action of $\lambda(\mathbb{C})$, 
  so $\varphi(\varphi(C)) = \varphi(C)$.
  Let $M'$ be the open subset of $M_{0,0}^\circ(X_{I'}, \bar{\eta})$, where
  $\varphi_*$ is defined. Then $\varphi_*(M') = \varphi_*(M)$.
  As $M_{0,0}^\circ(X_{I'}, \bar{\eta})$ is irreducible,
  the image $\varphi_*(M')$ is irreducible as well and so is $\varphi_*(M)$.

  Next, we look at the map $\varphi_*: M \to \varphi_*(M)$.
  Its image is irreducible, so there exists an irreducible component $M_0 \subset M$
  such that $\varphi_*(M_0) = \varphi_*(M)$.  
  Suppose there exists another irreducible component $M_1 \subset M$.
  Since $(\varphi_t(C))_{t\in \mathbb{C}}$ is a continuous family of curves, we get that 
  $\varphi_*(M_0) \subset \overline{M_0}$ and $\varphi_*(M_1) \subset \overline{M_1}$.
  Also we have $\varphi_*(M_1) \subset \varphi_*(M) = \varphi_*(M_0) \subset \overline{M_0}$.
  
  Then any point $x \in \varphi_*(M_1)$ lies in the closures $\overline{M_0}$ and $\overline{M_1}$
  of two irreducible components in $\overline{M}$. 
  As $x$ lies in the image $\varphi_*(M_1) \subset \varphi_*(M)$, 
  it is of the form $\varphi(C_0)$ for some $C_0 \in M \subset M_{0,0}^\circ(X_I,\eta)$. 
  By Proposition \ref{prop_curve_point_smooth_boundary_div} the point $x = \varphi(C_0)$ is smooth in both 
  $M_{0,0}(X_I, \eta)$ and $M_{0,0}(X_{I'}, \eta)$, then the two irreducible components meet 
  at the smooth point $x$ and therefore must coincide.
  As $M$ is dense in $M_{0,0}^\circ(X_I,\eta)$, the result follows.
  $\qed$
\end{proof}

\begin{theorem}\label{theorem_irreducibility_group_case}
  Let $X$ be the wonderful compactification of a semisimple adjoint group $G$. Let $\eta$ be a movable curve class on $X$.
  Then the space $\overline{M_{0,0}^\circ(X,\eta)}$ is irreducible.
\end{theorem}

\begin{proof}
  First we note that by Proposition \ref{corollary_wond_g_stable_var_curve_rep}
  $M_{0,0}^\circ(X, \eta)$ is nonempty.

  We have a chain of $G$-stable subvarieties in $X$:
  \[
    \cap_{i = 1}^l X_i \subset \cap_{i = 2}^l X_i \subset \dots \subset X_l \subset X,
  \]
  which can also be written as 
  \[
    X_{\{1,\dots, l\}} \subset X_{\{2,\dots,l\}} \subset \dots \subset X_{\{l-1,l\}} \subset X_{\{l\}} \subset X
  \]
  where every pair $X_{I'} = X_{I\cup \{j\}} \subset X_{I}$ with $I = \{j + 1,\dots, l\}$ 
  satisfies the conditions of Lemma \ref{lemma_irreducibility_induction_step} for any $j \in \{1, \dots, l\}$.
  
  By Proposition \ref{prop_homogeneous_variety_mod_space_irreducible} the moduli space 
  $M_{0,0}^\circ(\cap_{i = 1}^l X_i, \eta)$ is irreducible and by Lemma \ref{lemma_irreducibility_induction_step}
  the same holds for $M_{0,0}^\circ(\cap_{i = 2}^l X_i, \eta)$. 
  By trivial induction on $j \in \{1, \dots, l\}$ we get that $M_{0,0}^\circ (\cap_{i=j+1}^l X_i, \eta)$ 
  is irreducible for every $j \ni \{1, \dots, l\}$, in particular 
  for $j = l$, which yields the irreducibility of $M_{0,0}^\circ(X, \eta)$.
  $\qed$
\end{proof}

\bibliographystyle{alpha}
\bibliography{these}

\end{document}